\documentclass{article}
\usepackage{graphicx} 
\usepackage{fullpage}
\usepackage{amsfonts}
\usepackage{amsmath}
\usepackage{amsthm}
\usepackage{amssymb}
\usepackage{color}
\usepackage{xcolor}
\usepackage[colorlinks,linkcolor=blue, citecolor=blue]{hyperref}
\usepackage[section]{algorithm}
\usepackage{algpseudocode}
\usepackage{framed}
\usepackage{footnote}
\usepackage{geometry}
\usepackage{xfrac}
\usepackage{todonotes}
\usepackage{tabularx,booktabs}
\usepackage{svg,cite}

\usepackage{psfrag}
\usepackage{pdflscape}
\usepackage{multirow}
\usepackage{blindtext} 
\usepackage{enumitem}
\usepackage{subcaption}
\usepackage[capitalise]{cleveref}
\crefname{subsection}{Subsection}{Subsections}

\usepackage{asb}

\numberwithin{equation}{section}

\newtheorem{theorem}{Theorem}[section]

\newtheorem{lemma}[theorem]{Lemma}

\newtheorem{assumption}{Assumption}[section]

\newtheorem{condition}{Condition}[section]
\newtheorem{corollary}[theorem]{Corollary}

\numberwithin{algorithm}{section}

\usepackage{array}
\usepackage{arydshln}
\usepackage{stackengine}

\newcommand{\papertitle}{On the Convergence and Complexity of Proximal Gradient and Accelerated Proximal Gradient Methods under Adaptive Gradient Estimation}

\newcommand{\paperauthorb}{Raghu Bollapragada}
\newcommand{\paperauthorbaffiliation}{Operations Research and Industrial Engineering Program, University of Texas at Austin}
\newcommand{\paperauthorbemail}{raghu.bollapragada@utexas.edu}
\newcommand{\paperauthorc}{Shagun Gupta}

\newcommand{\paperauthorcemail}{shagungupta@utexas.edu}

\begin{document}
\title{\papertitle}

\author{
    \paperauthorb\footnotemark[2]\ 
    \footnotemark[3]
    \and \paperauthorc\footnotemark[2]
}

\maketitle

\renewcommand{\thefootnote}{\fnsymbol{footnote}}
\footnotetext[2]{\paperauthorbaffiliation. (\url{\paperauthorbemail},\url{\paperauthorcemail})}
\footnotetext[3]{Corresponding author.}
\renewcommand{\thefootnote}{\arabic{footnote}}

\begin{abstract}
   In this paper, we propose a proximal gradient method and an accelerated proximal gradient method for solving composite optimization problems, where the objective function is the sum of a smooth and a convex, possibly nonsmooth, function.
We consider settings where the smooth component is either a finite-sum function or an expectation of a stochastic function, making it computationally expensive or impractical to evaluate its gradient.
To address this, we utilize gradient estimates within the proximal gradient framework. Our methods dynamically adjust the accuracy of these estimates, increasing it as the iterates approach a solution, thereby enabling high-precision solutions with minimal computational cost.
We analyze the methods when the smooth component is nonconvex, convex, or strongly convex, using a biased gradient estimate. In all cases, the methods achieve the optimal iteration complexity for first-order methods.
When the gradient estimate is unbiased, we further refine the analysis to show that the methods simultaneously achieve optimal iteration complexity and optimal complexity in terms of the number of stochastic gradient evaluations. Finally, we validate our theoretical results through numerical experiments.
\end{abstract}
\section{Introduction}

In this paper, we consider composite optimization problems of the form
\begin{equation} \label{eq:intro_composite_problem}
    \min_{x \in \Rmbb^d} \; \phi(x) = f(x) + h(x),
\end{equation}
where $f : \Rmbb^d \rightarrow \Rmbb$ is a continuously differentiable function and $h : \Rmbb^d \rightarrow \Rmbb \cup \{+\infty\}$ is a closed, convex, proper, and possibly nonsmooth function. We focus on settings where $h(x)$ admits a simple structure that enables efficient computation of the proximal operator,
\begin{equation} \label{eq:prox_operator}
    \prox_{\alpha, h} (y) = \argmin_{x \in \Rmbb^d} \; h(x) + \frac{1}{2\alpha} \|x - y\|^2, \quad \mbox{with } \alpha>0,
\end{equation}
which allows for the use of proximal gradient methods \cite{bertsekas2015convex,duchi2010composite} to efficiently solve \eqref{eq:intro_composite_problem}.
Examples of such $h(x)$ include the $l_1$-norm penalty and the 
indicator function of a set to enforce 
simple convex constraints such as box constraints, norm-ball constraints, or boundary conditions.
This class of problems has been extensively studied due to its wide range of applications, including image processing \cite{fadili2010total,hansen2006deblurring}, data science \cite{chen2010graph,jenatton2010proximal}, and inverse problems \cite{bauschke2011fixed,tropp2010computational}.  
When $f(x)$ is a convex function, Nesterov's acceleration \cite{nesterov1983method} can be applied to proximal gradient methods to achieve improved convergence rates, as shown in \cite{beck2009fast,scheinberg2014fast,schmidt2011convergence,nesterov2013gradient}.

We analyze problems where the smooth component $f(x)$ takes one of the following forms:
\newcounter{eqnAP}
\newcounter{eqnBP}
\refstepcounter{equation}
\setcounter{eqnAP}{\value{equation}}
\label{eq:deter_error_obj}
\refstepcounter{equation}
\setcounter{eqnBP}{\value{equation}}
\label{eq:stoch_error_obj}
\begin{align*}
    f(x) = \dfrac{1}{N} \sum_{i =1}^N F(x, \xi_i)  \qquad \qquad \text{(\ref*{eq:deter_error_obj}),} \qquad \text{or} \qquad  \qquad
    f(x) = \Embb [F(x, \xi)], \qquad \qquad \text{(\ref*{eq:stoch_error_obj})}
\end{align*}
where \eqref{eq:deter_error_obj} is a finite-sum function resulting in a finite-sum problem over the dataset $\Scal = \{\xi_1, \xi_2, \ldots, \xi_N\}$ with $F : \Rmbb^d \times \Scal \rightarrow \Rmbb$, and \eqref{eq:stoch_error_obj} is an expectation function resulting in an expectation problem over the random variable $\xi$ with the associated probability space $(\Xi,\Omega, \Pcal)$, $F : \Rmbb^d \times \Xi \rightarrow \Rmbb$ and $\Embb[\cdot]$ denotes the expectation with respect to $\Pcal$.
Since computing the exact gradient of $f(x)$ in these settings is computationally prohibitive, proximal gradient methods rely on gradient estimates \cite{beiser2023adaptive,xie2024constrained,nguyen2024stochastic,ghadimi2016accelerated,kulunchakov2019generic,hu2009accelerated,pham2020proxsarah,lan2012optimal}.
The main computational costs in such methods are: (1) the number of proximal operator evaluations (i.e., iterations) and, (2) the number of stochastic gradient evaluations.
Most existing methods use unbiased gradient estimates of fixed accuracy \cite{ghadimi2016accelerated,hu2009accelerated,pham2020proxsarah,lan2012optimal}, typically achieving optimal performance with respect to only one of the two costs. Adaptive gradient estimation methods, by contrast, dynamically adjust the accuracy of the gradient estimate based on the quality of the current iterate. By using low-accuracy estimates in the early stages and gradually increasing the accuracy of the estimates, they can achieve high-accuracy solutions with minimal computational effort. To this end, we propose proximal gradient methods with adaptive gradient estimation to optimize both cost metrics.
Furthermore, in many practical scenarios, the bias in the gradient estimate is intrinsic and cannot be fully eliminated, such as federated learning with non-IID data \cite{li2019convergence} and Bayesian optimization using surrogate models \cite{snoek2012practical}.

As a result, we propose and analyze a proximal gradient method and an accelerated proximal gradient method for the finite-sum problem \eqref{eq:deter_error_obj} and the expectation problem \eqref{eq:stoch_error_obj} that adaptively control the accuracy of the estimate and allow for biased gradient estimates.
These methods achieve the optimal iteration complexity for first-order methods when the objective function is nonconvex, convex and strongly convex while using biased gradient estimates. When the gradient estimate is unbiased, the methods additionally attain the optimal complexity for the number of stochastic gradient evaluations in all three settings. A summary of these complexity results is provided in \cref{tab:Complexity_summary}.

\subsection{Literature Review}

Proximal gradient methods for deterministic composite optimization problems are well studied; see \cite{bertsekas2015convex,duchi2010composite} and references therein. Accelerated variants for deterministic convex problems, such as Nesterov’s acceleration \cite{nesterov2013gradient,schmidt2011convergence}, the fast iterative
shrinkage-thresholding algorithm (FISTA) \cite{beck2009fast}, and its backtracking extension \cite{scheinberg2014fast}, are also well established in the literature. Several works \cite{sun2019convergence,schmidt2011convergence,rebegoldi2022scaled,bonettini2018inertial} have analyzed proximal and accelerated proximal gradient methods with inexact gradients, where gradient accuracy is controlled through a predetermined deterministic sequence, in settings similar to the finite-sum problem \eqref{eq:deter_error_obj}. In contrast, our work employs Nesterov's acceleration and adaptively adjusts the accuracy of the gradient estimates based on the current iterate for both the finite-sum \eqref{eq:deter_error_obj} and expectation \eqref{eq:stoch_error_obj} problems.

Many proximal and accelerated proximal gradient methods have been proposed that employ unbiased stochastic gradient estimates for the finite-sum problem \eqref{eq:deter_error_obj} and the expectation problem \eqref{eq:stoch_error_obj}.
In \cite{kulunchakov2019generic,hu2009accelerated,lan2012optimal}, the authors proposed accelerated methods that achieve the optimal complexity in terms of the number of stochastic gradient evaluations for the expectation problem \eqref{eq:stoch_error_obj} when the objective function is convex or strongly convex. 
In \cite{ghadimi2016accelerated}, accelerated proximal gradient methods using Nesterov's acceleration were analyzed, establishing optimal complexity results for the number of proximal operator evaluations and the number of  stochastic gradient evaluations for the expectation problem when the objective function is nonconvex.
In \cite{pham2020proxsarah,j2016proximal,li2018simple}, the authors employ accelerated proximal gradient methods with variance reduction techniques and achieve the optimal complexity in terms of stochastic gradient evaluations (in expectation) for the finite-sum problem \eqref{eq:deter_error_obj}. In \cite{nguyen2024stochastic}, the authors extend the FISTA based method \cite{scheinberg2014fast} to the expectation problem \eqref{eq:stoch_error_obj} while allowing for biased gradient estimates. 
The method replaced the backtracking line search in \cite{scheinberg2014fast} with a step search mechanism to adaptively determine the step size, and assumed control over the error in the gradient estimate in probability, unlike the other works discussed above, which control the expected error in the gradient estimate. Although this method establishes stronger convergence guarantees than convergence in expectation under its assumptions, the analysis is limited to general convex objectives and results in suboptimal complexity for the number of stochastic gradient evaluations.

In \cite{beiser2023adaptive,xie2024constrained}, the authors employed unbiased gradient estimates in proximal gradient methods and used adaptive sampling strategies to control the accuracy of the gradient estimate for the expectation problem \eqref{eq:stoch_error_obj}. They established theoretical convergence guarantees for various objective functions but did not provide complexity results for the number of stochastic gradient evaluations.
We adopt similar conditions to \cite{xie2024constrained} to control the accuracy of the gradient estimate, while allowing for biased estimates and extending the approach to accelerated proximal gradient methods, achieving optimal iteration complexity for first-order methods \cite{nesterov2018lectures,carmon2020lower}.
When unbiased gradient estimates are available via sample average approximations, our conditions guide sample size selection in a way similar to \cite{beiser2023adaptive,xie2024constrained}, and also yield the optimal complexity for the number of stochastic gradient evaluations for the expectation problem \eqref{eq:stoch_error_obj} for first-order methods \cite{ghadimi2012optimal,lan2012optimal,arjevani2023lower,ghadimi2016accelerated}.  
These results are summarized in \cref{tab:Complexity_summary}. Similar complexity guarantees for the finite-sum problem \eqref{eq:deter_error_obj} can also be established following the same procedure and are omitted for brevity. 

Lastly, while our methods assume exact solutions to the proximal operator, several works have explored algorithms with inexact proximal updates \cite{bonettini2018inertial,schmidt2011convergence}. Additionally, approaches that go beyond the structured nonsmoothness in problem \eqref{eq:intro_composite_problem} have been investigated in \cite{wang2023SQP} for constrained settings, and in \cite{jalilzadeh2022smoothed,zhang2024private,marrinan2023zeroth} where smoothing techniques are employed to develop zeroth-order methods. 

\begin{table}
\caption{Summary of the best-established complexity results in this paper.} 
\label{tab:Complexity_summary}
\begin{center}
\begin{tabular}{lc*{4}{>{\centering\arraybackslash}p{2.5cm}}}
    \toprule
    \multirow{3}{*}{$f(x)$} &\multicolumn{2}{c}{Proximal Gradient} &\multicolumn{2}{c}{Accelerated Proximal Gradient} \\ \cmidrule{2-5}
    &Proximal&Gradients&Proximal&Gradients\\
    &Operators&(Unbiased)&Operators&(Unbiased)\\\midrule
    \addstackgap[15pt]
    Nonconvex & $\Ocal\left(\dfrac{1}{\epsilon}\right)$ & $\Ocal\left(\dfrac{1}{\epsilon^2}\right)$ & - & - \\
    \hdashline \addstackgap[15pt]
    Convex & $\Ocal\left(\dfrac{1}{\epsilon}\right)$ & $\Ocal\left(\dfrac{1}{\epsilon^2}\right)$ & $\Ocal\left(\dfrac{1}{\sqrt{\epsilon}}\right)$ & $\Ocal\left(\dfrac{1}{\epsilon^{1.5}}\right)$ \\
    \hdashline \addstackgap[15pt]
    Strongly Convex & $\Ocal\left(\kappa \log \dfrac{1}{\epsilon}\right)$ & $\Ocal\left(\dfrac{\kappa}{\epsilon}\right)$ & $\Ocal\left(\sqrt{\kappa} \log \dfrac{1}{\epsilon}\right)$ & $\Ocal\left( \dfrac{\sqrt{\kappa}}{\epsilon}\right)$ \\
    \bottomrule
\end{tabular}
\end{center}

Note: The solution accuracy $\epsilon$ depends on the nature of the objective function: for nonconvex functions, it refers to the squared norm of the gradient; for convex and strongly convex functions, it corresponds to the optimality gap in function value. Here, $\kappa$ denotes the condition number.
\end{table}












    

\subsection{Contributions}
We summarize our main contributions as follows:
\begin{enumerate}
    \item We propose a proximal gradient method and an accelerated proximal gradient method for both the finite-sum problem \eqref{eq:deter_error_obj} and the expectation problem \eqref{eq:stoch_error_obj}. These methods employ gradient estimates whose accuracy is adaptively controlled using deterministic and stochastic generalizations of the well-known “norm condition” \cite{bollapragada2018adaptive}.

    \item We show that the proposed methods achieve the optimal complexity in terms of the number of proximal operator evaluations for first-order methods, even when using biased gradient estimates, for nonconvex, convex, and strongly convex objective functions in both the finite-sum and expectation problems; see \cref{tab:Complexity_summary}.
    
    \item We further show that, when the gradient estimates are unbiased, the methods simultaneously achieve the optimal complexity for the number of stochastic gradient evaluations for the expectation problem \eqref{eq:stoch_error_obj} for first-order methods, across nonconvex, convex and strongly convex objective functions; see \cref{tab:Complexity_summary}. 
\end{enumerate}

\subsection{Paper Organization}
The paper is organized as follows. In \cref{sec:method}, we describe the proposed algorithm, the adaptive conditions used to control the accuracy of the gradient estimate, and the preliminary assumptions. In \cref{sec:theory}, we present the theoretical analysis, 
covering the nonconvex case in \cref{sec:non_convexity}, the convex case in \cref{sec:convexity}, and the strongly convex case in \cref{sec:strong_convexity}. 
We illustrate the empirical performance of the proposed algorithm in \cref{sec:exp} and provide concluding remarks in \cref{sec:conc}.

\subsection{Notation}
Let $\Rmbb$ denote the set of real numbers and $\Rmbb^d$ denote the set of $d$ dimensional real vectors. Unless otherwise specified, $|\cdot|$ denotes the Euclidean norm of a vector, and $|\cdot|$ denotes either the absolute value of a real number or the cardinality of a set, depending on context. The ceiling function is denoted by $\lceil\cdot\rceil$. Expectation and variance with respect to the distribution $\Pcal$ are denoted as $\Embb[\cdot]$ and $\Var[\cdot]$, respectively. We denote the optimal value for problem \eqref{eq:intro_composite_problem} as $\phi^*$.
\section{Proposed Algorithm} \label{sec:method}

In this section, we present the preliminary assumptions and describe the proposed (accelerated) proximal gradient method, along with the conditions used to adaptively control the accuracy of the gradient estimates. 
For the case of unbiased gradient estimates, we also provide a sequence of sample average approximations that satisfy these conditions.

We begin with the following assumption about the objective function. 

\begin{assumption} \label{ass:base_smoothness}
    The function $f : \Rmbb^d \rightarrow \Rmbb$ is continuously differentiable and has $L$-Lipschitz continuous gradients (i.e., $f$ is $L$-smooth). The function $h : \Rmbb^d \rightarrow \Rmbb \cup \{+\infty\}$ is closed, convex, and proper.
\end{assumption}
A well-known result for functions with Lipschitz continuous gradients (see \cite{nesterov2018lectures}) is 
\begin{equation} \label{eq:smoothness}
    f(a) \leq f(b) + \nabla f(a)^T (b - a) + \tfrac{L}{2} \|b - a\|^2 \quad \forall a, b \in \Rmbb^d.
\end{equation}
Under \cref{ass:base_smoothness}, since $h(x)$ is a convex function, computing the proximal operator \eqref{eq:prox_operator} is well-defined. Specifically, it corresponds to the unique solution of a strongly convex optimization problem.

We now describe the proposed algorithm. At each iteration $k \geq 0$, the algorithm maintains two iterates: the decision variable $x_k$ and an auxiliary variable $y_k$. It computes $g_k$, an estimate of 
$\nabla f(y_k)$, which may be biased. 
The accuracy of this estimate is controlled through 
adaptive conditions described later in this section.
The next iterate $x_{k+1}$ is obtained by taking a step from $y_k$ in the direction of $g_k$, followed by evaluating the proximal operator \eqref{eq:prox_operator}, i.e. $x_{k+1} = \prox_{\alpha_k, h}(y_{k} - \alpha_k g_k)$, where $\alpha_k > 0$ is the step size. Under \cref{ass:base_smoothness}, this is equivalent to
\begin{equation} \label{eq:prox_operator_alternate}
    x_{k+1} = \argmin_{x \in \Rmbb^d} \; f(y_k) + g_k^T (x - y_k) + \tfrac{1}{2 \alpha_k} \|x - y_k\|^2 + h(x), \quad \text{with} \,\, \alpha_k \in \left(0, \tfrac{1}{L}\right].
\end{equation}
The auxiliary variable is then updated in one of two ways. Under the proximal gradient method (referred to as \textbf{Option I}), it is set as $y_{k+1} = x_{k+1}$ and under the accelerated proximal gradient method (\textbf{Option II}), it is updated using the rule $y_{k+1} = x_{k+1} + \beta_{k+1}(x_{k+1} - x_k)$, where $\{ \beta_k \}$ is a user-defined sequence.
The complete procedure is summarized in \cref{alg:stochastic_proximal_gradient} and the following remark.

\begin{algorithm}[H]
    Inputs : Initial iterate $x_{0}$, initial auxiliary iterate $y_0 = x_0$, step size sequence $\{\alpha_k\}$, and acceleration parameter sequence $\{\beta_k\}$.
    \caption{(Accelerated) Proximal Gradient Method with Adaptive Gradient Estimation}
    \begin{algorithmic}[1]
        \For{k = 0, 1, 2, ...}
            \State Compute a gradient estimate $g_k$ of $\nabla f(y_k)$
            \State Proximal step: $x_{k+1} = \prox_{\alpha_k, h}\left(y_k - \alpha_k g_k\right)$
            \State \textbf{Option I:} Set $y_{k+1} = x_{k+1}$
            \State \textbf{Option II:} Set $y_{k+1} = x_{k+1} + \beta_{k+1}(x_{k+1} - x_k)$
        \EndFor
    \end{algorithmic}
    \label{alg:stochastic_proximal_gradient}
\end{algorithm}
\bremark
    We make the following remarks regarding \cref{alg:stochastic_proximal_gradient}.
    \begin{itemize}
        \item \textbf{Gradient Estimate (Line 2):} At each iteration, the algorithm computes a possibly biased gradient estimate $g_k$ for $\nabla f(y_k)$.
        The accuracy of this estimate is adaptively controlled to ensure sufficient progress and is described below.
        \item \textbf{Proximal Step (Line 3):} The proximal step computes the proximal operator defined in \eqref{eq:prox_operator}. We assume that the function $h(x)$ is simple enough such that the proximal operator can be evaluated efficiently.
        \item \textbf{Acceleration Option (Lines 4-5):} The iterate update $y_{k+1}$ follows either \textbf{Option I}, corresponding to the standard proximal gradient method, or \textbf{Option II}, which incorporates acceleration based on \cite{nesterov2013gradient,schmidt2011convergence}, where $\{\beta_k\}$ is a predetermined user-defined sequence.
    \end{itemize}
\eremark

We now introduce some quantities to describe the conditions controlling the accuracy of the gradient estimate.
The reduced gradient at iteration $k \geq 0$ is defined as
\begin{equation} \label{eq:projected_gradient}
    R_{\alpha_k}(y_k) = \tfrac{1}{\alpha_k}\left(y_{k} - x_{k+1}\right).
\end{equation}
The true step computed using $\nabla f(y_k) $ and the corresponding true reduced gradient at iteration $k \geq 0$ are defined as
\begin{equation} \label{eq:true_quantities}
    \hat{x}_{k+1} = \prox_{\alpha_k, h}(y_{k} - \alpha_k \nabla f(y_k)) \quad \text{and} \quad R_{\alpha_k}^{true}(y_k) = \tfrac{1}{\alpha_k}\left(y_{k} - \hat{x}_{k+1}\right),
\end{equation}
respectively. If $\left\| R_{\alpha_k}^{true}(y_k)\right\| = 0$, then $y_k$ is a stationary point for $\phi(x)$ \cite{beiser2023adaptive}. The error in the reduced gradient can be bounded in terms of the gradient estimation error using the contraction property of the proximal operator as follows:
\begin{align*}
    \|R_{\alpha}(y_k) - R_{\alpha}^{true}(y_k)\| &= \tfrac{1}{\alpha_k}\|\prox_{\alpha_k, h}(y_{k} - \alpha_k g_k) - \prox_{\alpha_k, h}(y_{k} - \alpha_k \nabla f(y_k))\| \\
    & \leq \|g_k - \nabla f(y_k)\| \numberthis \label{eq:projected_grad_error}.
\end{align*}
We also define a nested sequence of $\sigma$-algebras $\{\Gcal_k\}$ where $\Gcal_0 = \{x_0\}$ and $\Gcal_k = \{x_0, g_0, g_1, \dots, g_{k-1}\}$ $\forall k \geq 1$. Hence, both $x_k$ and $y_k$ are specified under $\Gcal_k$. We denote the conditional expectation given $\Gcal_k$ as $\Embb_k[\cdot] = \Embb[\cdot | \Gcal_k]$ and the total expectation, (i.e., the expectation given the initial conditions) as $\Embb[\cdot] = \Embb[\cdot|\Gcal_0]$. 

We next introduce the conditions controlling the accuracy of the gradient estimate that guarantee fast convergence, analogous to deterministic methods. We adapt the conditions proposed in \cite{xie2024constrained} for proximal gradient methods, originally based on the well-known \emph{norm condition} \cite{bollapragada2018adaptive,o2024fast,byrd2012sample,carter1991global} for smooth optimization, to our (accelerated) proximal setting, and further introduce relaxations to improve efficiency. The proposed conditions are presented below.

\begin{condition} \label{cond:sampling}
    The gradient estimate $g_k$ in \cref{alg:stochastic_proximal_gradient}, $\forall k \geq 0$, is chosen such that:
    \begin{enumerate}
        \item For the finite-sum problem \eqref{eq:deter_error_obj}: With constants $\eta_k \in [0, 1)$ and $\iota_0, \delta_k \geq 0$, 
        \begin{align*}
            \|g_k - \nabla f(y_{k})\| &\leq \tfrac{\eta_k}{2} \left\|R_{\alpha_k}(y_k)\right\| + \iota_0 \delta_k.
        \end{align*}
        \item For the expectation problem \eqref{eq:stoch_error_obj}: With constants $\tilde{\eta}_k \in [0, 1)$ and $\tilde{\iota}_0, \tilde{\delta}_k \geq 0$, 
        \begin{align*}
            \Embb_k\left[\|g_k - \nabla f(y_k)\|^2\right] \leq \tfrac{\tilde{\eta}_k^2}{4} \left\|\Embb_k\left[R_{\alpha_k}(y_k)\right]\right\|^2 + \tilde{\iota}_0^2\tilde{\delta}_k^2.
        \end{align*}
    \end{enumerate}
\end{condition}
\cref{cond:sampling}, while utilizing sampled quantities on the right-hand side of the inequality, controls the accuracy of the gradient estimate using the true reduced gradient \eqref{eq:true_quantities} as the optimality measure, as shown in \cref{lem:condtion_rearrange} using \eqref{eq:projected_grad_error}. To relax restrictions on the gradient estimation error, we introduce the sequences $\delta_k$ and $\tilde{\delta}_k$ in addition to the optimality measures on the right-hand side of \cref{cond:sampling}; these are appropriately chosen to ensure good performance of the proposed algorithms. We show that \cref{alg:stochastic_proximal_gradient} achieves the optimal complexity for the number of proximal operator evaluations when employing a biased gradient estimate that satisfies \cref{cond:sampling}. When the gradient estimate is unbiased,
we construct sample average approximations, where \cref{cond:sampling} governs the sample size to ensure fast convergence.
To determine the complexity of the number of stochastic gradient evaluations in this case, we introduce an additional assumption regarding the variance (or error) of the stochastic gradients. We then present a lemma that constructs a sequence of unbiased gradient estimates satisfying \cref{cond:sampling}.

\begin{assumption} \label{ass:bounded_variance}
    For any run of \cref{alg:stochastic_proximal_gradient}, $\exists \sigma \geq 0$ such that $\forall k \geq 0$,
    \begin{enumerate}
        \item For the finite-sum problem \eqref{eq:deter_error_obj}: $\|\nabla F(y_k, \xi_i) - \nabla f(y_k)\|^2 \leq \sigma^2$, $\forall i = 1, 2, \dots, N$.
        \item For the expectation problem \eqref{eq:stoch_error_obj}: $\Var\left[\nabla F(y_k, \xi) | \Gcal_k\right] \leq \sigma^2$.
    \end{enumerate}
\end{assumption}
\cref{ass:bounded_variance} is frequently employed to characterize the complexity of the number of stochastic gradient evaluations, as done for finite-sum problems in \cite{pham2020proxsarah,j2016proximal,li2018simple} and for expectation problems in \cite{ghadimi2016accelerated,lan2012optimal,beiser2023adaptive,ghadimi2012optimal}.
Since \cref{alg:stochastic_proximal_gradient} performs a proximal step over the closed and proper function $h(x)$ each iteration, it is reasonable to expect a bounded deviation of the components gradients for the finite-sum problem \eqref{eq:deter_error_obj} and the variance to remain bounded for the expectation problem \eqref{eq:stoch_error_obj}, over the set of iterates.
We now present a sequence of unbiased gradient estimates in the form of sample average approximations that satisfy \cref{cond:sampling} for the expectation problem \eqref{eq:stoch_error_obj}.

\begin{lemma} \label{lem:unbiased_sample_set_bound}
    Suppose Assumptions \ref{ass:base_smoothness} and \ref{ass:bounded_variance} hold in \cref{alg:stochastic_proximal_gradient} for the expectation problem \eqref{eq:stoch_error_obj}. Let $g_k = \tfrac{1}{|S_k|}\sum_{\xi \in S_k} \nabla F(y_k, \xi)$, where $S_k$ is a set of i.i.d. samples from $\Pcal$, independent of $\Gcal_k$, $\forall k \geq 0$. Then, \cref{cond:sampling} is satisfied $\forall k \geq 0$ if
    \begin{equation} \label{eq:unbiased_batch_size}
        |S_k| = \left\lceil\tfrac{\sigma^2}{\tfrac{\tilde{\eta}_k^2}{4} \left\|\Embb_k\left[R_{\alpha_k}(y_k)\right]\right\|^2 + \tilde{\iota}_0^2\tilde{\delta}_k^2} \right\rceil.
    \end{equation}
\end{lemma}
\bproof
In iteration $k \geq 0$, from \cref{ass:bounded_variance} and the definition of $g_k$, the gradient error can be bounded as,
\begin{align*}
    \Embb_k\left[\|g_k - \nabla f(y_k)\|^2\right] = \tfrac{\Var\left[\nabla F(y_k, \xi) | \Gcal_k\right]}{|S_k|} \leq \tfrac{\sigma^2}{|S_k|}.
\end{align*}
From \eqref{eq:unbiased_batch_size}, we have $|S_k| \geq \tfrac{\sigma^2}{\tfrac{\tilde{\eta}_k^2}{4} \left\|\Embb_k\left[R_{\alpha_k}(y_k)\right]\right\|^2 + \tilde{\iota}_0^2\tilde{\delta}_k^2}$. 
Therefore, the gradient error is bounded as,
\begin{equation*}
    \Embb_k\left[\|g_k - \nabla f(y_k)\|^2\right] \leq \tfrac{\sigma^2}{|S_k|} \leq \tfrac{\tilde{\eta}_k^2}{4} \left\|\Embb_k\left[R_{\alpha_k}(y_k)\right]\right\|^2 + \tilde{\iota}_0^2\tilde{\delta}_k^2,
\end{equation*}
satisfying \cref{cond:sampling}.
\eproof

Similar to \cref{lem:unbiased_sample_set_bound}, for the finite-sum problem \eqref{eq:deter_error_obj} under \cref{ass:bounded_variance}, 
at iteration $k \geq 0$, if $g_k = \tfrac{1}{|S_k|}\sum_{\xi \in S_k} \nabla F(y_k, \xi)$, where $S_k \subseteq \Scal$, \cref{cond:sampling} is satisfied if
\begin{equation*}
    |S_k| = \left\lceil \tfrac{N}{\left(1 + \tfrac{\eta_k \left\|R_{\alpha_k}(y_k)\right\| + 2\iota_0 \delta_k}{2\sigma}\right)}\right\rceil,
\end{equation*}
as shown in \cref{lem:deter_sample_set_bound}. While we only present the complexity of the number of stochastic gradient evaluations for the expectation problem \eqref{eq:stoch_error_obj} using unbiased gradient estimates, results for the finite-sum problem \eqref{eq:deter_error_obj} can be established through a similar procedure.
\section{Theoretical Analysis} \label{sec:theory}
In this section, we present the theoretical results for  \cref{alg:stochastic_proximal_gradient}.
We begin with some preliminary results, followed by the analysis for three classes of objective functions: nonconvex (\cref{sec:non_convexity}), general convex (\cref{sec:convexity}), and strongly convex (\cref{sec:strong_convexity}).
For each class, we consider both the finite-sum problem \eqref{eq:deter_error_obj} and the expectation problem \eqref{eq:stoch_error_obj}, and derive the iteration complexity, i.e., the number of proximal operator evaluations required to obtain an $\epsilon > 0$ accurate solution when using biased gradient estimates.
We then refine these results for the expectation problem when the gradient estimate is unbiased and provide the complexity in terms of the number of stochastic gradient evaluations.
All results naturally extend to smooth unconstrained optimization as a special case of composite optimization with $h(x) = 0$.

We begin with a technical descent lemma for problem \eqref{eq:intro_composite_problem} under \cref{ass:base_smoothness}.
\begin{lemma} \label{lem:general_descent_lemma}
    Suppose \cref{ass:base_smoothness} holds. Then, $\forall x \in \Rmbb^d$ and $\forall k \geq 0$, the iterates generated by \cref{alg:stochastic_proximal_gradient} satisfy,
    \begin{align*}
        \phi(x_{k+1}) &\leq f(y_k) + \nabla f(y_k)^T (x - y_k) + (g_k - \nabla f(y_k))^T (x - y_k) + \tfrac{1}{2\alpha_k} \|x - y_k\|^2 + h(x)  \\
        &\quad + (\nabla f(y_k) - g_k)^T (x_{k+1} - y_k) - \left(\tfrac{1}{2 \alpha_k} - \tfrac{L}{2}\right) \|x_{k+1} - y_{k}\|^2. 
    \end{align*}
\end{lemma}
\bproof
    From \eqref{eq:smoothness},
    \begin{align*}
        \phi(x_{k+1}) &\leq f(y_k) + \nabla f(y_k)^T (x_{k+1} - y_k) + \tfrac{L}{2} \|x_{k+1} - y_k\|^2 + h(x_{k+1}) \\
        &= f(y_k) + g_k^T (x_{k+1} - y_k) + \tfrac{1}{2 \alpha_k} \|x_{k+1} - y_k\|^2 + h(x_{k+1})  \\
        &\quad + (\nabla f(y_k) - g_k)^T (x_{k+1} - y_k) - \left(\tfrac{1}{2 \alpha_k} - \tfrac{L}{2}\right) \|x_{k+1} - y_{k}\|^2  \\
        &\leq f(y_k) + g_k^T (x - y_k) + \tfrac{1}{2 \alpha_k} \|x - y_k\|^2 + h(x)  \\
        &\quad + (\nabla f(y_k) - g_k)^T (x_{k+1} - y_k) - \left(\tfrac{1}{2 \alpha_k} - \tfrac{L}{2}\right) \|x_{k+1} - y_{k}\|^2,
    \end{align*}
    where the last inequality follows $\forall x \in \Rmbb^d$ from \eqref{eq:prox_operator_alternate}. Adding and subtracting $\nabla f(y_k)^T (x - y_k)$ on the right-hand side completes the proof.
\eproof

Next, we present a collection of inequalities for the gradient error under \cref{cond:sampling}, which will be frequently used in the subsequent analysis.
\begin{lemma} \label{lem:gradient_error_bounds}
    Suppose \cref{ass:base_smoothness} holds and the gradient estimate $g_k$ satisfies \cref{cond:sampling}.
    Then, $\forall k \geq 0$ in \cref{alg:stochastic_proximal_gradient}:
    \begin{enumerate}
        \item For the finite-sum problem \eqref{eq:deter_error_obj}: 
        \begin{align}
            \|g_k - \nabla f(y_{k})\|^2 &\leq \tfrac{\eta_k^2}{2} \left\|R_{\alpha_k}(y_k)\right\|^2 + 2\iota_0^2 \delta_k^2 , \label{eq:deter_error_alternate}\\
            (\nabla f(y_k) - g_k)^T (x_{k+1} - y_k) &\leq \tfrac{\alpha_k\eta_k}{2} \|R_{\alpha_k}(y_k)\|^2 + \alpha_k \iota_0\delta_k \|R_{\alpha_k}(y_k)\|, \label{eq:deter_error_prod}\\
            (\nabla f(y_k) - g_k)^T (x_{k+1} - y_k) &\leq \tfrac{\alpha_k\left(\eta_k + \iota_0^2\right)}{2} \|R_{\alpha_k}(y_k)\|^2 + \tfrac{\alpha_k \delta_k^2}{2}, \label{eq:deter_error_prod_squared}\\
            \|R_{\alpha_k}^{true}(y_k)\|^2 &\leq 2\left(1 + \tfrac{\eta_k}{2}\right)^2\|R_{\alpha_k}(y_k)\|^2 + 2 \iota_0^2\delta_k^2. \label{eq:deter_true_projected_grad_bound}
        \end{align}
        \item For the expectation problem \eqref{eq:stoch_error_obj}: 
        \begin{align}
            \Embb_k[\|g_k - \nabla f(y_k)\|] &\leq \tfrac{\tilde{\eta}_k}{2} \sqrt{\Embb_k\left[ \|R_{\alpha_k}(y_k)\|^2\right]} + \tilde{\iota}_0\tilde{\delta}_k, \label{eq:stoch_error_alternate}\\
            \Embb_k[(\nabla f(y_k) - g_k)^T (x_{k+1} - y_k)] &\leq \tfrac{\alpha_k \tilde{\eta}_k}{2} \Embb_k\left[\|R_{\alpha_k}(y_k)\|^2\right] + \alpha_k \tilde{\iota}_0\tilde{\delta}_k \sqrt{\Embb_k\left[\|R_{\alpha_k}(y_k)\|^2\right]}, \label{eq:stoch_error_prod}\\
            \Embb_k[(\nabla f(y_k) - g_k)^T (x_{k+1} - y_k)] &\leq \tfrac{\alpha_k \left(\tilde{\eta}_k + \tilde{\iota}_0^2\right)}{2} \Embb_k\left[\|R_{\alpha_k}(y_k)\|^2\right] + \tfrac{\alpha_k \tilde{\delta}_k^2}{2}, \label{eq:stoch_error_prod_squared}\\ 
            \left\|R_{\alpha_k}^{true}(y_k)\right\|^2 &\leq 2\left(1 + \tfrac{\tilde{\eta}_k^2}{4}\right) \left\|\Embb_k\left[R_{\alpha_k}(y_k)\right]\right\|^2 + 2 \tilde{\iota}_0^2 \tilde{\delta}_k^2. \label{eq:complexity_bound_help}
        \end{align}
    \end{enumerate}
\end{lemma}
\bproof
    For the finite-sum problem \eqref{eq:deter_error_obj}, from \cref{cond:sampling},
    \begin{align*}
        \|g_k - \nabla f(y_{k})\|^2 \leq \left(\tfrac{\eta_k}{2} \left\|R_{\alpha_k}(y_k)\right\| + \iota_0 \delta_k\right)^2 \leq \tfrac{\eta_k^2}{2} \left\|R_{\alpha_k}(y_k)\right\|^2 + 2 \iota_0^2 \delta_k^2,
    \end{align*}
    where the last inequality follows from the identity $(a + b)^2 \leq 2a^2 + 2b^2$, completing the proof of 
    \eqref{eq:deter_error_alternate}.
    For \eqref{eq:deter_error_prod}, by the Cauchy-Schwartz inequality,
    \begin{align*}
        (\nabla f(y_k) - g_k)^T (x_{k+1} - y_k) 
        &\leq \|\nabla f(y_k) - g_k\| \|x_{k+1} - y_k\| = \alpha_k \|\nabla f(y_k) - g_k\| \|R_{\alpha_k}(y_k)\| \\
        &\leq \tfrac{\alpha_k\eta_k}{2} \|R_{\alpha_k}(y_k)\|^2 + \alpha_k \iota_0 \delta_k\|R_{\alpha_k}(y_k)\|,
    \end{align*}
    where the second inequality follows from \cref{cond:sampling}, completing the proof. 
    Applying the identity $ab \leq \tfrac{a^2 + b^2}{2}$ to the right-hand side of the above inequality yields \eqref{eq:deter_error_prod_squared}.
    Finally for \eqref{eq:deter_true_projected_grad_bound},
    \begin{align*}
        \|R_{\alpha_k}^{true}(y_k)\| 
        \leq \|R_{\alpha_k}(y_k)\| + \|R_{\alpha_k}^{true}(y_k) - R_{\alpha_k}(y_k)\|
        \leq \|R_{\alpha_k}(y_k)\| + \|\nabla f(y_k) - g_k\|,
    \end{align*}
    where the second inequality follows from \eqref{eq:projected_grad_error}. 
    Substituting \cref{cond:sampling} in the above inequality, then squaring both sides and applying the identity $(a + b)^2 \leq 2a^2 + 2b^2$, completes the proof of  
    \eqref{eq:deter_true_projected_grad_bound}.

    For the expectation problem \eqref{eq:stoch_error_obj}, from \cref{cond:sampling},
    \begin{align*}
        (\Embb_k[\left\|g_k - \nabla f(y_k)\right\|])^2
        &\leq \Embb_k\left[\|g_k - \nabla f(y_k)\|^2\right] 
        \leq \tfrac{\tilde{\eta}_k^2}{4} \left\|\Embb_k\left[R_{\alpha_k}(y_k)\right]\right\|^2 + \tilde{\iota}_0^2\tilde{\delta}_k^2 \\
        &\leq \tfrac{\tilde{\eta}_k^2}{4} \Embb_k\left[\left\|R_{\alpha_k}(y_k)\right\|^2\right] + \tilde{\iota}_0^2\tilde{\delta}_k^2
        \leq \left(\tfrac{\tilde{\eta}_k}{2} \sqrt{\Embb_k\left[\left\|R_{\alpha_k}(y_k)\right\|^2\right]} + \tilde{\iota}_0\tilde{\delta}_k\right)^2,
    \end{align*} 
    yielding \eqref{eq:stoch_error_alternate}. For \eqref{eq:stoch_error_prod}, from the Cauchy-Schwartz inequality,
    \begin{align*}
        \Embb_k\left[(\nabla f(y_k) - g_k)^T (x_{k+1} - y_k)\right]
        &\leq \sqrt{\Embb_k\left[\|\nabla f(y_k) - g_k\|^2\right]}\sqrt{\Embb_k\left[\|x_{k+1} - y_k\|^2\right]} \\
        &= \alpha_k \sqrt{\Embb_k\left[\|\nabla f(y_k) - g_k\|^2\right]}\sqrt{\Embb_k\left[\|R_{\alpha_k}(y_k)\|^2\right]} \\
        &\leq \alpha_k \left(\tfrac{\tilde{\eta}_k}{2} \sqrt{\Embb_k\left[\|R_{\alpha_k}(y_k)\|^2\right]} + \tilde{\iota}_0\tilde{\delta}_k\right)\sqrt{\Embb_k\left[\|R_{\alpha_k}(y_k)\|^2\right]},
    \end{align*}
    where the equality follows from \eqref{eq:projected_gradient} and the second inequality follows from \cref{cond:sampling}, completing the proof.
    Applying the identity $ab \leq \tfrac{a^2 + b^2}{2}$ to \eqref{eq:stoch_error_prod} yields \eqref{eq:stoch_error_prod_squared}.
    Finally for \eqref{eq:complexity_bound_help}, using the identity $(a + b)^2 \leq 2a^2 + 2b^2$,
    \begin{align*}
        \left\|R_{\alpha_k}^{true}(y_k)\right\|^2 
        &\leq 2\left\|\Embb_k\left[R_{\alpha_k}(y_k)\right]\right\|^2 + 2\left\|\Embb_k\left[R_{\alpha_k}^{true}(y_k) - R_{\alpha_k}(y_k)\right]\right\|^2 \\
        &\leq 2\left\|\Embb_k\left[R_{\alpha_k}(y_k)\right]\right\|^2 + 2\Embb_k\left[\left\|R_{\alpha_k}^{true}(y_k) - R_{\alpha_k}(y_k)\right\|^2\right] \\
        &\leq 2\left\|\Embb_k\left[R_{\alpha_k}(y_k)\right]\right\|^2 + 2\Embb_k\left[\|\nabla f(y_k) - g_k\|^2\right],
    \end{align*}
    where the second inequality follows from Jenson's inequality, the third inequality follows from \eqref{eq:projected_grad_error}, and further using \cref{cond:sampling} completes the proof.
\eproof
\subsection{Nonconvex Objective Function} \label{sec:non_convexity}

In this section, we present the theoretical analysis for \cref{alg:stochastic_proximal_gradient} when the smooth function $f(x)$, and thus the composite function $\phi(x)$, is nonconvex.
The analysis is limited to \textbf{Option I} in \cref{alg:stochastic_proximal_gradient}, as \textbf{Option II} does not yield any improvements for nonconvex objective functions up to constant factors, as noted in \cite{ghadimi2016accelerated}.
We first establish the convergence of \cref{alg:stochastic_proximal_gradient} when using a biased gradient estimate, followed by the complexity of number of proximal operator evaluations. We then establish the complexity of the number of stochastic gradient evaluation for the expectation problem when using an unbiased gradient estimate.

\begin{theorem} \label{th:convergence_non_convex}
    Suppose \cref{ass:base_smoothness} holds and the gradient estimate $g_k$ satisfies \cref{cond:sampling}. Then, for \cref{alg:stochastic_proximal_gradient} with \textbf{Option I}:
    \begin{enumerate}
        \item For the finite-sum problem \eqref{eq:deter_error_obj},  if the parameters in \cref{cond:sampling} are chosen such that $\{\eta_k\} = \eta \in [0, 1)$, $\iota_0 \in \left[0, \sqrt{\tfrac{1-\eta}{2}}\right)$ and $\sum_{k=0}^{\infty} \delta_k^2 < \infty$, and the step size is chosen such that $\{\alpha_k\} = \alpha \leq \tfrac{1 - \eta}{2L}$, then $\{x_k\}$ converges to a stationary point with
        $\min_{k = 0, \dots, K-1} \|R_{\alpha_k}^{true}(x_k)\|^2 = \Ocal\left(\tfrac{1}{K}\right)$, $\forall K \geq 1$.
        \item For the expectation problem \eqref{eq:stoch_error_obj},  if the parameters in \cref{cond:sampling} are chosen such that $\{\tilde{\eta}_k\} = \tilde{\eta} \in [0, 1)$, $\tilde{\iota}_0 \in \left[0, \sqrt{\tfrac{1-\tilde{\eta}}{2}}\right)$ and $\sum_{k=0}^{\infty} \tilde{\delta}_k^2 < \infty$, and the step size is chosen such that $\{\alpha_k\} = \alpha \leq \tfrac{1 - \tilde{\eta}}{2L}$, then $\{x_k\}$ converges to a stationary point in expectation with $\min_{k = 0, \dots, K-1} \Embb\left[\left\|R_{\alpha_k}^{true}(x_k)\right\|^2\right] = \Ocal\left(\tfrac{1}{K}\right)$, $\forall K \geq 1$.
    \end{enumerate}
\end{theorem}
\bproof
With $y_k = x_k$ under \textbf{Option I}, consider the result from \cref{lem:general_descent_lemma}. Substituting
$x = x_k$ and using \eqref{eq:projected_gradient} yields,
\begin{equation} \label{eq:non_convex_unbiased_1}
    \phi(x_{k+1}) \leq \phi(x_k) + (\nabla f(x_k) - g_k)^T (x_{k+1} - x_k) - \alpha_k^2 \left(\tfrac{1}{2 \alpha_k} - \tfrac{L}{2}\right) \|R_{\alpha_k}(x_k)\|^2.
\end{equation}
For the finite-sum problem \eqref{eq:deter_error_obj}, under the defined parameters, substituting \eqref{eq:deter_error_prod_squared} from \cref{lem:gradient_error_bounds} into \eqref{eq:non_convex_unbiased_1} yields,
\begin{align*}
    \phi(x_{k+1}) &\leq \phi(x_k) + \tfrac{\alpha\left(\eta + \iota_0^2\right)}{2} \|R_{\alpha}(x_k)\|^2 + \tfrac{\alpha \delta_k^2}{2} - \alpha^2 \left(\tfrac{1}{2 \alpha} - \tfrac{L}{2}\right) \|R_{\alpha}(x_k)\|^2 \\
    &= \phi(x_k) - \alpha \left[\tfrac{1}{2} - \tfrac{\alpha L}{2} - \tfrac{\eta}{2} - \tfrac{\iota_0^2}{2}\right]\|R_{\alpha}(x_k)\|^2  + \tfrac{\alpha \delta_k^2}{2}.
\end{align*}
Rearranging the above inequality and substituting $c(\eta, \alpha) = \tfrac{1}{2} - \tfrac{\alpha L}{2} - \tfrac{\eta}{2} - \tfrac{\iota_0^2}{2}$, we get,
\begin{align*}
    \|R_{\alpha}(x_k)\|^2
    &\leq \tfrac{\phi(x_{k}) - \phi(x_{k+1})}{\alpha c(\alpha, \eta)} + \tfrac{\delta_k^2}{2c(\alpha, \eta)}.
\end{align*}
Under the defined parameters, $c(\eta, \alpha) \geq \tfrac{1}{2} - \tfrac{(1 - \eta)}{4} - \tfrac{\eta}{2} - \tfrac{\iota_0^2}{2} = \tfrac{1 - \eta}{4} - \tfrac{\iota_0^2}{2} > \tfrac{1 - \eta}{4} - \tfrac{1 - \eta}{4} = 0$. Thus, substituting the above bound into \eqref{eq:deter_true_projected_grad_bound} from \cref{lem:gradient_error_bounds} yields,
\begin{align*}
    \|R_{\alpha}^{true}(x_k)\|^2 &\leq 2\left(1 + \tfrac{\eta}{2}\right)^2\left[\tfrac{\phi(x_{k}) - \phi(x_{k+1})}{\alpha c(\alpha, \eta)} + \tfrac{\delta_k^2}{2c(\alpha, \eta)}\right] + 2 \iota_0^2\delta_k^2.
\end{align*}
The telescoping sum of the above inequality for $k = 0, \dots, K-1$ yields,
\begin{align*}
    \sum_{k = 0}^{K-1} \|R_{\alpha}^{true}(x_k)\|^2 &\leq 2\left(1 + \tfrac{\eta}{2}\right)^2\left[\tfrac{\phi(x_{0}) - \phi(x_{K})}{\alpha c(\alpha, \eta)} + \sum_{k=0}^{K-1} \tfrac{\delta_k^2}{2c(\alpha, \eta)}\right] + 2 \iota_0^2 \sum_{k=0}^{K-1}\delta_k^2.
\end{align*}
Rearranging the above inequality and using $\phi(x_K) \geq \phi^*$, we get, 
\begin{align*}
    \min_{k = 0, \dots,K-1} \|R_{\alpha}^{true}(x_k)\|^2 &\leq \tfrac{1}{K} \left\{2\left(1 + \tfrac{\eta}{2}\right)^2\left[\tfrac{\phi(x_{0}) - \phi^*}{\alpha c(\alpha, \eta)}\right] + \left[ \tfrac{\left(2 + \eta\right)^2}{4c(\alpha, \eta)} + 2 \iota_0^2\right] \sum_{k=0}^{K-1}\delta_k^2\right\},
\end{align*}
where all terms within the curly brackets on the right-hand side are bounded due to the condition $\sum_{k=0}^{\infty} \delta_k^2 < \infty$, completing the proof for the finite-sum problem \eqref{eq:deter_error_obj}.

For the expectation problem \eqref{eq:stoch_error_obj}, under the defined parameters, taking a conditional expectation of \eqref{eq:non_convex_unbiased_1} given $\Gcal_k$ and substituting \eqref{eq:stoch_error_prod_squared} from \cref{lem:gradient_error_bounds} yields,
\begin{align*}
    \Embb_k\left[\phi(x_{k+1})\right] &\leq \phi(x_k) + \tfrac{\alpha \left(\tilde{\eta} + \tilde{\iota}_0^2\right)}{2} \Embb_k\left[\|R_{\alpha}(x_k)\|^2\right] + \tfrac{\alpha \tilde{\delta}_k^2}{2} - \alpha^2 \left(\tfrac{1}{2 \alpha} - \tfrac{L}{2}\right) \Embb_k \left[\|R_{\alpha}(x_k)\|^2\right] \\
    &= \phi(x_k) + \tfrac{\alpha \tilde{\delta}_k^2}{2} - \alpha \left[\tfrac{1}{2} - \tfrac{\alpha L}{2} - \tfrac{\tilde{\eta}}{2} - \tfrac{\tilde{\iota}_0^2}{2}\right]\Embb_k \left[\left\|R_{\alpha}(x_k)\right\|^2\right].
\end{align*}
Rearranging the above inequality and substituting $\tilde{c}(\alpha, \tilde{\eta}) =  \tfrac{1}{2} - \tfrac{\alpha L}{2} - \tfrac{\tilde{\eta}}{2} - \tfrac{\tilde{\iota}_0^2}{2}$, we get,
\begin{align*}
    \Embb_k \left[\left\|R_{\alpha}(x_k)\right\|^2\right] &\leq \tfrac{\phi(x_{k}) - \Embb_k[\phi(x_{k+1})]}{\alpha \tilde{c}(\alpha, \tilde{\eta})} + \tfrac{\tilde{\delta}_k^2}{2\tilde{c}(\alpha, \tilde{\eta})}.
\end{align*}
Under the defined parameters, $\tilde{c}(\alpha, \tilde{\eta}) > 0$.
Since $\left\|\Embb_k[R_{\alpha_k}(x_k)]\right\|^2 \leq \Embb_k\left[\left\|R_{\alpha_k}(x_k)\right\|^2\right]$ by Jenson's inequality, substituting the above bound into \eqref{eq:complexity_bound_help} from \cref{lem:gradient_error_bounds} yields,
\begin{align*}
    \left\|R_{\alpha}^{true}(x_k)\right\|^2 &\leq 2\left(1 + \tfrac{\tilde{\eta}^2}{4}\right) \left[\tfrac{\phi(x_{k}) - \Embb_k[\phi(x_{k+1})]}{\alpha \tilde{c}(\alpha, \tilde{\eta})} + \tfrac{\tilde{\delta}_k^2}{2\tilde{c}(\alpha, \tilde{\eta})}\right] + 2 \tilde{\iota}_0^2 \tilde{\delta}_k^2.
\end{align*}

Taking the total expectation of the above inequality and summing telescopically over $k = 0, \ldots, K-1$ yields,
\begin{align*}
    \sum_{k=0}^{K-1} \Embb\left[\left\|R_{\alpha}^{true}(x_k)\right\|^2\right] &\leq 2\left(1 + \tfrac{\tilde{\eta}^2}{4}\right) \left[\tfrac{\phi(x_{0}) - \Embb\left[\phi(x_{K})\right]}{\alpha \tilde{c}(\alpha, \tilde{\eta})} + \sum_{k=0}^{K-1} \tfrac{\tilde{\delta}_k^2}{2\tilde{c}(\alpha, \tilde{\eta})}\right] + \sum_{k=0}^{K-1} 2 \tilde{\iota}_0^2 \tilde{\delta}_k^2.
\end{align*}
Rearranging the above inequality and using $\phi(x_K) \geq \phi^*$, we get,
\begin{align*}
    \min_{k=0, \dots, K-1} \Embb\left[\left\|R_{\alpha}^{true}(x_k)\right\|^2\right] &\leq
    \tfrac{1}{K} \left\{2\left(1 + \tfrac{\tilde{\eta}^2}{4}\right) \left[\tfrac{\phi(x_{0}) - \phi^*}{\alpha \tilde{c}(\alpha, \tilde{\eta})}\right] + 
    \left[\tfrac{4 + \tilde{\eta}^2}{4\tilde{c}(\alpha, \tilde{\eta})} + 2 \tilde{\iota}_0^2 \right]\sum_{k=0}^{K-1}  \tilde{\delta}_k^2\right\}.
\end{align*}
where all terms within the curly brackets on the right-hand side are bounded due to the condition $\sum_{k=0}^{\infty} \tilde{\delta}_k^2 <  \infty$, completing the proof for the expectation problem \eqref{eq:stoch_error_obj}.
\eproof

\cref{th:convergence_non_convex} establishes a sublinear convergence rate for \cref{alg:stochastic_proximal_gradient} with \textbf{Option I} for nonconvex objective functions, when using a possibly biased gradient estimate satisfying \cref{cond:sampling}. 
Thus, an $\epsilon > 0$ accurate solution, i.e., $\left\|R_{\alpha_k}^{true}(x_k)\right\|^2 \leq \epsilon$ for the finite-sum problem \eqref{eq:deter_error_obj} and $\Embb\left[\left\|R_{\alpha_k}^{true}(x_k)\right\|^2\right] \leq \epsilon$ for the expectation problem \eqref{eq:stoch_error_obj}, 
can be achieved in $\Ocal\left(\tfrac{1}{\epsilon}\right)$ iterations (proximal operator evaluations), matching the complexity bounds for deterministic first-order methods \cite{carmon2020lower, carmon2021lower}.
We now present the complexity for the number of stochastic gradient evaluations for \cref{alg:stochastic_proximal_gradient} for the expectation problem \eqref{eq:stoch_error_obj} with a nonconvex objective function when the gradient estimate is unbiased.

\begin{theorem} \label{th:sample_complexity_non_convex}
    Suppose Assumptions \ref{ass:base_smoothness} and \ref{ass:bounded_variance} hold, and \cref{cond:sampling} is satisfied for the expectation problem \eqref{eq:stoch_error_obj} via the unbiased gradient estimate in \cref{lem:unbiased_sample_set_bound}.
    Then, for \cref{alg:stochastic_proximal_gradient} with \textbf{Option I}, if the parameters in \cref{cond:sampling} are chosen such that $\{\tilde{\eta}_k\} = \tilde{\eta} \in (0, 1)$, $\tilde{\iota}_0 \in \left[0, \sqrt{\tfrac{1-\tilde{\eta}}{2}}\right)$, $\tilde{\delta}_k = \tfrac{1}{(k+1)^{1 + \nu}}$ $\forall k \geq 0$ with $\nu > 0$, and the step size is chosen such that $\{\alpha_k\} = \alpha \leq \tfrac{1 - \tilde{\eta}}{2L}$, 
    a solution satisfying 
    $\min\left\{\Embb\left[\left\|R^{true}_{\alpha_k}(x_k)\right\|^2\right], \left\|R_{\alpha_k}^{true}(x_k)\right\|^2\right\} \leq \epsilon$
    with $\epsilon > 0$ is achieved in $\Ocal \left(\epsilon^{-(2 + \nu)}\right)$ stochastic gradient evaluations. 
    If $\{\tilde{\delta}_k\} = 0$, this improves to $\Ocal \left(\epsilon^{-2}\right)$.
\end{theorem}
\bproof
Let $K_\epsilon \geq 1$ be the first iteration that achieves the desired solution accuracy. Hence, $\left\|R_{\alpha_k}^{true}(x_k)\right\|^2 > \epsilon$ $\forall k \leq K_\epsilon - 1$ and from \eqref{eq:complexity_bound_help}, 
$\left\|\Embb_k\left[R_{\alpha_k}(x_k)\right]\right\|^2 > \left(1 + \tfrac{\tilde{\eta}^2}{4}\right)^{-1}\left[\tfrac{\epsilon}{2} - \tilde{\iota}_0^2\tilde{\delta}_k^2\right]$ $\forall k \leq K_\epsilon - 1$.
Thus, the total number of stochastic gradient evaluations can be bounded using \cref{lem:unbiased_sample_set_bound} as, 
\begin{align*}
    \sum_{k=0}^{K_\epsilon - 1}|S_k| 
    &= \sum_{k=0}^{K_\epsilon - 1} \left\lceil\tfrac{\sigma^2}{\tfrac{\tilde{\eta}_k^2}{4} \left\|\Embb_k\left[R_{\alpha_k}(x_k)\right]\right\|^2 + \tilde{\iota}_0^2\tilde{\delta}_k^2} \right\rceil \leq \sum_{k=0}^{K_\epsilon - 1} \tfrac{2\sigma^2 (4 + \tilde{\eta}^2)}{\tilde{\eta}^2 \epsilon + 8\tilde{\iota}_0^2\tilde{\delta}_k^2}  + 1 
    \leq \sum_{k=0}^{K_\epsilon - 1} \tfrac{2\sigma^2(4 + \tilde{\eta}^2)}{\tilde{\eta}^2 \epsilon} + \tfrac{\sigma^2 (4 + \tilde{\eta}^2)}{4\tilde{\iota}_0^2\tilde{\delta}_k^2}  + K_{\epsilon}.
\end{align*}
For \cref{alg:stochastic_proximal_gradient} with \textbf{Option I}, $K_\epsilon$ is at most $\Ocal\left(\tfrac{1}{\epsilon}\right)$ from \cref{th:convergence_non_convex}, yielding,
\begin{align*}
    \sum_{k=0}^{K_\epsilon - 1}|S_k| 
    &\leq \tfrac{2\sigma^2 (4 + \tilde{\eta}^2) }{\tilde{\eta}^2 \epsilon}K_{\epsilon}  + \tfrac{\sigma^2 (4 + \tilde{\eta}^2)}{4\tilde{\iota}_0^2}K_{\epsilon}^{2 + \nu} + K_{\epsilon} = \Ocal\left(\tfrac{1}{\epsilon^{2 + \nu}}\right).
\end{align*}
Following the same procedure, if $\{\tilde{\delta}_k\} = 0$, $\sum_{k=0}^{K_\epsilon - 1}|S_k| \leq \tfrac{2\sigma^2 (4 + \tilde{\eta}^2) }{\tilde{\eta}^2 \epsilon}K_{\epsilon}  + K_{\epsilon} = \Ocal\left(\tfrac{1}{\epsilon^{2}}\right)$.

\eproof

\cref{th:sample_complexity_non_convex} matches the optimal complexity for the number of stochastic gradients for the expectation problem \eqref{eq:stoch_error_obj} over nonconvex objective functions \cite{ghadimi2012optimal}.
We conclude this section with a corollary to \cref{th:sample_complexity_non_convex}, using a definition of an $\epsilon$-accurate solution similar to that in \cite{ghadimi2016accelerated},  under the parameter setting $\{\tilde{\eta}_k\} = 0$.

\begin{corollary} \label{cor:sample_complexity_non_convex}
    Suppose the conditions in \cref{th:sample_complexity_non_convex} hold.
    If the parameters in \cref{cond:sampling} are chosen as $\{\tilde{\eta_k}\} = 0$, $\tilde{\iota}_0 \in \left(0, \tfrac{1}{\sqrt{2}}\right)$ and $\tilde{\delta}_k = \tfrac{1}{(k+1)^{1 + \nu}}$ $\forall k \geq 0$ with $\nu > 0$, and the step size is chosen such that $\{\alpha_k\} = \alpha \leq \tfrac{1}{2L}$, a solution satisfying $\Embb\left[\left\|R^{true}_{\alpha_k}(x_k)\right\|^2\right] \leq \epsilon$ with $\epsilon > 0$ is achieved in $\Ocal \left(\epsilon^{-(2 + \nu)}\right)$ stochastic gradient evaluations. 
\end{corollary}
\bproof
The proof follows from the same procedure as \cref{th:sample_complexity_non_convex}.
\eproof
The parameter settings in \cref{cor:sample_complexity_non_convex} reduce \cref{cond:sampling} to maintaining a predetermined sequence of 
gradient estimation errors, similar to \cite{schmidt2011convergence}.
\subsection{General Convex Objective Function} \label{sec:convexity}

In this section, we provide the theoretical analysis of \cref{alg:stochastic_proximal_gradient} when the smooth function $f(x)$, and thus the composite function $\phi(x)$, is convex.
We begin by stating the basic assumptions and definitions, along with some mathematical identities that will be used throughout the analysis. 

\begin{assumption} \label{ass:convexity}
    The function $f : \Rmbb^d \rightarrow \Rmbb$ is $L$-smooth and convex, and the function $h : \Rmbb^d \rightarrow \Rmbb \cup \{+\infty\}$ is closed, convex, and proper.
\end{assumption}
Under \cref{ass:convexity}, let $x^* \in \Rmbb^d$ denote an optimal solution where the optimal value $\phi^*$ is attained. Since $f(x)$ is differentiable and convex, from \cite{nesterov2018lectures},
\begin{equation} \label{eq:convexity_alternate}
    f(b) \geq f(a) + \nabla f(a)^T (b - a) \quad \forall a, b \in \Rmbb^d.
\end{equation}

For \cref{alg:stochastic_proximal_gradient} with \textbf{Option II}, under \cref{ass:convexity}, we define the sequence $\{\beta_k\}$ and two additional sequences $\{\theta_k\}$ and $\{v_k\}$ for the analysis as,
\begin{equation} \label{eq:acc_defintions_convex}
    \beta_k = \tfrac{k-1}{k+2} \quad \forall k \geq 1, \quad \theta_k = \tfrac{2}{k+1} \quad \forall k \geq 0, \text{ and} \quad v_{k} = x_{k-1} + \tfrac{1}{\theta_{k}} (x_{k} - x_{k-1}) \quad \forall k \geq 1,
\end{equation}
with $v_0 = x_0$. Under these definitions, $y_k$ can be expressed as,
\begin{align*}
    y_{k} 
    &= x_k + \tfrac{k-1}{k+2}\left(x_k - x_{k-1}\right)
    = \left(1 - \tfrac{2}{k+2}\right)x_k + \tfrac{2}{k + 2}\left(x_{k-1} + \tfrac{k+1}{2} (x_{k} - x_{k-1})\right) \\
    &= (1 - \theta_{k+1})x_k + \theta_{k+1}v_{k} \quad \forall k \geq 0. \numberthis \label{eq:y_k_alternate_convex}
\end{align*}
From \eqref{eq:acc_defintions_convex} and \eqref{eq:y_k_alternate_convex}, an alternative update form for $\{v_k\}$ can be derived as,
\begin{align*}
    v_{k+1} 
    & = x_k + \tfrac{1}{\theta_{k+1}} \left(x_{k+1} - x_k \right) = v_k - \tfrac{1}{\theta_{k+1}} \left((1 - \theta_{k+1})x_k + \theta_{k+1}v_{k} - x_{k+1}\right) \\
    &= v_k - \tfrac{1}{\theta_{k+1}} \left(y_{k} - x_{k+1}\right) \quad \forall k \geq 0. \numberthis \label{eq:v_k_alternate_convex}
\end{align*}
Finally, we also introduce a useful identity for $\{\theta_k\}$ as,
\begin{equation} \label{eq:acc_theta_identity_convex}
    \tfrac{(1 - \theta_k)}{\theta_k^2} = \tfrac{(k- 1)(k+1)}{4} = \tfrac{k^2 - 1}{4} \leq \tfrac{k^2}{4} = \tfrac{1}{\theta_{k-1}^2} \quad \forall k \geq 1.
\end{equation}

We now establish a general descent lemma under \cref{ass:convexity}.

\begin{lemma} \label{lem:general_descent_convex}
    Suppose \cref{ass:convexity} holds. Then, $\forall z \in \Rmbb^d$ and $\forall k \geq 0$, the iterates generated by \cref{alg:stochastic_proximal_gradient} satisfy,
    \begin{align*}
        \phi(x_{k+1}) &\leq \phi(z) + (g_k - \nabla f(y_k))^T (z - y_{k}) + (\nabla f(y_k) - g_k)^T (x_{k+1} - y_k)\\
        &\quad + \left[\tfrac{L}{2} - \tfrac{1}{\alpha_k}\right] \|x_{k+1} - y_k\|^2 + \tfrac{1}{\alpha_k}(y_k - x_{k+1} )^T(y_k - z).
    \end{align*}
\end{lemma}
\bproof
From \eqref{eq:smoothness},
\begin{align*}
    \phi(x_{k+1}) 
    &\leq f(y_k) + \nabla f(y_k)^T (x_{k+1} - y_k) + \tfrac{L}{2} \|x_{k+1} - y_k\|^2 + h(x_{k+1}) \\
    &\leq f(z) - \nabla f(y_k)^T (z - y_k) + \nabla f(y_k)^T (x_{k+1} - y_k) + \tfrac{L}{2} \|x_{k+1} - y_k\|^2 + h(x_{k+1}) \\
    &\leq f(z) - \nabla f(y_k)^T (z - y_{k}) + \nabla f(y_k)^T (x_{k+1} - y_k) + \tfrac{L}{2} \|x_{k+1} - y_k\|^2 \\
    &\quad + h(z) - \left(\tfrac{y_k - x_{k+1}}{\alpha_k} - g_k \right)^T(z - x_{k+1}) \\
    &= \phi(z) - \nabla f(y_k)^T (z - y_{k}) + \nabla f(y_k)^T (x_{k+1} - y_k) + \tfrac{L}{2} \|x_{k+1} - y_k\|^2 \\
    &\quad - \left(\tfrac{y_k - x_{k+1}}{\alpha_k} - g_k \right)^T(z - y_k + y_k - x_{k+1}),
\end{align*}
where the second inequality follows from \eqref{eq:convexity_alternate} $\forall z \in \Rmbb^d$, and the third inequality follows from the convexity of $h(x)$ and the definition \eqref{eq:prox_operator_alternate} for $x_{k+1}$ which yields $0 \in g_k +  \partial h(x_{k+1}) + \tfrac{x_{k+1} - y_k}{\alpha_k}$. Rearranging the terms in the last equality completes the proof.
\eproof

Using \cref{lem:general_descent_convex}, we establish recursive bounds on the optimality gap in function value for \cref{alg:stochastic_proximal_gradient} with \textbf{Option I} and \textbf{Option II}, under \cref{ass:convexity}.

\begin{lemma} \label{lem:general_descent_convex_refined}
    Suppose \cref{ass:convexity} holds.
    \begin{enumerate}
        \item For \cref{alg:stochastic_proximal_gradient} with \textbf{Option I}, $\forall k\geq 0$, the iterates satisfy,
        \begin{align*}
            \phi(x_{k+1})-\phi^*
            &\leq (g_k - \nabla f(x_k))^T (x^* - x_{k}) + (\nabla f(x_k) - g_k)^T (x_{k+1} - x_k) \\
            &\quad + \tfrac{1}{2\alpha_k}\left[ \|x_k - x^*\|^2 - \|x_{k+1} - x^*\|^2 - \left(1 - \alpha_k L\right)\|x_k - x_{k+1}\|^2\right].
        \end{align*}
        \item For \cref{alg:stochastic_proximal_gradient} with \textbf{Option II}, $\forall k\geq 0$, the iterates satisfy,
        \begin{align*}
            \phi(x_{k+1}) - \phi^* 
            &\leq (1 - \theta_{k+1})(\phi(x_k) - \phi^*) + \theta_{k+1}(g_k - \nabla f(y_k))^T (x^* -  v_k) \\
            &\quad + (\nabla f(y_k) - g_k)^T (x_{k+1} - y_k) \\
            &\quad + \tfrac{\theta_{k+1}^2}{2\alpha_k}\left[\|v_k - x^*\|^2 - \|v_{k+1} - x^*\|^2 - \left(1 - \alpha_k L\right)\|v_k - v_{k+1}\|^2\right].
        \end{align*}
    \end{enumerate}
\end{lemma}
\bproof
    For \cref{alg:stochastic_proximal_gradient} with \textbf{Option I}, consider the result from \cref{lem:general_descent_convex} with $y_k = x_k$ and $z = x^*$,
    \begin{align*}
        \phi(x_{k+1})-\phi^*
        &\leq (g_k - \nabla f(x_k))^T (x^* - x_{k}) + (\nabla f(x_k) - g_k)^T (x_{k+1} - x_k) \\
        &\quad + \left[\tfrac{L}{2} - \tfrac{1}{\alpha_k}\right] \|x_{k+1} - x_k\|^2 + \tfrac{1}{\alpha_k}(x_k - x_{k+1} )^T(x_k - x^*).
    \end{align*}
    Applying \cref{lem:a_b_Appendix} with $a_1 = x_k - x_{k+1}$, $a_2 = x_k - x^*$ and $c = \tfrac{\alpha_k L}{2}$ completes the proof for \textbf{Option I}.

    For \cref{alg:stochastic_proximal_gradient} with \textbf{Option II}, consider the result from \cref{lem:general_descent_convex} with $z = \theta_{k+1}x^* + (1 - \theta_{k+1})x_k$. Since  
    $\theta_{k+1} \in [0, 1]$ $\forall k \geq 0$, from \cref{ass:convexity},
    \begin{align*}
        \phi(x_{k+1}) - \phi^* 
        &\leq (1 - \theta_{k+1})(\phi(x_k) - \phi^*) + (g_k - \nabla f(y_k))^T (\theta_{k+1}x^* + (1 - \theta_{k+1})x_k - y_{k}) \\
        &\quad + (\nabla f(y_k) - g_k)^T (x_{k+1} - y_k) + \left[\tfrac{L}{2} - \tfrac{1}{\alpha_k}\right] \|x_{k+1} - y_k\|^2 \\
        &\quad + \tfrac{1}{\alpha_k}(y_k - x_{k+1} )^T(y_k - \theta_{k+1}x^* - (1 - \theta_{k+1})x_k) \\
        &= (1 - \theta_{k+1})(\phi(x_k) - \phi^*) + \theta_{k+1}(g_k - \nabla f(y_k))^T (x^* -  v_k) + (\nabla f(y_k) - g_k)^T (x_{k+1} - y_k) \\
        &\quad + \theta_{k+1}^2 \left[\tfrac{L}{2} - \tfrac{1}{\alpha_k}\right] \|v_{k+1} - v_k\|^2 + \tfrac{\theta_{k+1}^2}{\alpha_k}(v_k - v_{k+1} )^T(v_k - x^*),
    \end{align*}
    where the equality follows from \eqref{eq:y_k_alternate_convex} and \eqref{eq:v_k_alternate_convex}. 
    Applying \cref{lem:a_b_Appendix} with $a_1 = v_k - v_{k+1}$, $a_2 = v_k - x^*$ and $c = \tfrac{\alpha_k L}{2}$ completes the proof for \textbf{Option II}.
\eproof


We now present the convergence of \cref{alg:stochastic_proximal_gradient} under \cref{ass:convexity} when using a biased gradient estimate, followed by the complexity of number of proximal operator evaluations.

\begin{theorem} \label{th:convergenve_convex}
    Suppose \cref{ass:convexity} holds and the gradient estimate $g_k$ satisfies \cref{cond:sampling}.
    \begin{enumerate}
        \item For \cref{alg:stochastic_proximal_gradient} with \textbf{Option I}:
        \begin{enumerate}
            \item For the finite-sum problem \eqref{eq:deter_error_obj}, if the parameters in \cref{cond:sampling} are chosen such that $\{\eta_k\} \searrow 0$, $\{\eta_k\} \leq \eta < \tfrac{1}{2}$, $\sum_{k=0}^{\infty} \eta_k^2 < \infty$, $\iota_0^2 \in \left[0, \tfrac{1}{2} - \eta\right)$, $\sum_{k=0}^{\infty} \delta_k < \infty$ and $\sum_{k=0}^{\infty} \delta_k^2 < \infty$, and the step size is chosen such that 
            $\{\alpha_k\} = \alpha \leq \tfrac{1 - 2(\eta_k + \iota_0^2)}{2L}$ $\forall k \geq 0$, then $\{\phi(x_k)\}$ converges to the optimal value with $\min_{k = 0, \dots, K-1} \phi(x_{k}) - \phi^* = \Ocal\left(\tfrac{1}{K}\right)$  $\forall K \geq 1$.

            \item For the expectation problem \eqref{eq:deter_error_obj}, if the parameters in \cref{cond:sampling} are chosen such that $\{\tilde{\eta}_k\} \searrow 0$, $\{\tilde{\eta}_k\} \leq \tilde{\eta} < \tfrac{1}{2}$, $\sum_{k=0}^{\infty} \tilde{\eta}_k^2 < \infty$, $\tilde{\iota}_0^2 \in \left[0, \tfrac{1}{2} - \tilde{\eta}\right)$, $\sum_{k=0}^{\infty} \tilde{\delta}_k < \infty$ and $\sum_{k=0}^{\infty} \tilde{\delta}_k^2 < \infty$, and the step size is chosen such that 
            $\{\alpha_k\} = \alpha \leq \tfrac{1 - 2(\tilde{\eta}_k + \tilde{\iota}_0^2)}{2L}$ $\forall k \geq 0$, then $\{\phi(x_k)\}$ converges to the optimal value in expectation with $\min_{k = 0, \dots, K-1} \Embb\left[\phi(x_{k}) - \phi^*\right] = \Ocal\left(\tfrac{1}{K}\right)$ $\forall K \geq 1$.
        \end{enumerate}
            
        \item For \cref{alg:stochastic_proximal_gradient} with \textbf{Option II}:
        \begin{enumerate}
            \item For the finite-sum problem \eqref{eq:deter_error_obj}, if the parameters in \cref{cond:sampling} are chosen such that $\eta_k = \hat{\eta}t_k$ and $\delta_k = \hat{\delta} u_k$ $\forall k \geq 0$, where $\{\eta_k\} \leq \eta < \tfrac{1}{2}$, $\{t_k\} \searrow 0$, $\sum_{k=0}^{\infty} t_k^2 < \infty$, $\sum_{k=0}^{\infty} k t_k^2 < \infty$, $\sum_{k=0}^{\infty} (k+2)^2 u_k < \infty$, $\sum_{k=0}^{\infty} (k + 2)^2 u_k^2 < \infty$, 
            $\iota_0^2 \in \left[0, \tfrac{1}{2} - \eta\right)$ and $\hat{\eta}, \hat{\delta} \geq 0$ are sufficiently small, and the step size is chosen such that 
            $\{\alpha_k\} = \alpha \leq \tfrac{1 - 2(\eta_k + \iota_0^2)}{2L}$ $\forall k \geq 0$, then $\{\phi(x_k)\}$ converges to the optimal value with $\phi(x_{k}) - \phi^* = \Ocal\left(\tfrac{1}{(k + 1)^2}\right)$ $\forall k \geq 0$.

            \item For the expectation problem \eqref{eq:stoch_error_obj}, if the parameters in \cref{cond:sampling} are chosen such that $\tilde{\eta}_k = \hat{\eta}\tilde{t}_k$ and $\tilde{\delta}_k = \hat{\delta} \tilde{u}_k$ $\forall k \geq 0$, where $\{\tilde{\eta}_k\} \leq \tilde{\eta} < \tfrac{1}{2}$, $\{\tilde{t}_k\} \searrow 0$, $\sum_{k=0}^{\infty} \tilde{t}_k^2 < \infty$, $\sum_{k=0}^{\infty} k \tilde{t}_k^2 < \infty$, $\sum_{k=0}^{\infty} (k+2)^2 \tilde{u}_k < \infty$, $\sum_{k=0}^{\infty} (k + 2)^2 \tilde{u}_k^2 < \infty$, 
            $\tilde{\iota}_0^2 \in \left[0, \tfrac{1}{2} - \tilde{\eta}\right)$ and $\hat{\eta}, \hat{\delta} \geq 0$ are sufficiently small, and
            the step size is chosen such that 
            $\{\alpha_k\} = \alpha \leq \tfrac{1 - 2(\tilde{\eta}_k + \tilde{\iota}_0^2)}{2L}$ $\forall k \geq 0$, then $\{\phi(x_k)\}$ converges to the optimal value in expectation with $\Embb\left[\phi(x_{k}) - \phi^*\right] = \Ocal\left(\tfrac{1}{(k + 1)^2}\right)$ $\forall k \geq 0$.
        \end{enumerate}
    \end{enumerate}
\end{theorem}
\bproof
For \cref{alg:stochastic_proximal_gradient} with \textbf{Option I}, applying Cauchy-Schwarz inequality to the result in \cref{lem:general_descent_convex_refined}, we get, 
\begin{align*}
    \phi(x_{k+1})-\phi^*
    &\leq \|g_k - \nabla f(x_k)\| \|x_k - x^*\| + (\nabla f(x_k) - g_k)^T (x_{k+1} - x_k) \numberthis \label{eq:conv_biased_opt1_1}\\
    &\quad +  \tfrac{1}{2\alpha_k}\left[ \|x_k - x^*\|^2 - \|x_{k+1} - x^*\|^2 - \left(1 - \alpha_k L\right)\|x_k - x_{k+1}\|^2\right].
\end{align*}

For the finite-sum problem \eqref{eq:deter_error_obj} for \cref{alg:stochastic_proximal_gradient} with \textbf{Option I}, substituting \cref{cond:sampling} and \eqref{eq:deter_error_prod} from \cref{lem:gradient_error_bounds} into \eqref{eq:conv_biased_opt1_1}, 
\begin{align*}
    \phi(x_{k+1})-\phi^*
    &\leq \left(\tfrac{\eta_k}{2\alpha_k} \left\|x_k - x_{k+1}\right\| + \iota_0 \delta_k\right) \|x_k - x^*\| + \left(\tfrac{\eta_k}{2\alpha_k} \|x_k - x_{k+1}\|^2 + \iota_0\delta_k \|x_k - x_{k+1}\|\right)\\
    &\quad + \tfrac{1}{2\alpha_k}\left[ \|x_k - x^*\|^2 - \|x_{k+1} - x^*\|^2 - \left(1 - \alpha_k L\right)\|x_k - x_{k+1}\|^2\right] \\
    &\leq \tfrac{\eta_k^2}{4\alpha_k} \|x_{k} - x^*\|^2 + \tfrac{1}{4\alpha_k} \left\|x_k - x_{k+1}\right\|^2 + \iota_0 \delta_k \|x_{k} - x^*\| + \tfrac{\alpha_k\delta_k^2}{2}  + \tfrac{\iota_0^2 \|x_k - x_{k+1}\|^2}{2\alpha_k} \\
    &\quad + \tfrac{1}{2\alpha_k}\left[ \|x_k - x^*\|^2 - \|x_{k+1} - x^*\|^2 - \left(1 - \alpha_k L - \eta_k\right)\|x_k - x_{k+1}\|^2\right] \\
    &= \tfrac{\eta_k^2}{4\alpha_k} \|x_k - x^*\|^2 + \iota_0 \delta_k \|x_{k} - x^*\| + \tfrac{\alpha_k\delta_k^2}{2}\\
    &\quad + \tfrac{1}{2\alpha_k}\left[ \|x_k - x^*\|^2 - \|x_{k+1} - x^*\|^2 - \left(\tfrac{1}{2} - \alpha_k L - \eta_k - \iota_0^2\right)\|x_k - x_{k+1}\|^2\right],
\end{align*}
where the second inequality follows from the identity $ab \leq \tfrac{a^2 + b^2}{2}$ and young's inequality as $\iota_0\delta_k \|x_k - x_{k+1}\| \leq \tfrac{\alpha_k\delta_k^2}{2}  + \tfrac{\iota_0^2 \|x_k - x_{k+1}\|^2}{2\alpha_k}$
with $\alpha_k > 0$. Under the defined parameters, the constant $\tfrac{1}{2} - \alpha_k L - \eta_k - \iota_0^2 \geq 0$, and hence the last term in the upper bound can be omitted. Applying young's inequality again, $\|x_{k} - x^*\| \leq \tfrac{\|x_{k} - x^*\|^2}{2\alpha_k} + \tfrac{\alpha_k}{2}$, we get,
\begin{align*}
    \phi(x_{k+1}) - \phi^* 
    &\leq \left(\tfrac{\eta_k^2}{4\alpha_k} + \tfrac{\iota_0\delta_k}{2\alpha_k}\right) \|x_k - x^*\|^2 + \tfrac{\alpha_k(\iota_0\delta_k + \delta_k^2)}{2} + \tfrac{1}{2\alpha_k}\left[\|x_k - x^* \|^2 - \|x_{k+1} - x^*\|^2 \right]\\
    &\leq \left(\tfrac{\eta_k^2}{4\alpha} + \tfrac{\iota_0\delta_k}{2\alpha}\right) \|x_k - x^*\|^2 + \tfrac{(\iota_0\delta_k + \delta_k^2)}{2L} + \tfrac{1}{2\alpha}\left[\|x_k - x^* \|^2 - \|x_{k+1} - x^*\|^2 \right] \numberthis \label{eq:convex_biased_op1_finite_recursion},
\end{align*}
where the second inequality follows from $\alpha = \{\alpha_k\} \leq \tfrac{1}{L}$.
In the above bound, since $\phi(x_{k+1}) - \phi^* \geq 0$, 
one can create a recursive relation,
\begin{align*}
    \|x_{k+1} - x^*\|^2
    &\leq \left(\tfrac{\eta_k^2}{2} + \iota_0\delta_k + 1\right) \|x_k - x^*\|^2 + \tfrac{\alpha}{L}\left(\iota_0\delta_k + \delta_k^2\right).
\end{align*}
From \cref{lem:bounded_sequence_appendix_nonacc}, with $T_k = \| x_{k} - x^*\|^2$, $a_k = \tfrac{\eta_k^2}{2} + \iota_0\delta_k$ and $s_k = \tfrac{\alpha}{L}\left(\iota_0\delta_k + \delta_k^2\right)$, under the defined parameters, the sequence $\left\{\|x_k - x^*\|^2\right\}$ is bounded.
Taking a telescopic sum of \eqref{eq:convex_biased_op1_finite_recursion} for $k = 0, \ldots, K - 1$ yields,
\begin{align*}
    \sum_{k = 0}^{K - 1}\phi(x_{k+1}) - \phi^* 
    &\leq \sum_{k = 0}^{K - 1}\left(\tfrac{\eta_k^2}{4\alpha} + \tfrac{\iota_0\delta_k}{2\alpha}\right) \|x_k - x^*\|^2 + \sum_{k = 0}^{K - 1}\tfrac{(\iota_0\delta_k + \delta_k^2)}{2L} + \tfrac{1}{2\alpha} \|x_0 - x^* \|^2.
\end{align*}
Rearranging the terms in the above inequality,
\begin{align*}
    \min_{k = 0, \dots, K - 1}\phi(x_{k+1}) - \phi^* 
    &\leq \tfrac{1}{K}\left\{\sum_{k = 0}^{\infty}\left(\tfrac{\eta_k^2}{4\alpha} + \tfrac{\iota_0\delta_k}{2\alpha}\right) \|x_k - x^*\|^2 + \sum_{k = 0}^{\infty}\tfrac{(\iota_0\delta_k + \delta_k^2)}{2L} + \tfrac{1}{2\alpha}\|x_0 - x^* \|^2\right\},
\end{align*}
where all terms within the curly brackets on the right-hand side are bounded, completing the proof for the finite-sum problem \eqref{eq:deter_error_obj} for \cref{alg:stochastic_proximal_gradient} with \textbf{Option I}.

For the expectation problem \eqref{eq:stoch_error_obj} for \cref{alg:stochastic_proximal_gradient} with \textbf{Option I}, consider the conditional expectation of \eqref{eq:conv_biased_opt1_1} given $\Gcal_k$. Substituting \eqref{eq:stoch_error_alternate} and \eqref{eq:stoch_error_prod} from \cref{lem:gradient_error_bounds} yields,
\begin{align*}
    \Embb_k\left[\phi(x_{k+1})-\phi^*\right]
    &\leq \left(\tfrac{\tilde{\eta}_k}{2\alpha_k} \sqrt{\Embb_k\left[ \|x_{k} - x_{k+1}\|^2\right]} + \tilde{\iota}_0\tilde{\delta}_k\right) \|x_{k} - x^*\| \\
    &\quad + \left(\tfrac{ \tilde{\eta}_k}{2\alpha_k} \Embb_k\left[\|x_k - x_{k+1}\|^2\right] + \tilde{\iota}_0\tilde{\delta}_k \sqrt{\Embb_k\left[\|x_k - x_{k+1}\|^2\right]}\right) \\
    &\quad +  \tfrac{1}{2\alpha_k}\Embb_k\left[ \|x_k - x^*\|^2 - \|x_{k+1} - x^*\|^2 - \left(1 - \alpha_k L\right)\|x_k - x_{k+1}\|^2\right] \\
    &\leq \tfrac{\tilde{\eta}_k^2}{4\alpha_k} \|x_{k} - x^*\|^2 + \tilde{\iota}_0\tilde{\delta}_k \|x_{k} - x^*\| + \tfrac{\alpha_k\tilde{\delta}_k^2}{2}  \\
    &\quad +  \tfrac{1}{2\alpha_k}\Embb_k\left[ \|x_k - x^*\|^2 - \|x_{k+1} - x^*\|^2 - \left(\tfrac{1}{2} - \alpha_k L - \tilde{\eta}_k - \tilde{\iota}_0^2 \right)\|x_k - x_{k+1}\|^2\right],
\end{align*}
where the second inequality follows from young's inequality. 
Under the defined parameters, the constant $\tfrac{1}{2} - \alpha_k L - \tilde{\eta}_k - \tilde{\iota}_0^2 \geq 0$, and hence the last term in the upper bound can be omitted.
Taking the total expectation of the reduced bound and applying young's inequality yields,
\begin{align*}
    \Embb[\phi(x_{k+1}) - \phi^*] 
    &\leq \left(\tfrac{\tilde{\eta}_k^2}{4\alpha_k} + \tfrac{\tilde{\iota}_0\tilde{\delta}_k}{2\alpha_k}\right) \Embb\left[\|x_{k} - x^*\|^2\right] + \tfrac{\alpha_k (\tilde{\iota}_0\tilde{\delta}_k + \tilde{\delta}_k^2)}{2} + \tfrac{1}{2\alpha_k}\Embb\left[\|x_k - x^*\|^2 - \|x_{k+1} - x^*\|\right]\\
    &\leq \left(\tfrac{\tilde{\eta}_k^2}{\alpha} + \tfrac{\tilde{\iota}_0\tilde{\delta}_k}{2\alpha}\right) \Embb\left[\|x_{k} - x^*\|^2\right] + \tfrac{(\tilde{\iota}_0\tilde{\delta}_k + \tilde{\delta}_k^2)}{2L} + \tfrac{1}{2\alpha}\Embb\left[\|x_k - x^*\|^2 - \|x_{k+1} - x^*\|\right],
\end{align*}
where the second inequality follows from $\alpha = \{\alpha_k\} \leq \tfrac{1}{L}$. 
Similar to the finite-sum problem, one can create a recursive relation and use \cref{lem:bounded_sequence_appendix_nonacc} to show that $\left\{\Embb\left[\|x_{k} - x^*\|^2\right]\right\}$ is a bounded sequence. Taking a telescopic sum of the bound for $k = 0, \ldots, K-1$ yields,
\begin{align*}
    \sum_{k = 0}^{K-1}\Embb[\phi(x_{k+1}) - \phi^*] 
    &\leq \sum_{k = 0}^{K-1}\left(\tfrac{\tilde{\eta}_k^2}{4\alpha} + \tfrac{\tilde{\iota}_0\tilde{\delta}_k}{2\alpha}\right) \Embb\left[\|x_{k} - x^*\|^2\right] + \sum_{k = 0}^{K-1}\tfrac{(\tilde{\iota}_0\tilde{\delta}_k + \tilde{\delta}_k^2)}{2L} + \tfrac{1}{2\alpha}\Embb\left[\|x_0 - x^*\|^2\right].
\end{align*}
Rearranging the terms in the above inequality,
\begin{align*}
    \min_{k = 0, \dots, K - 1}\Embb[\phi(x_{k+1}) - \phi^*] 
    &\leq \tfrac{1}{K}\left\{\sum_{k = 0}^{\infty}\left(\tfrac{\tilde{\eta}_k^2}{4\alpha} + \tfrac{\tilde{\iota}_0\tilde{\delta}_k}{2\alpha}\right) \Embb\left[\|x_{k} - x^*\|^2\right] + \sum_{k = 0}^{\infty}\tfrac{(\tilde{\iota}_0\tilde{\delta}_k + \tilde{\delta}_k^2)}{2L} + \tfrac{1}{2\alpha}\Embb\left[\|x_0 - x^*\|^2\right]\right\},
\end{align*}
where all terms within the curly brackets on the right-hand side are bounded, completing the proof for the expectation problem \eqref{eq:stoch_error_obj} for \cref{alg:stochastic_proximal_gradient} with \textbf{Option I}.

For \cref{alg:stochastic_proximal_gradient} with \textbf{Option II}, applying Cauchy-Schwartz inequality to the result in \cref{lem:general_descent_convex_refined} yields,
\begin{align*}
    \phi(x_{k+1}) - \phi^* 
    &\leq (1 - \theta_{k+1})(\phi(x_k) - \phi^*) + \theta_{k+1}\|g_k - \nabla f(y_k)\|\|v_k - x^*\| \numberthis \label{eq:cnv_biased_opt2} \\
    &\quad + (\nabla f(y_k) - g_k)^T (x_{k+1} - y_k) - \tfrac{\theta_{k+1}^2}{2\alpha_k}\left(1 - \alpha_k L\right)\|v_k - v_{k+1}\|^2 \\
    &\quad + \tfrac{\theta_{k+1}^2}{2\alpha_k}\left[\|v_k - x^*\|^2 - \|v_{k+1} - x^*\|^2\right].
\end{align*}

For the finite-sum problem \eqref{eq:deter_error_obj} for \cref{alg:stochastic_proximal_gradient} with \textbf{Option II}, substituting \cref{cond:sampling}, \eqref{eq:deter_error_prod} from \cref{lem:gradient_error_bounds}, and \eqref{eq:v_k_alternate_convex} into \eqref{eq:conv_biased_opt1_1} yields,
\begin{align*}
    \phi(x_{k+1}) - \phi^* 
    &\leq (1 - \theta_{k+1})(\phi(x_k) - \phi^*) + \theta_{k+1}\left(\tfrac{\theta_{k+1}\eta_k}{2\alpha_k} \left\|v_k - v_{k+1}\right\| + \iota_0 \delta_k\right)\|v_k - x^*\| \\
    &\quad + \left(\tfrac{\eta_k\theta_{k+1}^2}{2\alpha_k} \|v_k - v_{k+1}\|^2 +  \theta_{k+1}\iota_0\delta_k \|v_k - v_{k+1}\|\right) - \tfrac{\theta_{k+1}^2}{2\alpha_k}\left(1 - \alpha_k L\right)\|v_k - v_{k+1}\|^2 \\
    &\quad + \tfrac{\theta_{k+1}^2}{2\alpha_k}\left[\|v_k - x^*\|^2 - \|v_{k+1} - x^*\|^2\right] \\
    &\leq (1 - \theta_{k+1})(\phi(x_k) - \phi^*) + \theta_{k+1}\iota_0 \delta_k\|v_k - x^*\| + \tfrac{\theta_{k+1}^2\eta_k^2}{4\alpha_k} \|v_k - x^*\|^2 + \tfrac{\alpha_k\delta_k^2}{2}\\
    &\quad  - \tfrac{\theta_{k+1}^2}{2\alpha_k}\left(\tfrac{1}{2} - \alpha_k L - \eta_k - \iota_0^2\right)\|v_k - v_{k+1}\|^2 + \tfrac{\theta_{k+1}^2}{2\alpha_k}\left[\|v_k - x^*\|^2 - \|v_{k+1} - x^*\|^2\right],
\end{align*}
where the second inequality follows from young's inequality. 
Under the defined parameters, the constant $\tfrac{1}{2} - \alpha_k L - \eta_k - \iota_0^2 \geq 0$, and hence the fifth term in the upper bound can be ignored. Applying young's inequality to the reduced bound yields,
\begin{align*}
    \phi(x_{k+1}) - \phi^* 
    &\leq (1 - \theta_{k+1})(\phi(x_{k}) - \phi^*) + \tfrac{\theta_{k+1}^2}{2\alpha_k}\left[\iota_0\delta_k + \tfrac{\eta_k^2}{2}\right]\|v_k - x^*\|^2 + \tfrac{\alpha_k (\delta_k^2 + \iota_0\delta_k)}{2}\\
    &\quad + \tfrac{\theta_{k+1}^2}{2\alpha_k}\left[ \|v_{k} - x^*\|^2 - \|v_{k+1} - x^*\|^2\right].
\end{align*}
Dividing the above inequality by $\theta_{k+1}^2$ and applying \eqref{eq:acc_theta_identity_convex}, we get,
\begin{align*}
    \tfrac{\phi(x_{k+1}) - \phi^*}{\theta_{k+1}^2} 
    &\leq \tfrac{1}{\theta_{k}^2}(\phi(x_{k}) - \phi^*) + \tfrac{1}{2\alpha_k}\left[\iota_0\delta_k + \tfrac{\eta_k^2}{2}\right]\|v_k - x^*\|^2 + \tfrac{\alpha_k (\delta_k^2 + \iota_0\delta_k)}{2\theta_{k+1}^2}\\
    &\quad + \tfrac{1}{2\alpha_k}\left[ \|v_{k} - x^*\|^2 - \|v_{k+1} - x^*\|^2\right] \\
    &\leq \tfrac{1}{\theta_{k}^2}(\phi(x_{k}) - \phi^*) + \tfrac{1}{2\alpha}\left[\iota_0\delta_k + \tfrac{\eta_k^2}{2}\right]\|v_k - x^*\|^2 + \tfrac{\delta_k^2 + \iota_0\delta_k}{2\theta_{k+1}^2L}\\
    &\quad + \tfrac{1}{2\alpha}\left[ \|v_{k} - x^*\|^2 - \|v_{k+1} - x^*\|^2\right],
\end{align*}
where the second inequality follows from $\alpha = \{\alpha_k\} \leq \tfrac{1}{L}$. 
By \cref{lem:bounded_sequence_appendix_acc}, with $R_k = \tfrac{1}{\theta_{k}^2}(\phi(x_{k}) - \phi^*)$, $T_k = \tfrac{1}{2\alpha}\|v_{k} - x^*\|^2$, $a_k = \iota_0\delta_k + \tfrac{\eta_k^2}{2}$, and $s_k = \tfrac{\delta_k^2 + \iota_0\delta_k}{2\theta_{k+1}^2L}$, the sequence $\{R_k\}$ is 
bounded 
under the specified parameters 
provided that $\hat{\eta}$ and $\hat{\delta}$ are chosen to be sufficiently small.
Thus, $\exists C \geq 0$ such that, 
$\phi(x_k) - \phi^* \leq C \theta_k^2 = \tfrac{4C}{(k+1)^2}$,
completing the proof for the finite-sum problem \eqref{eq:deter_error_obj} for \cref{alg:stochastic_proximal_gradient} with \textbf{Option II}.

For the expectation problem \eqref{eq:stoch_error_obj} for \cref{alg:stochastic_proximal_gradient} with \textbf{Option II}, consider the conditional expectation of \eqref{eq:cnv_biased_opt2} given $\Gcal_k$. Substituting \eqref{eq:stoch_error_alternate} and \eqref{eq:stoch_error_prod} from \cref{lem:gradient_error_bounds}, and \eqref{eq:v_k_alternate_convex} yields,
\begin{align*}
    \Embb_k\left[\phi(x_{k+1}) - \phi^* \right]
    &\leq (1 - \theta_{k+1})(\phi(x_k) - \phi^*) + \theta_{k+1}\left(\tfrac{\theta_{k+1}\tilde{\eta}_k}{2\alpha_k} \sqrt{\Embb_k\left[ \|v_k - v_{k+1}\|^2\right]} + \tilde{\iota}_0\tilde{\delta}_k\right)\|v_k - x^*\| \\
    &\quad + \tfrac{\theta_{k+1}^2 \tilde{\eta}_k}{2\alpha_k} \Embb_k\left[\|v_k - v_{k+1}\|^2\right] + \theta_{k+1} \tilde{\iota}_0\tilde{\delta}_k \sqrt{\Embb_k\left[\|v_k - v_{k+1}\|^2\right]} \\
    &\quad - \tfrac{\theta_{k+1}^2}{2\alpha_k}\left(1 - \alpha_k L\right)\Embb_k\left[\|v_k - v_{k+1}\|^2\right]  + \tfrac{\theta_{k+1}^2}{2\alpha_k}\Embb_k\left[\|v_k - x^*\|^2 - \|v_{k+1} - x^*\|^2\right] \\
    &\leq (1 - \theta_{k+1})(\phi(x_k) - \phi^*) + \theta_{k+1}\tilde{\iota}_0\tilde{\delta}_k\|v_k - x^*\|  + \tfrac{\theta_{k+1}^2\tilde{\eta}_k^2}{4\alpha_k} \|v_k - x^*\|^2 + \tfrac{\alpha_k\tilde{\delta}^2_k}{2} \\
    &\quad - \tfrac{\theta_{k+1}^2}{2\alpha_k}\left(\tfrac{1}{2} - \alpha_k L - \tilde{\eta}_k - \tilde{\iota}_0^2\right)\Embb_k\left[\|v_k - v_{k+1}\|^2\right]  + \tfrac{\theta_{k+1}^2}{2\alpha_k}\Embb_k\left[\|v_k - x^*\|^2 - \|v_{k+1} - x^*\|^2\right],
\end{align*}
where the second inequality follows from young's inequality. 
Under the defined parameters, the constant $\tfrac{1}{2} - \alpha_k L - \tilde{\eta}_k - \tilde{\iota}_0^2 \geq 0$, and hence the fifth term in the upper bound can be ignored. Applying young's inequality to the reduced bound yields,
\begin{align*}
    \Embb_k\left[\phi(x_{k+1}) - \phi^*\right] 
    &\leq (1 - \theta_{k+1})(\phi(x_{k}) - \phi^*) + \tfrac{\theta_{k+1}^2}{2\alpha_k}\left[\tilde{\iota}_0\tilde{\delta}_k + \tfrac{\tilde{\eta}_k^2}{2}\right]\|v_k - x^*\|^2 + \tfrac{\alpha_k (\tilde{\delta}_k^2 + \tilde{\iota}_0\tilde{\delta}_k)}{2}\\
    &\quad + \tfrac{\theta_{k+1}^2}{2\alpha_k}\Embb_k\left[ \|v_{k} - x^*\|^2 - \|v_{k+1} - x^*\|^2\right].
\end{align*}
Taking the total expectation of the above bound, dividing by $\theta_{k+1}^2$ and using \eqref{eq:acc_theta_identity_convex}, we get,
\begin{align*}
    \Embb\left[\tfrac{\phi(x_{k+1}) - \phi^*}{\theta_{k+1}^2}\right] 
    &\leq \tfrac{1}{\theta_{k}^2}\Embb\left[\phi(x_{k}) - \phi^*\right] + \tfrac{1}{2\alpha_k}\left[\tilde{\iota}_0\tilde{\delta}_k + \tfrac{\tilde{\eta}_k^2}{2}\right]\Embb\left[\|v_k - x^*\|^2\right] + \tfrac{\alpha_k (\tilde{\delta}_k^2 + \tilde{\iota}_0\tilde{\delta}_k)}{2\theta_{k+1}^2}\\
    &\quad + \tfrac{1}{2\alpha_k}\Embb\left[\|v_{k} - x^*\|^2 - \|v_{k+1} - x^*\|^2\right] \\
    &\leq \tfrac{1}{\theta_{k}^2}\Embb\left[\phi(x_{k}) - \phi^*\right] + \tfrac{1}{2\alpha}\left[\tilde{\iota}_0\tilde{\delta}_k + \tfrac{\tilde{\eta}_k^2}{2}\right]\Embb\left[\|v_k - x^*\|^2\right] + \tfrac{\tilde{\delta}_k^2 + \tilde{\iota}_0\tilde{\delta}_k}{2\theta_{k+1}^2 L}\\
    &\quad + \tfrac{1}{2\alpha}\Embb\left[\|v_{k} - x^*\|^2 - \|v_{k+1} - x^*\|^2\right],
\end{align*}
where the final inequality results from $\alpha = \{\alpha_k\} \leq \tfrac{1}{L}$. 
By \cref{lem:bounded_sequence_appendix_acc}, with $R_k = \tfrac{1}{\theta_{k}^2}\Embb\left[\phi(x_{k}) - \phi^*\right]$, $T_k = \tfrac{1}{2\alpha}\Embb\left[\|v_{k} - x^*\|^2\right]$, $a_k = \tilde{\iota}_0\tilde{\delta}_k + \tfrac{\tilde{\eta}_k^2}{2}$, and $s_k = \tfrac{\tilde{\delta}_k^2 + \tilde{\iota}_0\tilde{\delta}_k}{2\theta_{k+1}^2L}$, the sequence $\{R_k\}$ is bounded under the specified parameters provided that $\hat{\eta}$ and $\hat{\delta}$ are chosen to be sufficiently small. Thus, $\exists C \geq 0$ such that, $\Embb\left[\phi(x_k) - \phi^*\right] \leq C \theta_k^2 = \tfrac{4C}{(k+1)^2}$,
completing the proof for the expectation problem \eqref{eq:stoch_error_obj} for \cref{alg:stochastic_proximal_gradient} with \textbf{Option II}.
\eproof

\cref{th:convergenve_convex} establishes a sublinear rate of convergence for \cref{alg:stochastic_proximal_gradient} with \textbf{Option I} and \textbf{Option II} for a general convex objective function when using a biased gradient estimate satisfying \cref{cond:sampling}. 
Thus, an $\epsilon > 0$ accurate solution, i.e., $\phi(x_k) - \phi^* \leq \epsilon$ for the finite-sum problem \eqref{eq:deter_error_obj} and $\Embb\left[\phi(x_k) - \phi^* \right] \leq \epsilon$ for the expectation problem \eqref{eq:stoch_error_obj}, which are stronger definitions than those defined in \cref{sec:non_convexity}, can be achieved in $\Ocal\left(\tfrac{1}{\epsilon}\right)$ iterations (proximal operator evaluations) with \textbf{Option I} and in $\Ocal\left(\tfrac{1}{\sqrt{\epsilon}}\right)$ iterations (proximal operator evaluations) with \textbf{Option II}, matching the results for deterministic first-order methods \cite{nesterov2018lectures}.

When gradient estimate is unbiased, the required conditions on the parameters can be simplified for the expectation problem \eqref{eq:stoch_error_obj} as shown in the following corollary.

\begin{corollary} \label{corollary:convergence_unbiased_convex}
    Suppose the conditions in \cref{th:convergenve_convex} hold for the expectation problem \eqref{eq:stoch_error_obj} 
    and the gradient estimate is unbiased, i.e. $\Embb_k[g_k] = \nabla f(y_k)$.
    \begin{enumerate}
        \item For \cref{alg:stochastic_proximal_gradient} with \textbf{Option I}, if the parameters in \cref{cond:sampling} are chosen such that $\{\tilde{\eta}_k\} = \tilde{\eta} \in [0, 1)$, $\tilde{\iota}_0^2 \in [0, 1 - \tilde{\eta})$, and $\sum_{k=0}^{\infty}\tilde{\delta}_k^2 < \infty$, 
        and the step size is chosen such that $\{\alpha_k\} = \alpha \leq \tfrac{1 - \tilde{\eta} - \tilde{\iota}_0^2}{L}$, then $\{\phi(x_k)\}$ converges to the optimal value in expectation with $\min_{k=0, \dots, K-1} \Embb[\phi(x_{k}) - \phi^*] = \Ocal\left(\tfrac{1}{K}\right)$ $\forall K \geq 1$.

        \item For \cref{alg:stochastic_proximal_gradient} with \textbf{Option II}, if the parameters in \cref{cond:sampling} are chosen such that $\{\tilde{\eta}_k\} = \tilde{\eta} \in [0, 1) $, $\tilde{\iota}_0^2 \in [0, 1 - \tilde{\eta})$, and $\sum_{k=0}^{\infty}\tfrac{\tilde{\delta}_k^2}{\theta_{k+1}^2} < \infty$, and the step size is chosen such that $\{\alpha_k\} = \alpha \leq \tfrac{1 - \tilde{\eta} - \tilde{\iota}_0^2}{L}$, then $\{\phi(x_k)\}$ converges to the optimal value in expectation with 
        $\Embb[\phi(x_{K}) - \phi^*] = \Ocal\left(\tfrac{1}{K^2}\right)$ $\forall K \geq 1$.
        \end{enumerate}
\end{corollary}
\bproof
For \cref{alg:stochastic_proximal_gradient} with \textbf{Option I}, consider the conditional expectation of the result in \cref{lem:general_descent_convex_refined} given $\Gcal_k$ under the defined parameters.
From the assumption of an unbiased gradient estimate and substituting \eqref{eq:stoch_error_prod_squared} from \cref{lem:gradient_error_bounds}, the bound reduces to, 
\begin{align*}
    \Embb_k \left[\phi(x_{k+1})-\phi^*\right]
    &\leq \tfrac{\left(\tilde{\eta} + \tilde{\iota}_0^2\right)}{2\alpha} \Embb_k\left[\|x_k - x_{k+1}\|^2\right] + \tfrac{\alpha \tilde{\delta}_k^2}{2} \\
    &\quad + \tfrac{1}{2\alpha}\left[ \Embb_k\left[\|x_k - x^*\|^2\right] - \Embb_k\left[\|x_{k+1} - x^*\|^2\right] - \left(1 - \alpha L\right)\Embb_k\left[\|x_k - x_{k+1}\|^2\right]\right] \\
    &= \tfrac{1}{2\alpha}\left[ \Embb_k\left[\|x_k - x^*\|^2\right] - \Embb_k\left[\|x_{k+1} - x^*\|^2\right]\right]  + \tfrac{\alpha \tilde{\delta}_k^2}{2} \\
    &\quad - \tfrac{1}{2\alpha}\left[1 - \alpha L - \tilde{\eta} - \tilde{\iota}_0^2\right]\Embb_k\left[\|x_k - x_{k+1}\|^2\right],
\end{align*}
where the constant $1 - \alpha L - \tilde{\eta} - \tilde{\iota}_0^2 \geq 0$ under the defined parameters. Taking a telescopic sum for $k = 0, \ldots, K - 1$ of the total expectation of the reduced bound yields,
\begin{align*}
    \sum_{k=0}^{K-1} \Embb \left[\phi(x_{k+1})-\phi^*\right]
    &\leq \tfrac{1}{2\alpha}\left[ \Embb\left[\|x_0 - x^*\|^2\right] - \Embb\left[\|x_{K} - x^*\|^2\right]\right]  + \sum_{k=0}^{K-1} \tfrac{\alpha \tilde{\delta}_k^2}{2}.
\end{align*}
Rearranging the terms of the above inequality, we get
\begin{align*}
    \min_{k=0, \dots, K-1} \Embb \left[\phi(x_{k+1})-\phi^*\right]
    &\leq \tfrac{1}{K} \left\{\tfrac{1}{2\alpha}\Embb\left[\|x_0 - x^*\|^2\right] + \tfrac{\alpha}{2}\sum_{k=0}^{K-1} \tilde{\delta}_k^2 \right\},
\end{align*}
where all terms within the curly brackets on the right-hand side are bounded due to the condition $\sum_{k=0}^{\infty}\tilde{\delta}_k^2 < \infty$, 
completing the proof for \textbf{Option I}.

For \cref{alg:stochastic_proximal_gradient} with \textbf{Option II}, consider the conditional expectation of the result in \cref{lem:general_descent_convex_refined} given $\Gcal_k$ under the defined parameters.
From the assumption of an unbiased gradient estimate and substituting \eqref{eq:stoch_error_prod_squared} from \cref{lem:gradient_error_bounds}, the bound reduces to,
\begin{align*}
    \Embb_k\left[\phi(x_{k+1}) - \phi^*\right] 
    &\leq (1 - \theta_{k+1})(\phi(x_k) - \phi^*) + \tfrac{\left(\tilde{\eta} + \tilde{\iota}_0^2\right)}{2\alpha} \Embb_k\left[\|x_{k+1} - y_k\|^2\right] + \tfrac{\alpha \tilde{\delta}_k^2}{2} \\
    &\quad + \tfrac{\theta_{k+1}^2}{2\alpha}\Embb_k \left[\|v_k - x^*\|^2 - \|v_{k+1} - x^*\|^2 - \left(1 - \alpha L\right)\|v_k - v_{k+1}\|^2\right] \\
    &= (1 - \theta_{k+1})(\phi(x_k) - \phi^*) + \tfrac{\theta_{k+1}^2}{2\alpha}\Embb_k \left[\|v_k - x^*\|^2 - \|v_{k+1} - x^*\|^2 \right] + \tfrac{\alpha \tilde{\delta}_k^2}{2} \\
    &\quad - \tfrac{\theta_{k+1}^2}{2\alpha} \left[1 - \alpha L - \tilde{\eta} - \tilde{\iota}_0^2 \right] \Embb_k\left[\|v_k - v_{k+1}\|^2\right]
\end{align*}
where the equality follows from \eqref{eq:v_k_alternate_convex} and the constant $1 - \alpha L - \tilde{\eta} - \tilde{\iota}_0^2 \geq 0$ under the defined parameters. 
Dividing the total expectation of the reduced bound by $\theta_{k+1}^2$ and applying \eqref{eq:acc_theta_identity_convex} yields,
\begin{align*}
    \tfrac{1}{\theta_{k+1}^2}\Embb\left[\phi(x_{k+1}) - \phi^*\right] 
    &\leq \tfrac{1}{\theta_k^2}\Embb\left[\phi(x_k) - \phi^*\right] + \tfrac{1}{2\alpha}\Embb \left[\|v_k - x^*\|^2 - \|v_{k+1} - x^*\|^2 \right] + \tfrac{\alpha \tilde{\delta}_k^2}{2\theta_{k+1}^2}.
\end{align*}
Taking a telescoping sum of the above for $k = 0, \dots, K - 1$ yields,
\begin{align*}
    \tfrac{1}{\theta_{K}^2}\Embb\left[\phi(x_{K}) - \phi^*\right] 
    &\leq \tfrac{1}{\theta_0^2}\left[\phi(x_0) - \phi^*\right] + \tfrac{1}{2\alpha}\Embb \left[\|v_0 - x^*\|^2 - \|v_{K} - x^*\|^2 \right] + \sum_{k=0}^{K-1} \tfrac{\alpha \tilde{\delta}_k^2}{2\theta_{k+1}^2} 
\end{align*}
Multiplying the above inequality by $\theta_{K}^2$ yields,
\begin{align*}
    \Embb\left[\phi(x_{K}) - \phi^*\right] 
    &\leq \tfrac{4}{(K+1)^2} \left\{\tfrac{1}{4}\left[\phi(x_0) - \phi^*\right] + \tfrac{1}{2\alpha}\Embb \left[\|v_0 - x^*\|^2 - \|v_{K} - x^*\|^2 \right] + \sum_{k=0}^{K-1} \tfrac{\alpha \tilde{\delta}_k^2}{2\theta_{k+1}^2} \right\},
\end{align*}
where all terms within the curly brackets on the right-hand side are bounded due to the defined $\tilde{\delta}_k$, 
completing the proof for \textbf{Option II}.
\eproof
Compared to \cref{th:convergenve_convex}, the parameter settings in \cref{corollary:convergence_unbiased_convex} for \cref{alg:stochastic_proximal_gradient} and \cref{cond:sampling} are less restrictive, as the use of an unbiased gradient estimate simplifies the analysis by reducing the number of error terms involved.

We now present the complexity for number of stochastic gradient evaluations for \cref{alg:stochastic_proximal_gradient} when an unbiased gradient estimate is used for the expectation problem \eqref{eq:stoch_error_obj} under \cref{ass:convexity}.

\begin{theorem} \label{th:sample_complexity_convex}
    Suppose Assumptions \ref{ass:bounded_variance} and \ref{ass:convexity} hold, and \cref{cond:sampling} is satisfied for the expectation problem \eqref{eq:stoch_error_obj} via the unbiased gradient estimate in \cref{lem:unbiased_sample_set_bound}. Then, to achieve a solution satisfying 
    $\min\left\{\Embb\left[\phi(x_k) - \phi^*\right], \left\|R_{\alpha_k}^{true}(y_k)\right\|^2\right\} \leq \epsilon$
    with  $\epsilon > 0$, the number of stochastic gradient evaluations required is as follows:
    \begin{enumerate}
        \item For \cref{alg:stochastic_proximal_gradient} with \textbf{Option I}, if the parameters in \cref{cond:sampling} are chosen such that $\{\tilde{\eta}_k\} = \tilde{\eta} \in (0, 1)$, $\tilde{\iota}_0^2 \in [0, 1 - \tilde{\eta})$, and $\tilde{\delta}_k^2  = \tfrac{1}{(k + 1)^{1 + \nu}}$ $\forall k \geq 0$, where $\nu > 0$, and the step size is chosen such that $\{\alpha_k\} = \alpha \leq \tfrac{1 - \tilde{\eta} - \tilde{\iota}_0^2}{L}$, then the number of stochastic gradient evaluations required is $\Ocal \left(\epsilon^{-(2 + \nu)}\right)$. If $\{\tilde{\delta}_k\} = 0$, this improves to $\Ocal\left(\epsilon^{-2}\right)$. 

        \item For \cref{alg:stochastic_proximal_gradient} with \textbf{Option II}, if the parameters in \cref{cond:sampling} are chosen such that $\{\tilde{\eta}_k\} = \tilde{\eta} \in (0, 1) $, $\tilde{\iota}_0^2 \in [0, 1 - \tilde{\eta})$, and $\tilde{\delta}_k^2  = \tfrac{1}{(k + 1)^{3 + 2\nu}}$ $\forall k \geq 0$, where $\nu > 0$, and the step size is chosen such that $\{\alpha_k\} = \alpha \leq \tfrac{1 - \tilde{\eta} - \tilde{\iota}_0^2}{L}$, then the number of stochastic gradient evaluations required is $\Ocal \left(\epsilon^{-(2 + \nu)}\right)$. If $\{\tilde{\delta}_k\} = 0$, this improves to $\Ocal\left(\epsilon^{-\tfrac{3}{2}}\right)$. 
    \end{enumerate}
\end{theorem}
\bproof
Let $K_\epsilon \geq 1$ be the first iteration that achieves the desired solution accuracy. Hence, $\left\|R_{\alpha_k}^{true}(y_k)\right\|^2 > \epsilon$ $\forall k \leq K_\epsilon - 1$ and from \eqref{eq:complexity_bound_help}, 
$\left\|\Embb_k\left[R_{\alpha_k}(y_k)\right]\right\|^2 > \left(1 + \tfrac{\tilde{\eta}^2}{4}\right)^{-1}\left[\tfrac{\epsilon}{2} - \tilde{\iota}_0^2\tilde{\delta}_k^2\right]$ $\forall k \leq K_\epsilon - 1$.
Thus, the total number of stochastic gradient evaluations can be bounded using \cref{lem:unbiased_sample_set_bound} as, 
\begin{align*}
    \sum_{k=0}^{K_\epsilon - 1}|S_k| 
    &= \sum_{k=0}^{K_\epsilon - 1} \left\lceil\tfrac{\sigma^2}{\tfrac{\tilde{\eta}_k^2}{4} \left\|\Embb_k\left[R_{\alpha_k}(y_k)\right]\right\|^2 + \tilde{\iota}_0^2\tilde{\delta}_k^2} \right\rceil \leq \sum_{k=0}^{K_\epsilon - 1} \tfrac{2\sigma^2 (4 + \tilde{\eta}^2)}{\tilde{\eta}^2 \epsilon + 8\tilde{\iota}_0^2\tilde{\delta}_k^2}  + 1 
    \leq \sum_{k=0}^{K_\epsilon - 1} \tfrac{2\sigma^2(4 + \tilde{\eta}^2)}{\tilde{\eta}^2 \epsilon} + \tfrac{\sigma^2 (4 + \tilde{\eta}^2)}{4\tilde{\iota}_0^2\tilde{\delta}_k^2}  + K_{\epsilon}.
\end{align*}

For \cref{alg:stochastic_proximal_gradient} with \textbf{Option I}, $K_\epsilon$ is at most $\Ocal\left(\tfrac{1}{\epsilon}\right)$ from \cref{corollary:convergence_unbiased_convex}, yielding 
\begin{align*}
    \sum_{k=0}^{K_\epsilon-1}|S_k| 
    \leq \tfrac{2\sigma^2 (4 + \tilde{\eta}^2)}{\tilde{\eta}^2 \epsilon}K_{\epsilon}  + \tfrac{\sigma^2 (4 + \tilde{\eta}^2) }{4\tilde{\iota}_0^2}K_{\epsilon}^{2 + \nu}  + K_{\epsilon} = \Ocal\left(\tfrac{1}{\epsilon^{2 + \nu}}\right).
\end{align*}
Following the same procedure, if $\{\tilde{\delta}_k\} = 0$, $\sum_{k=0}^{K_\epsilon-1}|S_k| \leq \tfrac{2\sigma^2 (4 + \tilde{\eta}^2) }{\tilde{\eta}^2 \epsilon}K_{\epsilon}  + K_{\epsilon} = \Ocal\left(\tfrac{1}{\epsilon^{2}}\right)$. 

For \cref{alg:stochastic_proximal_gradient} with \textbf{Option II}, $K_\epsilon$ is at most $\Ocal\left(\tfrac{1}{\sqrt{\epsilon}}\right)$ from \cref{corollary:convergence_unbiased_convex}, yielding
\begin{align*}
    \sum_{k=0}^{K_\epsilon-1}|S_k| 
    \leq \tfrac{2\sigma^2(4 + \tilde{\eta}^2)}{\tilde{\eta}^2 \epsilon}K_\epsilon  + \tfrac{\sigma^2 (4 + \tilde{\eta}^2) }{4\tilde{\iota}_0^2}K_{\epsilon}^{4 + 2\nu} + K_{\epsilon} = \Ocal\left(\tfrac{1}{\epsilon^{2 + \nu}}\right).
\end{align*}
Following the same procedure, if $\{\tilde{\delta}_k\} = 0$, $\sum_{k=0}^{K_\epsilon-1}|S_k| \leq \tfrac{2\sigma^2 (4 + \tilde{\eta}^2)}{\tilde{\eta}^2 \epsilon} K_{\epsilon}  + K_{\epsilon} = \Ocal\left(\tfrac{1}{\epsilon^{3/ 2}}\right)$. 
\eproof

\cref{th:sample_complexity_convex} matches the optimal complexity for the number of stochastic gradient evaluations for the expectation problem \eqref{eq:stoch_error_obj} with a general convex objective function \cite{lan2012optimal}. 
We conclude this section with a corollary to \cref{th:sample_complexity_convex}, similar to \cref{cor:sample_complexity_non_convex}, using a definition of an $\epsilon$-accurate solution similar to that in \cite{lan2012optimal},  under the parameter setting $\{\tilde{\eta}_k\} = 0$.

\begin{corollary} \label{cor:sample_complexity_convex}
    Suppose the conditions in \cref{th:sample_complexity_convex} hold. Then, to achieve a solution satisfying $\Embb\left[\phi(x_k) - \phi^*\right] \leq \epsilon$ with $\epsilon > 0$:
    \begin{enumerate}
        \item For \cref{alg:stochastic_proximal_gradient} with \textbf{Option I}, if the parameters in \cref{cond:sampling}  are chosen as
        $\{\tilde{\eta}_k\} = 0$, $\tilde{\iota}_0^2 \in (0, 1)$, and $\tilde{\delta}_k^2  = \tfrac{1}{(k + 1)^{1 + \nu}}$ $\forall k \geq 0$, where $\nu > 0$, and the step size is chosen such that $\{\alpha_k\} = \alpha \leq \tfrac{1 - \tilde{\iota}_0^2}{L}$, then the number of stochastic gradient evaluations required is $\Ocal \left(\epsilon^{-(2 + \nu)}\right)$.

        \item For \cref{alg:stochastic_proximal_gradient} with \textbf{Option II}, if the parameters in \cref{cond:sampling} are chosen as $\{\tilde{\eta}_k\} = 0 $, $\tilde{\iota}_0^2 \in (0, 1)$, and $\tilde{\delta}_k^2  = \tfrac{1}{(k + 1)^{3 + 2\nu}}$ $\forall k \geq 0$, where $\nu > 0$, and the step size is chosen such that $\{\alpha_k\} = \alpha \leq \tfrac{1 - \tilde{\iota}_0^2}{L}$, then the number of stochastic gradient evaluations required is $\Ocal \left(\epsilon^{-(2 + \nu)}\right)$.
    \end{enumerate}
\end{corollary}
\bproof
The proof follows from the same procedure as \cref{th:sample_complexity_convex}.
\eproof
The conditions in \cref{cor:sample_complexity_convex} reduce \cref{cond:sampling} to using a predetermined error sequence for the gradient estimates, similar to \cite{schmidt2011convergence}.
\subsection{Strongly Convex Objective Function} \label{sec:strong_convexity}

In this section, we present the theoretical analysis of \cref{alg:stochastic_proximal_gradient} when the smooth function $f(x)$, and thus the composite function $\phi(x)$, is strongly convex.
We begin by stating the basic assumptions and definitions, along with some mathematical identities that will be used throughout the analysis.

\begin{assumption} \label{ass:strong_convexity}
The function $f : \Rmbb^d \rightarrow \Rmbb$ is $L$-smooth and $\mu$-strongly convex, i.e., 
    \begin{equation*}
        f(\gamma a + (1 - \gamma) b) \leq \gamma f(a) + (1 - \gamma) f(b) - \tfrac{\mu}{2} \gamma (1 - \gamma) \|a - b\|^2 \quad \forall a, b \in \Rmbb^d, \,\, \forall \gamma \in [0, 1],
    \end{equation*}
    and the function $h : \Rmbb^d \rightarrow \Rmbb \cup \{+\infty\}$ is closed, convex, and proper.
\end{assumption}
Under \cref{ass:strong_convexity}, let $x^* \in \Rmbb^d$ be the unique optimal solution. Since $f(x)$ is differentiable and strongly convex, from \cite{nesterov2018lectures},
\begin{equation} \label{eq:strong_convexity_alternate}
    f(b) \geq f(a) + \nabla f(a)^T (b - a) + \tfrac{\mu}{2} \|b - a\|^2 \quad \forall a, b \in \Rmbb^d.
\end{equation}

For \cref{alg:stochastic_proximal_gradient} with \textbf{Option II}, under \cref{ass:strong_convexity}, we define the sequence $\{\beta_k\}$ and two additional sequences $\{\theta_k\}$ and $\{v_k\}$ for the analysis as,
\begin{equation} \label{eq:acc_defintions_strong_convex}
    \beta_k = \tfrac{1 - \theta_k}{1 + \theta_k} \quad \text{and} \quad \theta_k = \sqrt{\mu\alpha_k} \quad \forall k\geq 0 \text{, and} \quad v_{k} = x_{k-1} + \tfrac{1}{\theta_{k}} (x_{k} - x_{k-1}) \quad \forall k \geq 1,
\end{equation}
with $v_0 = x_0$.
Under these definitions, $y_k$ can be expressed as,
\begin{align*}
    y_{k} 
    &= x_k + \tfrac{1-  \theta_k}{1 + \theta_k}\left(x_k - x_{k-1}\right) = \tfrac{1}{1 + \theta_k}\left(x_k + \theta_k x_{k-1} + x_{k} - x_{k-1}\right)  \\
    &= \tfrac{x_k + \theta_{k}v_{k}}{1 + \theta_{k}} \quad \forall k \geq 0. \numberthis \label{eq:y_k_alternate_strong_convex}
\end{align*}
When a constant step size is employed in \cref{alg:stochastic_proximal_gradient}, i.e., $\{\alpha_k\} = \alpha$ and $\{\theta_k\} = \theta = \sqrt{\mu\alpha}$, then  using \eqref{eq:acc_defintions_strong_convex} and \eqref{eq:y_k_alternate_strong_convex}, the update form for $\{v_k\}$ can be expressed as,
\begin{align*}
    v_{k+1}
    &= x_k + \tfrac{1}{\theta} \left(x_{k+1} - x_k \right) = v_k \left(1 + \tfrac{\theta^2}{1 + \theta} - \theta - \tfrac{1}{1 + \theta}\right)  + x_{k} \left(\tfrac{\theta}{1 + \theta} - \tfrac{1}{\theta(1 + \theta)} \right) + x_{k+1}\left(\tfrac{1}{\theta}\right) \\
    &= v_k  + \theta \left(\tfrac{x_k + \theta v_{k}}{1 + \theta} - v_k\right) - \tfrac{1}{\theta} \left(\tfrac{x_k + \theta v_{k}}{1 + \theta} - x_{k+1}\right) \\
    &= v_k  + \theta (y_k - v_k) - \tfrac{1}{\theta} (y_k - x_{k+1})  \quad \forall k \geq 0. \numberthis \label{eq:v_k_alternate_strong_convex}
\end{align*}

We now establish a descent lemma that further refines \cref{lem:general_descent_lemma} under \cref{ass:strong_convexity}.

\begin{lemma} \label{lem:prelim_descent_strong_convex}
    Suppose \cref{ass:strong_convexity} holds. Then, $\forall x \in \Rmbb^d$ and $\forall k \geq 0$, the iterates generated by \cref{alg:stochastic_proximal_gradient} satisfy,
    \begin{align*}
        \phi(x_{k+1}) 
        &\leq 
        \phi(x)  + \left[\tfrac{1}{2\alpha_k} - \tfrac{\mu}{2}\right] \|x - y_k\|^2  + (g_k - \nabla f(y_k))^T (x - y_k)\\
        &\quad + (\nabla f(y_k) - g_k)^T (x_{k+1} - y_k) - \left(\tfrac{1}{2 \alpha_k} - \tfrac{L}{2}\right) \|x_{k+1} - y_{k}\|^2.
    \end{align*}
\end{lemma}
\bproof
From \eqref{eq:strong_convexity_alternate}, the first two terms in the result from \cref{lem:general_descent_lemma} can be bounded as,
\begin{align*}
    \phi(x_{k+1}) 
    &\leq f(x) - \tfrac{\mu}{2}\|x - y_k\|^2 + (g_k - \nabla f(y_k))^T (x - y_k) + \tfrac{1}{2\alpha_k} \|x - y_k\|^2 + h(x)  \\
    &\quad + (\nabla f(y_k) - g_k)^T (x_{k+1} - y_k) - \left(\tfrac{1}{2 \alpha_k} - \tfrac{L}{2}\right) \|x_{k+1} - y_{k}\|^2,
\end{align*}
completing the proof.
\eproof

While the result in \cref{lem:prelim_descent_strong_convex} is sufficient for analyzing \cref{alg:stochastic_proximal_gradient} with \textbf{Option I}, a different Lyapunov function is required to analyze \cref{alg:stochastic_proximal_gradient} with \textbf{Option II} for a strongly convex objective function.
We first establish a recursion for the distance between $\{v_k\}$ and the optimal solution $x^*$ for \cref{alg:stochastic_proximal_gradient} with \textbf{Option II}.

\begin{lemma} \label{lem:v_alternate_strong_convex}
    Suppose \cref{ass:strong_convexity} holds. Then, $\forall k \geq 0$, the iterates generated by \cref{alg:stochastic_proximal_gradient} with \textbf{Option II} and $\{\alpha_k\} = \alpha \leq \tfrac{1}{L}$ satisfy,
    \begin{align*}
        \tfrac{\mu}{2}\|v_{k + 1} - x^*\|^2 
        &= \tfrac{\mu}{2}(1 - \theta)\|v_k - x^*\|^2 - \tfrac{\mu}{2}\theta(1 - \theta)\|v_k - y_k\|^2
        + \tfrac{\mu\theta}{2} \|y_k - x^*\|^2 \\
        &\quad + \tfrac{1}{2\alpha} \|y_k - x_{k+1}\|^2 
        + \tfrac{1}{\alpha} (y_k - x_{k+1})^T\left[\theta x^* + (1 - \theta)x_k - y_k \right],
    \end{align*}
     where $\theta = \sqrt{\mu\alpha}$.
\end{lemma}
\bproof
From \eqref{eq:v_k_alternate_strong_convex},
\begin{align*}
    v_{k+1} - x^* &= v_k - x^* + \theta (y_k - v_k) - \tfrac{1}{\theta} (y_k - x_{k+1}) = (1 - \theta)(v_k - x^*) + \theta (y_k - x^*) - \tfrac{1}{\theta} (y_k - x_{k+1}).
\end{align*}
Taking the euclidean norm of the above and squaring both sides yields,
\begin{align*}
    \|v_{k+1} - x^*\|^2 &= 
    (1 - \theta)^2\|v_k - x^*\|^2
    + \theta^2 \|y_k - x^*\|^2
    + \tfrac{1}{\theta^2} \|y_k - x_{k+1}\|^2  \\
    &\quad + 2\theta (1 - \theta) (v_k - x^*)^T(y_k - x^*) - \tfrac{2}{\theta} (y_k - x_{k+1})^T\left[(1 - \theta)(v_k - x^*) + \theta (y_k - x^*)\right].
\end{align*}
Multiplying the above equality by $\tfrac{\mu}{2} = \tfrac{\theta^2}{2\alpha}$ and using $(1 - \theta)^2 = 1 - \theta - \theta(1 - \theta)$, we get,
\begin{align*}
    \tfrac{\mu}{2}\|v_{k+1} - x^*\|^2 &= 
    \tfrac{\mu}{2}(1 - \theta)\|v_k - x^*\|^2
    - \tfrac{\mu}{2}\theta(1 - \theta)\|v_k - x^*\|^2
    + \tfrac{\mu\theta^2}{2} \|y_k - x^*\|^2
    + \tfrac{1}{2\alpha} \|y_k - x_{k+1}\|^2  \\
    &\quad + \mu\theta (1 - \theta) (v_k - x^*)^T(y_k - x^*) - \tfrac{\theta}{\alpha} (y_k - x_{k+1})^T\left[(1 - \theta)(v_k - x^*) + \theta (y_k - x^*)\right] \numberthis \label{eq:stong_convex_op2_v_k_x^*_equality}.
\end{align*}
The second, third and fifth terms on the right-hand side of \eqref{eq:stong_convex_op2_v_k_x^*_equality} can be simplified together as,
\begin{align*}
    &- \tfrac{\mu}{2}\theta(1 - \theta)\|v_k - x^*\|^2
    + \tfrac{\mu\theta^2}{2} \|y_k - x^*\|^2
    + \mu\theta (1 - \theta) (v_k - x^*)^T(y_k - x^*) \\
    &\quad = - \tfrac{\mu}{2}\theta(1 - \theta)\|v_k - x^*\|^2
    + \mu\theta (1 - \theta) (v_k - x^*)^T(y_k - x^*)
    - \tfrac{\mu\theta(1 - \theta)}{2} \|y_k - x^*\|^2
    + \tfrac{\mu\theta}{2} \|y_k - x^*\|^2 \\
    &\quad = - \tfrac{\mu}{2}\theta(1 - \theta)\|v_k - y_k\|^2
    + \tfrac{\mu\theta}{2} \|y_k - x^*\|^2,
\end{align*}
where the second equality follows from $\theta^2 = \theta - \theta(1 - \theta)$.
The last term on the right-hand side of \eqref{eq:stong_convex_op2_v_k_x^*_equality} can be simplified as,
\begin{align*}
    &\tfrac{1}{\alpha} (y_k - x_{k+1})^T\left[\theta(1 - \theta)(v_k - x^*) + \theta^2 (y_k - x^*)\right] \\
    &\qquad = \tfrac{1}{\alpha} (y_k - x_{k+1})^T\left[\theta(1 - \theta)\left(\tfrac{y_k (1 + \theta) - x_k}{\theta} - x^*\right) + \theta^2 (y_k - x^*)\right] \\
    &\qquad = \tfrac{1}{\alpha} (y_k - x_{k+1})^T\left[y_k  - (1 - \theta)x_k - \theta x^*\right],
\end{align*}
where the first equality follows from \eqref{eq:y_k_alternate_strong_convex}.
Substituting these simplified expressions into \eqref{eq:stong_convex_op2_v_k_x^*_equality} completes the proof.
\eproof

Using \cref{lem:v_alternate_strong_convex}, we now establish the Lyapunov function and the recursive relation to analyze \cref{alg:stochastic_proximal_gradient} with \textbf{Option II} under \cref{ass:strong_convexity}.

\begin{lemma} \label{lem:general_descent_acc_strong_convex}
    Suppose \cref{ass:strong_convexity} holds. Then, $\forall k \geq 0$, the iterates generated by \cref{alg:stochastic_proximal_gradient} with \textbf{Option II} and $\{\alpha_k\} = \alpha \leq \tfrac{1}{L}$ satisfy,
    \begin{align*}
        &\phi(x_{k+1}) -\phi^* + \tfrac{\mu}{2}\|v_{k + 1} - x^*\|^2 \\
        &\quad \leq (1 - \theta)\left[\phi(x_k)- \phi^* + \tfrac{\mu}{2}\|v_k - x^*\|^2\right] - \tfrac{\mu}{2}\theta(1 - \theta)\|v_k - y_k\|^2 - \left[\tfrac{1}{2\alpha} - \tfrac{L}{2}\right] \|x_{k+1} - y_k\|^2 \\
        &\qquad + (g_k - \nabla f(y_k))^T (\theta x^* + (1 - \theta) x_{k} - y_k) + (g_k - \nabla f(y_k))^T (y_k - x_{k+1}),
    \end{align*}
    where $\theta = \sqrt{\mu\alpha}$.
\end{lemma}
\bproof
From \eqref{eq:smoothness},
\begin{align*}
    \phi(x_{k+1}) 
    &\leq f(y_k) + \nabla f(y_k)^T (x_{k+1} - y_k) + \tfrac{L}{2} \|x_{k+1} - y_k\|^2 + h(x_{k+1}) \\
    &\leq f(y_k) + \nabla f(y_k)^T (x_{k+1} - y_k) + \tfrac{L}{2} \|x_{k+1} - y_k\|^2 + h(x) - \left(\tfrac{y_{k} - x_{k+1}}{\alpha_k} - g_k\right)^T (x - x_{k+1})\\
    &\leq f(x) - \nabla f(y_k)^T (x - y_k) - \tfrac{\mu}{2} \|x - y_k\|^2 + \nabla f(y_k)^T (x_{k+1} - y_k) + \tfrac{L}{2} \|x_{k+1} - y_k\|^2 \\
    &\quad + h(x) - \left(\tfrac{y_{k} - x_{k+1}}{\alpha_k} - g_k\right)^T (x - x_{k+1})\\
    &= \phi(x)  - \tfrac{\mu}{2} \|x - y_k\|^2  + \tfrac{L}{2} \|x_{k+1} - y_k\|^2 - \left(\tfrac{y_{k} - x_{k+1}}{\alpha_k}\right)^T (x - x_{k+1}) \\
    &\quad - (\nabla f(y_k) - g_k)^T (x - x_{k+1})
\end{align*}
where the second inequality follows, $\forall x\in\Rmbb^d$, from the convexity of $h(x)$ and the definition \eqref{eq:prox_operator_alternate} for $x_{k+1}$ which yields $0 \in g_k +  \partial h(x_{k+1}) + \tfrac{x_{k+1} - y_k}{\alpha_k}$, and the third inequality follows from \eqref{eq:strong_convexity_alternate}. Substituting $x = \theta x^* + (1 - \theta) x_{k}$ where $\theta \in [0, 1]$, from \cref{ass:strong_convexity},
\begin{align*}
    \phi(x_{k+1}) -\phi^*
    &\leq (1 - \theta)\left[\phi(x_k)- \phi^*\right] - \tfrac{\mu}{2}\left[\theta (1 - \theta) \|x_k - x^*\|^2 + \|x - y_k\|^2\right] \numberthis \label{eq:str_conv_lemma_1}\\
    &\quad + \tfrac{L}{2} \|x_{k+1} - y_k\|^2 - \left(\tfrac{y_{k} - x_{k+1}}{\alpha_k}\right)^T (x - x_{k+1}) - (\nabla f(y_k) - g_k)^T (x - x_{k+1}).
\end{align*}
The second term on the right-hand side of \eqref{eq:str_conv_lemma_1} can be simplified as,
\begin{align*}
    &\theta (1 - \theta) \|x_k - x^*\|^2 + \|x - y_k\|^2 \\
    &\quad = \theta (1 - \theta) \|x_k - x^*\|^2 + \|\theta x^* + (1 - \theta) x_{k} - y_k\|^2 \\
    &\quad = \theta (1 - \theta) \|x_k - x^*\|^2 
    + (1 - \theta)^2 \|x_k - x^*\|^2 
    + \|y_k - x^*\|^2 
    + 2(1 - \theta)(x_{k} - x^*)^T (x^* - y_k)  \\
    &\quad = (1 - \theta) \|x_k - x^*\|^2 
    + \|y_k - x^*\|^2 
    + 2(1 - \theta)(x_{k} - x^*)^T (x^* - y_k)  \\
    &\quad \geq (1 - \theta) \|x_k - x^*\|^2 
    + \|y_k - x^*\|^2 
    - (1 - \theta)\|x_k - x^*\|^2 
    - (1 - \theta)\|y_k - x^*\|^2  \\
    &\quad = \theta \|y_k - x^*\|^2,
\end{align*}
where the inequality follows from the identity $2a^Tb \geq -\|a\|^2 -\|b\|^2$, $\forall a,b\in \Rmbb^d$. 
Substituting this bound into \eqref{eq:str_conv_lemma_1} and adding the result from \cref{lem:v_alternate_strong_convex}, we get,
\begin{align*}
    &\phi(x_{k+1}) -\phi^* + \tfrac{\mu}{2}\|v_{k + 1} - x^*\|^2 \\
    &\quad \leq (1 - \theta)\left[\phi(x_k)- \phi^* + \tfrac{\mu}{2}\|v_k - x^*\|^2\right] - \tfrac{\mu}{2}\theta(1 - \theta)\|v_k - y_k\|^2 + \tfrac{L}{2} \|x_{k+1} - y_k\|^2 \\
    &\qquad - \left(\tfrac{y_{k} - x_{k+1}}{\alpha_k}\right)^T (\theta x^* + (1 - \theta) x_{k} - x_{k+1} - y_k + y_k)\\
    &\qquad - (\nabla f(y_k) - g_k)^T (\theta x^* + (1 - \theta) x_{k} - x_{k+1} - y_k + y_k) \\
    &\qquad + \tfrac{1}{2\alpha} \|y_k - x_{k+1}\|^2 
    + \tfrac{1}{\alpha} (y_k - x_{k+1})^T\left[\theta x^* + (1 - \theta)x_k - y_k \right],
\end{align*}
where simplifying the expression completes the proof.
\eproof

We now present the convergence of \cref{alg:stochastic_proximal_gradient} under \cref{ass:strong_convexity} when using a biased gradient estimate, and then discuss the corresponding complexity of number of proximal operator evaluations.

\begin{theorem} \label{th:convergence_strong_convex}
    Suppose \cref{ass:strong_convexity} holds and the gradient estimate $g_k$ satisfies \cref{cond:sampling}.
    \begin{enumerate}
        \item For \cref{alg:stochastic_proximal_gradient} with \textbf{Option I}:
        \begin{enumerate}
            \item For the finite-sum problem \eqref{eq:deter_error_obj}, if the parameters in \cref{cond:sampling} 
            are chosen such that $\{\eta_k\} = \eta \in \left[0, \tfrac{1}{2}\right)$, $\iota_0 < \infty $, and $\delta_k = \delta^k$ $\forall k \geq 0$, where $\delta \in [0, 1)$, and the step size is chosen such that $\{\alpha_k\} = \alpha \leq \tfrac{1 - 4 \eta^2}{2 L}$, then $\{x_k\}$ converges to $x^*$ at a linear rate as, 
            \begin{align*}
                \phi(x_{k}) - \phi^* 
                &\leq \max\left\{1 - \tfrac{\mu\alpha}{3}, \delta^2\right\}^{k+1} \max\left\{\phi(x_0) - \phi^*, \tfrac{12\iota_0^2}{\mu}\right\}.
            \end{align*}

            \item For the expectation problem \eqref{eq:stoch_error_obj}, if the parameters in \cref{cond:sampling} are chosen such that $\{\tilde{\eta}_k\} = \tilde{\eta} \in \left[0, \tfrac{1}{\sqrt{2}}\right)$, $\tilde{\iota}_0 < \infty$, and $\tilde{\delta}_k = \tilde{\delta}^k$ $\forall k \geq 0$, where $\tilde{\delta} \in [0, 1)$, and the step size is chosen such that $\{\alpha_k\} = \alpha \leq \tfrac{1 - 2\tilde{\eta}^2}{2 L}$, then $\{x_k\}$ converges to $x^*$ in expectation at a linear rate as,
            \begin{align*}
                \Embb[\phi(x_{k}) - \phi^*] 
                &\leq \max\left\{1 - \tfrac{\mu\alpha}{3}, \tilde{\delta}^2\right\}^{k+1} \max\left\{\phi(x_0) - \phi^*, \tfrac{6\tilde{\iota}_0^2}{\mu}\right\}.
            \end{align*}
        \end{enumerate}

        \item For \cref{alg:stochastic_proximal_gradient} with \textbf{Option II}:
        \begin{enumerate}
            \item For the finite-sum problem \eqref{eq:deter_error_obj}, if the parameters in \cref{cond:sampling} 
            are chosen such that \{$\eta_k\} = \eta \leq \sqrt{\tfrac{\hat{c}}{4(L + \hat{c})}}$ where $\hat{c} = \tfrac{\mu}{4}(1 - \sqrt{\tfrac{\mu}{L}})$,  $\iota_0 < \infty$, and $\delta_k = \delta^k$ $\forall k \geq 0$, where $\delta \in [0, 1) $, and the step size is chosen as $\{\alpha_k\} = \alpha = \tfrac{1}{2(L + \hat{c})}$, then $\{x_k\}$ converges to $x^*$ at a linear rate as, 
            \begin{align*}
                \phi(x_{k}) -\phi^* \leq 
                \max\left\{1 - \tfrac{\sqrt{\alpha\mu}}{4}, \delta^2\right\}^{k+1}\max\left\{\phi(x_0)- \phi^* + \tfrac{\mu}{2}\|x_0 - x^*\|^2, \tfrac{8\iota_0^2}{\hat{c}\sqrt{\mu\alpha}}\right\}. 
            \end{align*}

            \item For the expectation problem \eqref{eq:stoch_error_obj}, if the parameters in \cref{cond:sampling} are chosen such that \{$\tilde{\eta}_k\} = \tilde{\eta} \leq \sqrt{\tfrac{\hat{c}}{2(L + \hat{c})}}$ where $\hat{c} = \tfrac{\mu}{4}(1 - \sqrt{\tfrac{\mu}{L}})$, $\iota_0 < \infty$, and $\tilde{\delta}_k = \tilde{\delta}^k$ $\forall k \geq 0$, where $\delta \in [0, 1)  $, and the step size is chosen as $\{\alpha_k\} = \alpha = \tfrac{1}{2(L + \hat{c})}$, then $\{x_k\}$ converges to $x^*$ in expectation at a linear rate as, 
            \begin{align*}
                \Embb\left[\phi(x_{k}) - \phi^*\right]
                \leq 
                \max\left\{1 - \tfrac{\sqrt{\alpha\mu}}{4}, \tilde{\delta}^2\right\}^{k+1}\max\left\{\phi(x_0)- \phi^* + \tfrac{\mu}{2}\|x_0 - x^*\|^2, \tfrac{4\tilde{\iota}_0^2}{\hat{c}\sqrt{\mu\alpha}}\right\}.
            \end{align*}
        \end{enumerate}
    \end{enumerate}
\end{theorem}
\bproof
For \cref{alg:stochastic_proximal_gradient} with \textbf{Option I}, consider the result in \cref{lem:prelim_descent_strong_convex} with $y_k = x_k$. Applying young's inequality with constants $c_1, c_2 > 0$ as $(g_k - \nabla f(y_k))^T (x - y_k) \leq \tfrac{\alpha_k}{2 c_1}\|g_k - \nabla f(x_k)\|^2 + \tfrac{ c_1}{2\alpha_k}\|x - x_k\|^2$ and $(\nabla f(y_k) - g_k)^T (x_{k+1} - y_k)\leq \tfrac{\alpha_k}{2 c_2}\|\nabla f(x_k) - g_k\|^2 + \tfrac{c_2}{2 \alpha_k}\|x_{k+1} - x_k\|^2$ yields, 
\begin{align*}
    \phi(x_{k+1}) 
    &\leq 
    \phi(x)  + \left[\tfrac{1}{2\alpha_k} - \tfrac{\mu}{2}\right] \|x - x_k\|^2  + \tfrac{\alpha_k}{2 c_1}\|g_k - \nabla f(x_k)\|^2 + \tfrac{ c_1}{2\alpha_k}\|x - x_k\|^2\\
    &\quad + \tfrac{\alpha_k}{2 c_2}\|\nabla f(x_k) - g_k\|^2 + \tfrac{c_2}{2 \alpha_k}\|x_{k+1} - x_k\|^2 - \left(\tfrac{1}{2 \alpha_k} - \tfrac{L}{2}\right) \|x_{k+1} - x_{k}\|^2 \\
    &= \phi(x)  + \tfrac{1}{2\alpha_k}\left[1 - \mu\alpha_k + c_1\right] \|x - x_k\|^2 + \tfrac{\alpha_k}{2}\left[\tfrac{1}{c_1} + \tfrac{1}{c_2}\right] \|g_k - \nabla f(x_k)\|^2  \\
    &\quad - \tfrac{1}{2\alpha_k} \left(1 - \alpha_k L - c_2\right) \|x_{k+1} - x_{k}\|^2.
\end{align*}
Substituting $x = \nu_k x^* + (1 - \nu_k) x_k$, where $\nu_k \in [0, 1]$, setting $c_1 = c_2 = \tfrac{1}{2}$, and using \cref{ass:strong_convexity} in the above inequality, we get, 
\begin{align*}
    \phi(x_{k+1}) - \phi^*
    &\leq (1 - \nu_k) (\phi(x_k) - \phi^*) + \left[- \tfrac{\mu}{2} \nu_k(1 - \nu_k) + \tfrac{\nu_k^2}{4\alpha_k} \left(3 - 2 \mu \alpha_k \right)\right]\|x_k - x^*\|^2     \\
    &\quad + 2\alpha_k \|g_k - \nabla f(x_k)\|^2 - \tfrac{1}{4\alpha_k} \left(1 - 2 \alpha_k L\right)  \|x_{k+1} - x_{k}\|^2.
\end{align*}
Substituting $\nu_k = \tfrac{2\mu \alpha_k}{3}$ reduces the bound to
\begin{equation} \label{eq:str_cnv_biased_op1_1}
    \phi(x_{k+1}) - \phi^*
    \leq \left(1 -\tfrac{2\mu\alpha_k}{3}\right) (\phi(x_k) - \phi^*) + 2\alpha_k \|g_k - \nabla f(x_k)\|^2 - \tfrac{1}{4\alpha_k} \left(1 - 2 \alpha_k L\right)  \|x_{k+1} - x_{k}\|^2.
\end{equation}

For the finite-sum problem \eqref{eq:deter_error_obj} for \cref{alg:stochastic_proximal_gradient} with \textbf{Option I}, under the defined parameters, substituting \eqref{eq:deter_error_alternate} from \cref{lem:gradient_error_bounds} into \eqref{eq:str_cnv_biased_op1_1} yields,
\begin{align*}
    \phi(x_{k+1}) - \phi^*
    &\leq \left(1 -\tfrac{2\mu\alpha}{3}\right) (\phi(x_k) - \phi^*) + 2\alpha \left(\tfrac{\eta^2}{2\alpha^2} \left\|x_{k+1} - x_k\right\|^2 + 2\iota_0^2 \delta^{2k}\right)\\
    &\quad - \tfrac{1}{4\alpha} \left(1 - 2 \alpha L\right) \|x_{k+1} - x_{k}\|^2 \\
    &= \left(1 -\tfrac{2\mu\alpha}{3}\right) (\phi(x_k) - \phi^*) + 4\alpha \iota_0^2 \delta^{2k} - \tfrac{1}{4 \alpha} \left(1 - 2 \alpha L - 4 \eta^2 \right) \|x_{k+1} - x_{k}\|^2,
\end{align*}
where the constant $1 - 2 \alpha L - 4 \eta^2 \geq 0$
under the defined parameters. Hence, the bound reduces to,
\begin{align*}
    \phi(x_{k+1}) - \phi^*
    &\leq \left(1 -\tfrac{2\mu\alpha}{3}\right) (\phi(x_k) - \phi^*) + 4\alpha \iota_0^2 \delta^{2k},
\end{align*}
where applying \cref{lem:linear_convergence_induction} with $\omega = \tfrac{\mu\alpha}{3}$  completes the proof for the finite-sum problem \eqref{eq:deter_error_obj} for \cref{alg:stochastic_proximal_gradient} with \textbf{Option I}.

For the expectation problem \eqref{eq:stoch_error_obj} for \cref{alg:stochastic_proximal_gradient} with \textbf{Option I}, consider the conditional expectation of \eqref{eq:str_cnv_biased_op1_1} given $\Gcal_k$. From \cref{cond:sampling}, under the defined parameters, we get,
\begin{align*}
    \Embb_k\left[\phi(x_{k+1}) - \phi^*\right]
    &\leq \left(1 -\tfrac{2\mu\alpha}{3}\right) (\phi(x_k) - \phi^*) + 2\alpha \left(\tfrac{\tilde{\eta}^2}{4} \Embb_k\left[\|x_{k+1} - x_{k}\|^2\right] + \tilde{\iota}_0^2\tilde{\delta}^{2k}\right) \\
    &\quad - \tfrac{1}{4\alpha} \left(1 - 2 \alpha L\right) \Embb_k\left[\|x_{k+1} - x_{k}\|^2\right] \\
    &\leq \left(1 -\tfrac{2\mu\alpha}{3}\right) (\phi(x_k) - \phi^*) + 2\alpha\tilde{\iota}_0^2\tilde{\delta}^{2k} - \tfrac{1}{4\alpha}\left(1 - 2\alpha L - 2 \tilde{\eta}^2\right) \Embb_k\left[\|x_{k+1} - x_{k}\|^2\right],
\end{align*}
where the constant $1 - 2\alpha L - 2 \tilde{\eta}^2 \geq 0$ under the defined parameters. Hence, the total expectation of the above bound yields,
\begin{align*}
    \Embb\left[\phi(x_{k+1}) - \phi^*\right]
    &\leq \left(1 -\tfrac{2\mu\alpha}{3}\right) \Embb\left[\phi(x_k) - \phi^*\right] + 2\alpha\tilde{\iota}_0^2\tilde{\delta}^{2k},
\end{align*}
where applying \cref{lem:linear_convergence_induction} with $\omega = \tfrac{\mu\alpha}{3}$ completes the proof for the expectation problem \eqref{eq:stoch_error_obj} for \cref{alg:stochastic_proximal_gradient} with \textbf{Option I}.

For \cref{alg:stochastic_proximal_gradient} with \textbf{Option II}, consider the result in \cref{lem:general_descent_acc_strong_convex}. Substituting $v_k$ from \eqref{eq:y_k_alternate_strong_convex} yields,    
\begin{align*}
    &\phi(x_{k+1}) -\phi^* + \tfrac{\mu}{2}\|v_{k + 1} - x^*\|^2 \numberthis \label{eq:str_conv_biased_opt2_1}\\
    &\quad \leq (1 - \theta)\left[\phi(x_k)- \phi^* + \tfrac{\mu}{2}\|v_k - x^*\|^2\right]  - \tfrac{1}{2\alpha} \left[1 - \alpha L\right] \|x_{k+1} - y_k\|^2 \\
    &\qquad - \tfrac{\mu}{2}\theta(1 - \theta)\left\|\tfrac{(1 + \theta)y_k - x_k}{\theta} - y_k\right\|^2 + (g_k - \nabla f(y_k))^T (y_k - x_{k+1})\\
    &\qquad + (g_k - \nabla f(y_k))^T (\theta x^* + y_k(1 + \theta) - \theta v_k - \theta x_{k} - y_k).
\end{align*}
The last three terms on the right-hand side of \eqref{eq:str_conv_biased_opt2_1} can be simplified using young's inequality with constants $c_1, c_2 > 0$ as,
\begin{align*}
    &- \tfrac{\mu}{2\theta}(1 - \theta)\left\|y_k - x_k\right\|^2 + (g_k - \nabla f(y_k))^T (y_k - x_{k+1}) + \theta (g_k - \nabla f(y_k))^T ( y_k - x_k + x^* - v_k) \\
    &\quad \leq - \tfrac{\mu}{2\theta}(1 - \theta)\left\|y_k - x_k\right\|^2 + \tfrac{\alpha}{2c_1}\|g_k - \nabla f(y_k)\|^2 + \tfrac{c_1}{2\alpha}\|y_k - x_{k+1}\|^2 \\
    &\qquad + \tfrac{\alpha}{2c_2}\left\|g_k - \nabla f(y_k)\right\|^2 + \tfrac{c_2 \theta^2}{2\alpha}\left\| y_k - x_k + x^* - v_k\right\|^2 \\
    &\quad \leq - \left[\tfrac{\mu}{2\theta}(1 - \theta) - \tfrac{c_2 \theta^2}{\alpha}\right]\left\|y_k - x_k\right\|^2 + \tfrac{\alpha}{2}\left[\tfrac{1}{c_1} + \tfrac{1}{c_2}\right]\|g_k - \nabla f(y_k)\|^2 \\
    &\qquad + \tfrac{c_1}{2\alpha}\|y_k - x_{k+1}\|^2 + \tfrac{c_2 \theta^2}{\alpha}\left\| v_k - x^* \right\|^2,
\end{align*}
where the second inequality follows from the identity $\|a+b\|^2 \leq 2(\|a\|^2 + \|b\|^2)$ $\forall a,b \in \Rmbb^d$. Setting $c_2 = \tfrac{\theta(1 - \theta)}{4}$ yields the constant $\tfrac{\mu}{2\theta}(1 - \theta) - \tfrac{c_2 \theta^2}{\alpha} = \tfrac{\mu}{2\theta}(1 - \theta) - \tfrac{\theta (1 - \theta) \mu}{4} = \tfrac{\mu (1 - \theta)}{2}\left[\tfrac{1}{\theta}- \tfrac{\theta}{2}\right] \geq 0$, since $\theta \in (0, 1)$.
Therefore, the corresponding term can be dropped, and \eqref{eq:str_cnv_biased_op1_1} reduces to
\begin{align*}
    &\phi(x_{k+1}) -\phi^* + \tfrac{\mu}{2}\|v_{k + 1} - x^*\|^2 \\
    &\quad \leq (1 - \theta)\left[\phi(x_k)- \phi^*\right] + \left[\tfrac{\mu (1 - \theta)}{2} + \tfrac{c_2 \theta^2}{\alpha} \right]\|v_k - x^*\|^2  \\
    &\qquad - \tfrac{1}{2\alpha}\left[1 - \alpha L - c_1\right] \|x_{k+1} - y_k\|^2 + \tfrac{\alpha}{2}\left[\tfrac{1}{c_1} + \tfrac{1}{c_2}\right]\|g_k - \nabla f(y_k)\|^2 \\
    &\quad \leq \left(1 - \tfrac{\theta}{2}\right)\left[\phi(x_k)- \phi^* + \tfrac{\mu}{2} \|v_k - x^*\|^2 \right] - \tfrac{1}{2\alpha}\left[1 - \alpha L - c_1\right] \|x_{k+1} - y_k\|^2 \numberthis \label{eq:str_cnv_biased_opt2_2}\\
    &\qquad + \tfrac{\alpha}{2}\left[\tfrac{1}{c_1} + \tfrac{1}{c_2}\right]\|g_k - \nabla f(y_k)\|^2,
\end{align*}
where the second inequality follows from $\tfrac{\mu (1 - \theta)}{2} + \tfrac{c_2 \theta^2}{\alpha} = \tfrac{\mu (1 - \theta)}{2} + \tfrac{\theta (1 - \theta) \mu}{4} = \tfrac{\mu}{2}(1 - \theta)\left(1 + \tfrac{\theta }{2} \right) \leq \tfrac{\mu}{2}\left(1 - \tfrac{\theta}{2}\right)$ since $\theta \in (0, 1)$ and $\phi(x_k)- \phi^* \geq 0$.

For the finite-sum problem \eqref{eq:deter_error_obj} for \cref{alg:stochastic_proximal_gradient} with \textbf{Option II}, substituting \eqref{eq:deter_error_alternate} into \eqref{eq:str_cnv_biased_opt2_2} under the defined parameters yields,
\begin{align*}
    &\phi(x_{k+1}) -\phi^* + \tfrac{\mu}{2}\|v_{k + 1} - x^*\|^2 \\
    &\quad \leq \left(1 - \tfrac{\theta}{2}\right)\left[\phi(x_k)- \phi^* + \tfrac{\mu}{2} \|v_k - x^*\|^2 \right] - \tfrac{1}{2\alpha}\left[1 - \alpha L - c_1\right] \|x_{k+1} - y_k\|^2 \\
    &\qquad + \tfrac{\alpha}{2}\left[\tfrac{1}{c_1} + \tfrac{1}{c_2}\right]\left[\tfrac{\eta^2}{2} \left\|\tfrac{x_{k+1} - y_k}{\alpha}\right\|^2 + 2\iota_0^2 \delta^{2k}\right] \\
    &\quad \leq \left(1 - \tfrac{\theta}{2}\right)\left[\phi(x_k)- \phi^* + \tfrac{\mu}{2} \|v_k - x^*\|^2 \right] + \tfrac{\alpha}{2}\left[\tfrac{1}{c_1} + \tfrac{1}{c_2}\right]\left[ 2\iota_0^2 \delta^{2k}\right] \\
    &\qquad - \tfrac{1}{2\alpha}\left[1 - \alpha L - c_1 -\tfrac{\eta^2}{2} \left(\tfrac{1}{c_1} + \tfrac{1}{c_2}\right)\right] \|x_{k+1} - y_k\|^2.
\end{align*}
With $c_1 = \alpha \hat{c}$ and previously defined $c_2 = \tfrac{\sqrt{\alpha \mu}(1 - \sqrt{\alpha \mu})}{4}  \geq \tfrac{\alpha\mu}{4}\left(1 - \sqrt{\tfrac{\mu}{L}}\right) = \alpha \hat{c}$ under the defined parameters (since $\alpha \leq \tfrac{1}{L}$ and $\sqrt{\alpha \mu} \leq 1$), the constant multiplying the last term can be upper bounded as,
\begin{equation*}
    1 - \alpha L - c_1 -\tfrac{\eta^2}{2} \left(\tfrac{1}{c_1} + \tfrac{1}{c_2}\right) \geq 1 - \alpha (L + \hat{c}) - \tfrac{\eta^2 }{\alpha \hat{c}} \geq 1 - \tfrac{(L + \hat{c})}{2(L + \hat{c})} - \tfrac{2 (L + \hat{c}) }{\hat{c}} \tfrac{\hat{c}}{4 (L + \hat{c})} = 0.
\end{equation*}
Thus, the last term in the upper bound can be omitted, yielding,
\begin{align*}
    \phi(x_{k+1}) -\phi^* + \tfrac{\mu}{2}\|v_{k + 1} - x^*\|^2 \leq \left(1 - \tfrac{\theta}{2}\right)\left[\phi(x_k)- \phi^* + \tfrac{\mu}{2} \|v_k - x^*\|^2 \right] + \tfrac{2\iota_0^2 \delta^{2k}}{\hat{c}},
\end{align*}
where applying \cref{lem:linear_convergence_induction} with $\omega = \tfrac{\theta}{4}$ completes the proof for the finite-sum problem \eqref{eq:deter_error_obj} for \cref{alg:stochastic_proximal_gradient} with \textbf{Option II}.

For the expectation problem \eqref{eq:stoch_error_obj} for \cref{alg:stochastic_proximal_gradient} with \textbf{Option II}, consider the conditional expectation of \eqref{eq:str_cnv_biased_opt2_2} given $\Gcal_k$. From \cref{cond:sampling} under the defined parameters, we get,
\begin{align*}
    &\Embb_k\left[\phi(x_{k+1}) -\phi^* + \tfrac{\mu}{2}\|v_{k + 1} - x^*\|^2\right] \\
    &\quad \leq \left(1 - \tfrac{\theta}{2}\right)\left[\phi(x_k)- \phi^* + \tfrac{\mu}{2} \|v_k - x^*\|^2 \right] - \tfrac{1}{2\alpha}\left[1 - \alpha L - c_1\right] \Embb_k \left[\|x_{k+1} - y_k\|^2\right] \\
    &\qquad + \tfrac{\alpha}{2}\left[\tfrac{1}{c_1} + \tfrac{1}{c_2}\right]\left[\tfrac{\tilde{\eta}^2}{4\alpha^2}\Embb_k\left[\|x_{k+1} - y_k\|^2\right] + \tilde{\iota}_0^2\tilde{{\delta}}^{2k}\right] \\
    &\quad \leq \left(1 - \tfrac{\theta}{2}\right)\left[\phi(x_k)- \phi^* + \tfrac{\mu}{2} \|v_k - x^*\|^2 \right] + \tfrac{\alpha}{2}\left[\tfrac{1}{c_1} + \tfrac{1}{c_2}\right]\left[\tilde{\iota}_0^2\tilde{{\delta}}^{2k}\right] \\
    &\qquad - \tfrac{1}{2\alpha}\left[1 - \alpha L - c_1 - \tfrac{\tilde{\eta}^2}{4} \left(\tfrac{1}{c_1} + \tfrac{1}{c_2}\right)\right] \Embb_k \left[\|x_{k+1} - y_k\|^2\right].
\end{align*}
With $c_1 = \alpha \hat{c}$ and previously defined $c_2 = \tfrac{\sqrt{\alpha \mu}(1 - \sqrt{\alpha \mu})}{4}  \geq \tfrac{\alpha\mu}{4}\left(1 - \sqrt{\tfrac{\mu}{L}}\right) = \alpha \hat{c}$ under the defined parameters, the constant for the last term can be upper bounded as,
\begin{equation*}
    1 - \alpha L - c_1 - \tfrac{\tilde{\eta}^2}{4} \left(\tfrac{1}{c_1} + \tfrac{1}{c_2}\right) \geq 1 - \alpha (L + \hat{c}) - \tfrac{\tilde{\eta}^2 }{2\alpha \hat{c}} \geq 1 - \tfrac{(L + \hat{c})}{ 2 (L + \hat{c})} - \tfrac{2 (L + \hat{c}) }{2\hat{c}} \tfrac{\hat{c}}{2 (L + \hat{c})} = 0.
\end{equation*}
Thus, the last term in the upper bound can be omitted, yielding,
\begin{align*}
    &\Embb\left[\phi(x_{k+1}) -\phi^* + \tfrac{\mu}{2}\|v_{k + 1} - x^*\|^2\right] 
    \leq 
    \left(1 - \tfrac{\theta}{2}\right)\Embb\left[\phi(x_k)- \phi^* + \tfrac{\mu}{2}\|v_k - x^*\|^2\right] 
    + \tfrac{\tilde{\iota}_0^2\tilde{{\delta}}^{2k}}{\hat{c}},
\end{align*}
where applying \cref{lem:linear_convergence_induction} with $\omega = \tfrac{\theta}{4}$ completes the proof for the expectation problem \eqref{eq:stoch_error_obj} for \cref{alg:stochastic_proximal_gradient} with \textbf{Option II}.
\eproof

\cref{th:convergence_strong_convex} establishes linear rate of convergence for \cref{alg:stochastic_proximal_gradient} with \textbf{Option I} and \textbf{Option II} for strongly convex objective functions when using a biased gradient estimate satisfying \cref{cond:sampling}.
Thus, an $\epsilon > 0$ accurate solution, with the same definition as in \cref{sec:convexity}, can be achieved in $\Ocal\left(\kappa \log \tfrac{1}{\epsilon}\right)$ iterations (proximal operator evaluations) with \textbf{Option I} and in $\Ocal\left(\sqrt{\kappa} \log \tfrac{1}{\epsilon}\right)$ iterations (proximal operator evaluations) with \textbf{Option II}, where $\kappa = \tfrac{L}{\mu}$ is the condition number, matching the results for deterministic first-order methods \cite{nesterov2018lectures}.
We note that the presented analysis is significantly simpler than that of \cite{schmidt2011convergence}, where the authors analyze accelerated proximal gradient methods with predetermined deterministic errors in the gradient estimate and the solution to the proximal operator. In \cite{schmidt2011convergence}, the analysis establishes that the sequences $\{x_k\}$ and $\{v_k\}$ remain within a finite distance of $x^*$ to accommodate the errors resulting from using a biased gradient estimate. In contrast, the adaptive nature of the error in our gradient estimates allows us to incorporate this error directly into the Lyapunov function, thereby simplifying the analysis.

The parameter settings required in \cref{th:convergence_strong_convex} can be simplified when the gradient estimate for the expectation problem \eqref{eq:stoch_error_obj} is unbiased, as shown in the following corollary.

\begin{corollary} \label{corollary:convergence_strong_convex_unbiased}
    Suppose the conditions in \cref{th:convergence_strong_convex} hold for the expectation problem \eqref{eq:stoch_error_obj} and the gradient estimate is unbiased, i.e. $\Embb_k[g_k] = \nabla f(y_k)$.

    \begin{enumerate}
        \item For \cref{alg:stochastic_proximal_gradient} with \textbf{Option I}, if the parameters in \cref{cond:sampling} are chosen such that $\{\tilde{\eta}_k\} = \tilde{\eta} \in [0, 1)$ and $\tilde{\delta}_k = \tilde{\delta}^k$ $\forall k \geq 0$, where $\tilde{\delta} \in [0, 1)$, and the step size is chosen such that $\{\alpha_k\} = \alpha \leq \tfrac{2 - \tilde{\eta}^2}{2L}$, then $\{x_k\}$ converges to $x^*$ in expectation at a linear rate as,
        \begin{align*}
            \Embb[\phi(x_{k}) - \phi^*]
            \leq 
            \max\left\{1 - \tfrac{\mu\alpha}{2}, \tilde{\delta}^2\right\}^{k+1} \max\left\{\phi(x_0) - \phi^*, \tfrac{2\tilde{\iota}_0^2}{\mu}\right\}.
        \end{align*}
        \item For \cref{alg:stochastic_proximal_gradient} with \textbf{Option II}, if the parameters in \cref{cond:sampling} are chosen such that $\{\tilde{\eta}_k\} = \tilde{\eta} \in [0, 1)$, $\tilde{\iota}_0^2 \in [0, \tilde{\eta})$, $\tilde{\delta}_k = \tilde{\delta}^k$ $\forall k \geq 0$, where $\tilde{\delta} \in [0, 1)$, and the step size is chosen such that $\{\alpha_k\} = \alpha \leq \tfrac{1 - \tilde{\eta} - \tilde{\iota}_0^2}{L}$, then $\{x_k\}$ converges to $x^*$ in expectation at a linear rate as,
        \begin{align*}
            \Embb\left[\phi(x_{k}) -\phi(x^*)\right]
            &\leq 
            \max\left\{1 - \tfrac{\sqrt{\mu\alpha}}{2}, \tilde{\delta}^2\right\}^{k+1} 
            \max\left\{\phi(x_0) - \phi^* + \tfrac{\mu}{2}\|x_0 - x^*\|^2, \sqrt{\tfrac{\alpha}{\mu}}\right\}.
        \end{align*}
    \end{enumerate}
\end{corollary}
\bproof
For \cref{alg:stochastic_proximal_gradient} with \textbf{Option I}, consider the result in \cref{lem:prelim_descent_strong_convex} with $y_k = x_k$. Under the defined parameters and $x = \mu \alpha x^* + (1 - \mu \alpha)x_k$ where $\mu \alpha \leq 1$, using \cref{ass:strong_convexity}, we get,
\begin{align*}
    \phi(x_{k+1}) - \phi^*
    &\leq 
    (1 - \mu\alpha)(\phi(x_k) - \phi^*) 
    - \tfrac{\mu^2 \alpha(1 - \mu\alpha)}{2}\|x_k - x^*\|^2
    + \tfrac{\mu^2\alpha}{2}\left(1 - \mu \alpha\right) \|x_k - x^*\|^2 \\
    &\quad + \mu \alpha (g_k - \nabla f(x_k))^T (x^* - x_k) + (\nabla f(x_k) - g_k)^T (x_{k+1} - x_k) - \tfrac{1}{2\alpha}\left(1 - \alpha L\right) \|x_{k+1} - x_{k}\|^2 \\
    &= 
    (1 - \mu\alpha)(\phi(x_k) - \phi^*) + \mu \alpha (g_k - \nabla f(x_k))^T (x^* - x_k)\\
    &\quad + (\nabla f(x_k) - g_k)^T (x_{k+1} - x_{k}) - \tfrac{1}{2\alpha}\left(1 - \alpha L\right) \|x_{k+1} - x_{k}\|^2.
\end{align*}
Taking conditional expectation of the above inequality given $\Gcal_k$, the second term on the right-hand side is zero due to the assumption of an unbiased gradient estimate, yielding,
\begin{align*}
    \Embb_k\left[\phi(x_{k+1}) - \phi^*\right]
    &\leq 
    (1 - \mu\alpha)(\phi(x_k) - \phi^*) + \Embb_k\left[(\nabla f(x_k) - g_k)^T (x_{k+1} - x_{k})\right] \\
    &\quad - \tfrac{1}{2\alpha}\left(1 - \alpha L\right) \Embb_k\left[\|x_{k+1} - x_{k}\|^2\right] \\
    &\leq 
    (1 - \mu\alpha)(\phi(x_k) - \phi^*) + \tfrac{\tilde{\eta}^2}{4\alpha} \Embb_k\left[\left\|x_{k+1} - x_{k}\right\|^2\right] + \alpha \tilde{\iota}_0^2\tilde{\delta}^{2k}  \\
    &\quad - \tfrac{1}{2\alpha}\left(1 - \alpha L\right) \Embb_k\left[\|x_{k+1} - x_{k}\|^2\right] \\
    &= 
    (1 - \mu\alpha)(\phi(x_k) - \phi^*) + \alpha \tilde{\iota}_0^2\tilde{\delta}^{2k} - \tfrac{1}{2\alpha}\left(1 - \alpha L - \tfrac{\tilde{\eta}^2}{2}\right) \Embb_k\left[\|x_{k+1} - x_{k}\|^2\right],
\end{align*}
where the second inequality follows from Cauchy-Schwartz inequality and \cref{cond:sampling} as,
\begin{align*}
    &\Embb_k[(\nabla f(x_k) - g_k)^T (x_{k+1} - x_{k})] \\
    &\qquad= \Embb_k[(\nabla f(x_k) - g_k)^T (x_{k+1} - x_{k} - \hat{x}_{k+1} + \hat{x}_{k+1})] = \Embb_k[(\nabla f(x_k) - g_k)^T (x_{k+1} - \hat{x}_{k+1})] \\
    &\qquad \leq \Embb_k[\|\nabla f(x_k) - g_k\|\|x_{k+1} - \hat{x}_{k+1}\|] = \alpha\Embb_k[\|\nabla f(x_k) - g_k\|\|R_{\alpha}^{true}(x_k) - R_{\alpha}(x_{k})\|]] \\
    &\qquad \leq \alpha \Embb_k[\|\nabla f(x_k) - g_k\|^2] \leq \tfrac{\tilde{\eta}^2}{4\alpha} \Embb_k\left[\left\|x_{k+1} - x_{k}\right\|^2\right] + \alpha \tilde{\iota}_0^2\tilde{\delta}^{2k}.
\end{align*}
Under the defined parameters, $1 - \alpha L - \tfrac{\tilde{\eta}^2}{2} \geq 0$. Hence, the total expectation of the bound yields,
\begin{align*}
    \Embb\left[\phi(x_{k+1}) - \phi^*\right]
    &\leq 
    (1 - \mu\alpha)\Embb\left[\phi(x_k) - \phi^*\right] + \alpha \tilde{\iota}_0^2\tilde{\delta}^{2k},
\end{align*}
where applying \cref{lem:linear_convergence_induction} with $\omega = \tfrac{\mu\alpha}{2}$ completes the proof for \textbf{Option I}.

For \cref{alg:stochastic_proximal_gradient} with \textbf{Option II}, consider the result in \cref{lem:general_descent_acc_strong_convex} while simplifying the upper bound by ignoring the negative term $- \tfrac{\mu}{2}\theta(1 - \theta)\|v_k - y_k\|^2$. Taking conditional expectation given $\Gcal_k$, from the assumption of an unbiased gradient estimate, the bound reduces to,
\begin{align*}
    &\Embb_k\left[\phi(x_{k+1}) -\phi(x^*) + \tfrac{\mu}{2}\|v_{k + 1} - x^*\|^2\right] \\
    &\quad \leq (1 - \theta)\left[\phi(x_k)- \phi(x^*) + \tfrac{\mu}{2}\|v_k - x^*\|^2\right] - \tfrac{1}{2\alpha}\left(1 - \alpha L\right) \Embb_k\left[\|x_{k+1} - y_k\|^2\right]\\
    &\qquad + \Embb_k\left[(g_k - \nabla f(y_k))^T (y_k - x_{k+1})\right] \\
    &\quad \leq (1 - \theta)\left[\phi(x_k)- \phi(x^*) + \tfrac{\mu}{2}\|v_k - x^*\|^2\right] - \tfrac{1}{2\alpha}\left(1 - \alpha L\right) \Embb_k\left[\|x_{k+1} - y_k\|^2\right]\\
    &\qquad + \tfrac{\alpha (\tilde{\eta} + \tilde{\iota}_0^2)}{2} \Embb_k\left[\|R_{\alpha}(y_k)\|^2\right] + \tfrac{\alpha \tilde{\delta}^{2k}}{2} \\
    &\quad= (1 - \theta)\left[\phi(x_k)- \phi(x^*) + \tfrac{\mu}{2}\|v_k - x^*\|^2\right] + \tfrac{\alpha \tilde{\delta}^{2k}}{2} - \tfrac{\alpha}{2} (1 - \alpha L - \tilde{\eta} - \tilde{\iota}_0^2) \Embb_k\left[\|R_{\alpha}(y_k)\|^2\right],
\end{align*}
where the second inequality follows from \eqref{eq:stoch_error_prod_squared} from \cref{lem:gradient_error_bounds}. 
Under the defined parameters, $1 - \alpha L - \tilde{\eta}_k - \tilde{\iota}_0^2 \geq 0$. Hence, the total expectation of the bound yields,
\begin{align*}
    \Embb\left[\phi(x_{k+1}) -\phi(x^*) + 
    \tfrac{\mu}{2}\|v_{k + 1} - x^*\|^2\right]
    &\leq (1 - \theta)\Embb\left[\phi(x_k)- \phi(x^*) + \tfrac{\mu}{2}\|v_k - x^*\|^2\right] + \tfrac{\alpha \tilde{\delta}^{2k}}{2},
\end{align*}
where applying \cref{lem:linear_convergence_induction} with $\omega = \tfrac{\theta}{2}$ completes the proof for \textbf{Option II}.
\eproof

The parameter settings in \cref{corollary:convergence_strong_convex_unbiased} are less restrictive compared to those in \cref{th:convergence_strong_convex}, due to the unbiased nature of the gradient approximation, particularly with respect to the step size in \textbf{Option II}.

We now present the complexity for the number of stochastic gradient evaluations for \cref{alg:stochastic_proximal_gradient} when an unbiased gradient estimate is used for the expectation problem \eqref{eq:stoch_error_obj} under \cref{ass:strong_convexity}. Unlike in \cref{sec:non_convexity} and \cref{sec:convexity}, the parameter setting $\{\eta_k\} = 0$ is incorporated directly into the next theorem, as the optimal complexity for stochastic gradient evaluations is achieved under this parameter setting, along with a definition of an $\epsilon$-accurate solution similar to that in \cite{ghadimi2012optimal}.

\begin{theorem} \label{th:sample_complexity_strongly_convex}
    Suppose Assumptions \ref{ass:bounded_variance} and \ref{ass:strong_convexity} hold, and \cref{cond:sampling} is satisfied for the expectation problem \eqref{eq:stoch_error_obj} via the unbiased gradient estimate in \cref{lem:unbiased_sample_set_bound}. Then, to achieve a solution satisfying 
    $\min\left\{\Embb\left[\phi(x_k) - \phi^*\right], \left\|R_{\alpha_k}^{true}(y_k)\right\|^2\right\} \leq \epsilon$
    with $\epsilon > 0$, the number of stochastic gradient evaluations required is as follows:

    \begin{enumerate}
        \item For \cref{alg:stochastic_proximal_gradient} with \textbf{Option I}, if the parameters in \cref{cond:sampling} are chosen such that $\{\tilde{\eta}_k\} = \tilde{\eta} \in [0, 1)$  and $\tilde{\delta}_k  = \tilde{\delta}^k$ $\forall k \geq 0$, where $\tilde{\delta}^2 = 1 - \tfrac{\mu\alpha}{2} $, $\iota_0 > 0$, and the step size is chosen such that $\{\alpha_k\} = \alpha \leq \tfrac{2 - \tilde{\eta}^2}{2L}$, then 
        the number of stochastic gradient evaluations required is $\Ocal \left(\tfrac{\kappa}{\epsilon} \log \left(\tfrac{1}{\epsilon}\right)\right)$, where $\kappa = \tfrac{L}{\mu}$. If $\{\tilde{\eta}_k\} = 0$, then the number of stochastic gradient evaluations reduces to $\Ocal \left(\tfrac{\kappa}{\epsilon} \right)$ to achieve a solution satisfying $\Embb\left[\phi(x_k) - \phi^*\right] \leq \epsilon$ with $\epsilon > 0$. 

        \item For \cref{alg:stochastic_proximal_gradient} with \textbf{Option II}, if the parameters in \cref{cond:sampling} are chosen such that $\{\tilde{\eta}_k\} = \tilde{\eta} \in [0, 1)$, $\tilde{\iota}_0^2 \in (0, \tilde{\eta})$, $\tilde{\delta}_k = \tilde{\delta}^k$ $\forall k \geq 0$, where $\tilde{\delta}^2 = 1 - \tfrac{\sqrt{\mu\alpha}}{2} $, and the step size is chosen such that $\{\alpha_k\} = \alpha \leq \tfrac{1 - \tilde{\eta} - \tilde{\iota}_0^2}{L}$, then the number of stochastic gradient evaluations required is $\Ocal \left(\tfrac{\sqrt{\kappa}}{\epsilon} \log \left(\tfrac{1}{\epsilon}\right)\right)$, where $\kappa = \tfrac{L}{\mu}$. If $\{\tilde{\eta}_k\} = 0$, then the number of stochastic gradient evaluations reduces to $\Ocal\left(\tfrac{\sqrt{\kappa}}{\epsilon}\right)$ to achieve a solution satisfying $\Embb\left[\phi(x_k) - \phi^*\right] \leq \epsilon$ with $\epsilon > 0$.  
    \end{enumerate}
\end{theorem}

\bproof
Let $K_\epsilon \geq 1$ be the first iteration that achieves the desired solution accuracy. Hence, $\left\|R_{\alpha_k}^{true}(y_k)\right\|^2 > \epsilon$ $\forall k \leq K_\epsilon - 1$ and from \eqref{eq:complexity_bound_help}, 
$\left\|\Embb_k\left[R_{\alpha_k}(y_k)\right]\right\|^2 > \left(1 + \tfrac{\tilde{\eta}^2}{4}\right)^{-1}\left[\tfrac{\epsilon}{2} - \tilde{\iota}_0^2\tilde{\delta}_k^2\right]$ $\forall k \leq K_\epsilon - 1$.
Thus, the total number of stochastic gradient evaluations can be bounded using \cref{lem:unbiased_sample_set_bound} as, 
\begin{align*}
    \sum_{k=0}^{K_\epsilon - 1}|S_k| 
    &= \sum_{k=0}^{K_\epsilon - 1} \left\lceil\tfrac{\sigma^2}{\tfrac{\tilde{\eta}_k^2}{4} \left\|\Embb_k\left[R_{\alpha_k}(y_k)\right]\right\|^2 + \tilde{\iota}_0^2\tilde{\delta}_k^2} \right\rceil \leq \sum_{k=0}^{K_\epsilon - 1} \tfrac{2\sigma^2 (4 + \tilde{\eta}^2)}{\tilde{\eta}^2 \epsilon + 8\tilde{\iota}_0^2\tilde{\delta}_k^2}  + 1 
    \leq \sum_{k=0}^{K_\epsilon - 1} \tfrac{2\sigma^2(4 + \tilde{\eta}^2)}{\tilde{\eta}^2 \epsilon} + \tfrac{\sigma^2 (4 + \tilde{\eta}^2)}{4\tilde{\iota}_0^2\tilde{\delta}_k^2}  + K_{\epsilon}.
\end{align*}

For \cref{alg:stochastic_proximal_gradient} with \textbf{Option I}, $K_\epsilon$ is at most $\Ocal\left(\log_{\tfrac{1}{\tilde{\delta}^2}}\left(\tfrac{1}{\epsilon}\right)\right)$ from \cref{corollary:convergence_strong_convex_unbiased}, yielding, 
\begin{align*}
    \sum_{k=0}^{K_\epsilon-1}|S_k| 
    \leq \tfrac{2\sigma^2 (4 + \tilde{\eta}^2)}{\tilde{\eta}^2 \epsilon}K_\epsilon + \tfrac{\sigma^2(4 + \tilde{\eta}^2)}{4\tilde{\iota}_0^2} \tfrac{\tfrac{1}{\tilde{\delta}^{2(K_\epsilon + 1)}} - 1}{\tfrac{1}{\tilde{\delta}^{2}} - 1} + K_{\epsilon} = \Ocal \left(\tfrac{\kappa}{\epsilon} \log \left(\tfrac{1}{\epsilon}\right)\right).
\end{align*}
Following the same procedure, if $\{\tilde{\eta}_k\} = 0$,
a solution satisfying $\Embb\left[\phi(x_k) - \phi^*\right] \leq \epsilon$ with $\epsilon > 0$ is achieved in $\tilde{K}_\epsilon = \Ocal\left(\log_{\tfrac{1}{\tilde{\delta}^2}}\left(\tfrac{1}{\epsilon}\right)\right)$ iterations from \cref{corollary:convergence_strong_convex_unbiased}, and the total number of stochastic gradient evaluations is $\sum_{k=0}^{\tilde{K}_\epsilon-1}|S_k| \leq \tfrac{\sigma^2}{\tilde{\iota}_0^2} \tfrac{\tfrac{1}{\tilde{\delta}^{2(\tilde{K}_\epsilon + 1)}} - 1}{\tfrac{1}{\tilde{\delta}^{2}} - 1} + \tilde{K}_{\epsilon} = \Ocal \left(\tfrac{\kappa}{\epsilon} \right)$.

For \cref{alg:stochastic_proximal_gradient} with \textbf{Option II}, $K_\epsilon$ is at most $\Ocal\left(\log_{\tfrac{1}{\tilde{\delta}^2}}\left(\tfrac{1}{\epsilon}\right)\right)$ from \cref{corollary:convergence_strong_convex_unbiased}, yielding, 
\begin{align*}
    \sum_{k=0}^{K_\epsilon-1}|S_k| 
    \leq \tfrac{2\sigma^2(4 + \tilde{\eta}^2)}{\tilde{\eta}^2 \epsilon}K_\epsilon + \tfrac{\sigma^2(4 + \tilde{\eta}^2)}{4\tilde{\iota}_0^2} \tfrac{\tfrac{1}{\tilde{\delta}^{2(K_\epsilon + 1)}} - 1}{\tfrac{1}{\tilde{\delta}^{2}} - 1} + K_{\epsilon} = \Ocal \left(\tfrac{\sqrt{\kappa}}{\epsilon} \log \left(\tfrac{1}{\epsilon}\right)\right).
\end{align*}
Following the same procedure, if $\{\tilde{\eta}_k\} = 0$,
a solution satisfying $\Embb\left[\phi(x_k) - \phi^*\right] \leq \epsilon$ with $\epsilon > 0$ is achieved in $\tilde{K}_\epsilon = \Ocal\left(\log_{\tfrac{1}{\tilde{\delta}^2}}\left(\tfrac{1}{\epsilon}\right)\right)$ iterations from \cref{corollary:convergence_strong_convex_unbiased}, and the total number of stochastic gradient evaluations is $\sum_{k=0}^{\tilde{K}_\epsilon-1}|S_k| \leq \tfrac{\sigma^2}{\tilde{\iota}_0^2} \tfrac{\tfrac{1}{\tilde{\delta}^{2(\tilde{K}_\epsilon + 1)}} - 1}{\tfrac{1}{\tilde{\delta}^{2}} - 1} + \tilde{K}_{\epsilon} = \Ocal \left(\tfrac{\sqrt{\kappa}}{\epsilon} \right)$.
\eproof
\cref{th:sample_complexity_strongly_convex} matches the optimal complexity for the number of stochastic gradient evaluations for the expectation problem \eqref{eq:stoch_error_obj} with a strongly convex objective function \cite{ghadimi2012optimal}. Unlike the cases of nonconvex and general convex objective functions, a predefined sequence of 
decreasing gradient errors results in the optimal complexity for strongly convex objective functions. 
\section{Numerical Experiments} \label{sec:exp}
In this section, we illustrate the performance of \cref{alg:stochastic_proximal_gradient} 
on a synthetic strongly convex quadratic problem of the form
\begin{equation} \label{eq:quad_problem}
    \min_{x \in \Rmbb^{10}} \; \phi(x) = \frac{1}{N}  \sum_{i=1}^N \frac{1}{2} x^T Q_i x + b_i^T x  + \Ical_{\left[\|x\|^2 \leq 1\right]},
\end{equation}
where $Q_i \in \Rmbb^{10 \times 10}$ 
is positive definite ($Q_i \succ 0$) and $b_i \in \Rmbb^{10}$ are samples from a finite dataset with $N = 10^5$, generated via the process outlined in \cite{mokhtari2016network} with condition number $\kappa \approx 10^4$. 
The problem is constrained to the convex feasible region described by $\|x\|^2 \leq 1$, using the indicator function of a set as
\begin{equation*}
    \Ical_{\left[\|x\|^2 \leq 1\right]} = 
    \begin{cases}
    0 & \text{if } \|x\|^2 \leq 1, \\
    \infty & \text{otherwise.}
\end{cases}
\end{equation*}
Thus, the objective function is strongly convex\footnote{ Additional numerical experiments for $l_1$-regularized binary classification logistic regression problems are provided in \cref{app:logistic_Exps}.}.

The gradient estimate $g_k$ used in \cref{alg:stochastic_proximal_gradient} is a sample average approximation with sample set $S_k \subseteq \Scal$ $\forall k \geq 0$. We present results for five different sample selection strategies ($S_k$). The ``Deterministic" label corresponds to using the true problem gradient, i.e., the full dataset. The ``Stochastic" label corresponds to using 256 samples every iteration, i.e., $|S_k| = 256$ $\forall k \geq 0$.
The ``Geometric" label corresponds to starting from 32 samples and increasing the number of samples by 5\% each iteration, i.e., $|S_{k+1}| = \left\lceil1.05|S_k|\right\rceil$ $\forall k \geq 0$, equivalent to setting $\{\tilde{\eta}_k\} = 0$ in \cref{cond:sampling}.
The ``Adaptive" label corresponds to using \cref{cond:sampling} to control the accuracy of the gradient estimate, with $\{\tilde{\eta}_k\} = \tilde{\eta} = 0.1$, $\iota_0 = 0$, and the sample size selected using sampled estimates as described in \cite{xie2024constrained}. The ``Stochastic", ``Geometric" and ``Adaptive" strategies approximate problem \eqref{eq:quad_problem} as an expectation problem over a uniform distribution, as done in \cite{xie2024constrained,pham2020proxsarah,j2016proximal}, and use an unbiased gradient estimate by sampling each iteration with replacement independent of the current iterate. The ``Adaptive-biased" label corresponds to using \cref{cond:sampling} to control the accuracy of the gradient estimate by determining the sample size with the same parameters, while introducing bias by maintaining $S_{k} \subseteq S_{k+1}$ $\forall k \geq 0$.
We evaluate both \textbf{Option I} (Proximal Gradient) and \textbf{Option II} (Accelerated Proximal Gradient) for each sample selection strategy, with results labeled using ``-I" and ``-II", respectively. For \textbf{Option II}, we set ${\beta_k} = \beta = \tfrac{\sqrt{\kappa} - 1}{\sqrt{\kappa} + 1} \approx 0.98$.
The step size $\{\alpha_k\} = \alpha$ was tuned for each result over the set $\{10^{-i} | i = 0, 1, \ldots, 6  \}$.
The performance was measured by the optimality gap in function value, evaluated against both the number of proximal operator evaluations and the number of stochastic gradient evaluations.

When \cref{cond:sampling} is used to control the accuracy of the gradient estimate, the method achieves one of the most efficient performances in terms of gradient evaluations, while also outperforming constant sample size strategies in terms of proximal operator evaluations. Finally, although introducing bias in the gradient estimate degrades performance, the algorithm still converges and performs well under \textbf{Option II}.

\begin{figure}[H]
    \centering
    \includegraphics[width=\linewidth]{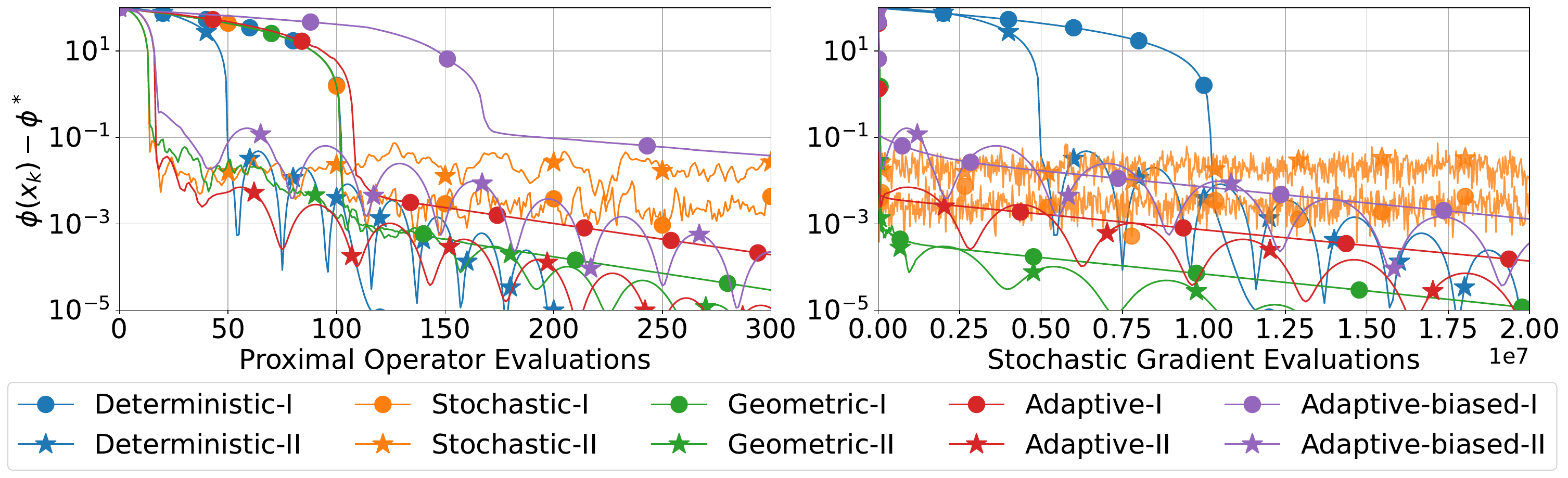}
    \caption{Optimality Gap ($\phi(x_k) - \phi^*$) with respect to the number of proximal operator evaluations and the number of stochastic gradient evaluations of \cref{alg:stochastic_proximal_gradient}, evaluated under \textbf{Option I} (Proximal Gradient) and \textbf{Option II} (Accelerated Proximal Gradient), using ``Deterministic", ``Stochastic", ``Geometric", ``Adaptive" and ``Adaptive-biased" sampling strategies on the strongly convex quadratic problem \eqref{eq:quad_problem}.}
    \label{fig:quad_exp}
\end{figure}

The results are summarized in \cref{fig:quad_exp}. First, for all sample selection strategies, \textbf{Option II}  yields better performance with respect to number of proximal operator evaluations during the initial phase. Second, a stochastic gradient estimate is initially more efficient than the deterministic approach with respect to the number of gradient evaluations, but less efficient in terms of 
number of proximal operator evaluations. When \cref{cond:sampling} is used to control the accuracy of the gradient estimate, the method achieves the most efficient performance in terms of gradient evaluations, while also outperforming constant sample size strategies in terms of proximal operator evaluations.
Finally, although introducing bias in the gradient estimate 
deteriorates performance, 
the algorithm still converges and 
performs well under \textbf{Option II}.

\section{Final Remarks} \label{sec:conc}

In this paper, we have proposed a proximal gradient method and an accelerated proximal gradient method for composite optimization problems, where the smooth component is either a finite-sum function or an expectation of a stochastic function.
The methods employed possibly biased estimates for the gradient of the smooth component, with the accuracy of these estimates adaptively adjusted via extensions of generalized ``norm" conditions tailored to the composite optimization setting to achieve computational efficiency.
For nonconvex, convex, and strongly convex objective functions, the methods achieved the optimal iteration complexity even with 
biased gradient estimates.
When the gradient estimate is unbiased, we refined the analysis which allowed for less restrictive parameter settings. In this case, the methods simultaneously achieved the optimal complexity for both the number of proximal operator evaluations and the number of stochastic gradient evaluations for nonconvex, convex, and strongly convex objective functions. Finally, we conducted preliminary numerical experiments that validated our theoretical results.

\bibliographystyle{plain}
\bibliography{references.bib}

\appendix

\section{Technical Results}
In this section, we present some technical results that have been used in the paper.

\begin{lemma} \label{lem:a_b_Appendix}
    Given $a_1, a_2 \in \Rmbb^d$ and $c \in \Rmbb$,
    \begin{equation*}
        a_1^T a_2 + (c - 1) \|a_1\|^2 = \tfrac{1}{2}(\|a_2\|^2 - \|a_1 - a_2\|^2 + (2c - 1) \|a_1\|^2).
    \end{equation*}
\end{lemma}
\bproof 
The proof follows as,
    \begin{align*}
        a_1^T a_2 + (c - 1) \|a_1\|^2 = \tfrac{1}{2}(2a_1^T a_2 - \|a_1\|^2 + (2c - 1) \|a_1\|^2) = \tfrac{1}{2}(\|a_2\|^2 - \|a_1 - a_2\|^2 + (2c - 1) \|a_1\|^2).
    \end{align*}
\eproof

\begin{lemma} \label{lem:linear_convergence_induction}
    Given a non-negative sequence $\{T_k\}$ such that $T_{k+1} \leq \rho_1 T_k + a \rho_2^k$ where $\rho_1, \rho_2 \in [0, 1)$ and $ 0 \leq a < \infty$, the sequence $\{T_k\} \rightarrow 0$ at a linear rate as,
    \begin{equation*}
        T_{k} \leq \max\{\rho_1 + \omega, \rho_2\}^{k+1} \max\left\{T_0, \tfrac{a}{\omega}\right\},
    \end{equation*}
    where $\omega > 0$ such that $\rho_1 + \omega < 1$.
\end{lemma}
\bproof
The proof follows by induction. For $T_0$, the result trivially holds. Then, if the result holds for $T_k$,
\begin{align*}
    T_{k+1} \leq \rho_1 T_k + a \rho_2^k &\leq \rho_1 \max\{\rho_1 + \omega, \rho_2\}^k \max\left\{T_0, \tfrac{a}{\omega}\right\} + a \rho_2^k \\
    &\leq \max\{\rho_1 + \omega, \rho_2\}^k \max\left\{T_0, \tfrac{a}{\omega}\right\} \left(\rho_1  + \omega\right) \\
    &\leq \max\{\rho_1 + \omega, \rho_2\}^{k+1} \max\left\{T_0, \tfrac{a}{\omega}\right\},
\end{align*}
thus completing the proof.
\eproof

\begin{lemma} \label{lem:bounded_sequence_appendix_nonacc}
    Given non-negative sequences $\{T_k\}$, $\{a_k\}$ and $\{s_k\}$ such that $T_{k+1} \leq (1 + a_k)T_k + s_k$, $T_0 < \infty$, $\sum_{k=0}^{\infty}a_k < \infty$, 
    $\sum_{k=0}^{\infty} s_k < \infty$, 
    then the sequence $\{T_k\}$ is bounded. 
\end{lemma}
\bproof
Unrolling the recursion, 
\begin{align*}
T_{k+1} \leq T_0\prod_{i=0}^{k} (1 + a_i) + \sum_{i=0}^{k} s_i \prod_{j=i+1}^{k}  (1 + a_i) \leq T_0\prod_{i=0}^{\infty} (1 + a_i) + \sum_{i=0}^{\infty} s_i \prod_{j=0}^{\infty}  (1 + a_j), 
\end{align*}
where the second inequality holds due to all the additive terms being non-negative and the product terms being at least one. From \cite[Theorem 1]{trench1999conditional}, since $\sum_{k=0}^{\infty}a_k < \infty$ with $\{a_k\} \geq 0$, the infinite product $\prod_{i=0}^{\infty} (1 + a_i) < \infty$. Thus, $T_{k+1} \leq \prod_{i=0}^{\infty} (1 + a_i) \left[T_0 + \sum_{i=0}^{\infty} s_i \right] < \infty$.
\eproof

\begin{lemma} \label{lem:bounded_sequence_appendix_acc}
    Given non-negative sequences $\{R_k\}$, $\{T_k\}$, $\{a_k\}$ and $\{s_k\}$ such that $R_{k+1} + T_{k+1} \leq R_k + (1 + a_k)T_k + s_k$, $T_0 < \infty$, $R_0 < \infty$, $\sum_{k=0}^{\infty} s_k < \infty$, and $a_k = \hat{\rho}\rho_k$ $\forall k \geq 0$ where $\sum_{k=0}^{\infty}\rho_k < \infty$, $\sum_{k=0}^{\infty} k \rho_k < \infty$ and $\hat{\rho}$ can be controlled to be sufficiently small, then the sequence $\{R_k\}$ is bounded. 
\end{lemma}
\bproof
Let $\sum_{k=0}^{\infty} s_k < \bar{s}$, $\sum_{k=0}^{\infty} a_k < \bar{a}$.
From \cite[Theorem 1]{trench1999conditional}, as $\sum_{k=0}^{\infty}a_k < \infty$ with $\{a_k\} \geq 0$, $\exists \hat{a} > 0$ such that $\prod_{i=0}^{\infty} (1 + a_i) < \hat{a}$.
We now unroll the recursion for $T_k$ as,
\begin{align*}
    T_{k+1} 
    &\leq (1 + a_k)T_k  + s_k + R_k - R_{k+1}
    \leq T_0 \prod_{i=0}^{k} (1 + a_i)  + \sum_{i=0}^{k} (s_i + R_i - R_{i+1}) \prod_{j=i+1}^{k}  (1 + a_j) \\
    &\leq T_0 \prod_{i=0}^{k} (1 + a_i)  + \sum_{i=0}^{k} (s_i + R_i) \prod_{j=i+1}^{k}  (1 + a_j)
    \leq \prod_{i=0}^{\infty} (1 + a_i)  \left[T_0 + \sum_{i=0}^{k} (s_i + R_i)\right] \\
    &\leq \hat{a} \left[T_0 + \sum_{i=0}^{k} (s_i + R_i)\right]. \numberthis \label{eq:appndx_Tk} 
\end{align*}
We unroll the recursion for $R_k$ using a telescopic sum as,
\begin{align*}
R_{k+1} - R_0 &\leq  T_0 - T_{k+1} + \sum_{i=0}^{k} \left(a_iT_i +  s_i\right).
\end{align*}
Further unrolling the above inequality and using \eqref{eq:appndx_Tk} yields,
\begin{align*}
    R_{k+1} &\leq R_0 + T_0 + \sum_{i=0}^{k} \left(a_iT_i +  s_i\right) \\
    &=   R_0 + T_0 + \hat{a}T_0\sum_{i=0}^{k} a_i + \hat{a} \sum_{i=0}^{k} \left(a_i \sum_{j=0}^{i-1} s_j\right)+ \hat{a} \sum_{i=0}^{k} \left(a_i \sum_{j=0}^{i-1}R_j\right)+ \sum_{i=0}^{k} s_i \\
    &\leq  R_0 + T_0 + \hat{a} T_0 \sum_{i=0}^{k} a_i + \hat{a} \left(\sum_{i=0}^{k} s_i\right) \left(\sum_{i=0}^{k} a_i\right) + \hat{a} \sum_{i=0}^{k} \left(a_i \sum_{j=0}^{i-1}R_j\right)+ \sum_{i=0}^{k} s_i \\
    &\leq  R_0 + T_0 + \hat{a} \bar{a} T_0 + \hat{a} \bar{a}\bar{s} + \hat{a} \sum_{i=0}^{k} \left(a_i \sum_{j=0}^{i-1}R_j\right) + \bar{s} = \hat{R} + \hat{a} \sum_{i=0}^{k} \left(a_i \sum_{j=0}^{i-1}R_j\right),
\end{align*}
where $\hat{R} = R_0 + T_0 + \hat{a} \bar{a} T_0 + \hat{a} \bar{a}\bar{s} + \bar{s}$. Let $\sum_{k=0}^{\infty} k \rho_k = \bar{\rho}$ and we define a constant $C = \tfrac{\hat{R}}{1-\hat{a}\hat{\rho}\bar{\rho}} > 0$ for sufficiently small $\hat{\rho}$ such that $1-\hat{a}\hat{\rho}\bar{\rho} \in (0, 1)$. We show via induction that $\{R_k\} \leq C$. First, $R_0 \leq \hat{R} \leq C$ as $1-\hat{a}\hat{\rho}\bar{\rho} \in (0, 1)$. Then, if the induction holds for $k \geq 0$,
\begin{align*}
    R_{k+1} \leq \hat{R} + \hat{a} \sum_{i=0}^{k} \left(a_i \sum_{j=0}^{i-1}R_j\right) \leq \hat{R} + \hat{a} \sum_{i=0}^{k} \left(a_i i C\right) \leq \hat{R} + \hat{a} \hat{\rho}\bar{\rho}C
    = \hat{R} + \hat{a} \hat{\rho}\bar{\rho} \tfrac{\hat{R}}{1-\hat{a}\hat{\rho}\bar{\rho}} = C,
\end{align*} 
completing the proof.
\eproof

\section{Additional Proofs}
In this section, we present proofs that have been omitted from the paper for brevity.

\begin{lemma} \label{lem:condtion_rearrange}
    Suppose \cref{ass:base_smoothness} holds and \cref{cond:sampling} is satisfied in \cref{alg:stochastic_proximal_gradient}. Then,
    \begin{enumerate}
        \item For the finite-sum problem \eqref{eq:deter_error_obj}: 
        \begin{align*}
            \left(1 - \tfrac{\eta_k}{2}\right)\|g_k - \nabla f(y_{k})\| \leq \tfrac{\eta_k}{2} \left\|R_{\alpha_k}^{true}(y_k)\right\| + \iota_0 \delta_k, \quad \forall k \geq 0.            
        \end{align*}
        \item For the expectation problem \eqref{eq:stoch_error_obj}: 
        \begin{align*}
            \left(1 - \tfrac{\tilde{\eta}^2_k}{2}\right)\Embb_k\left[\|g_k - \nabla f(y_k)\|^2\right] \leq \tfrac{\tilde{\eta}_k^2}{2} \Embb_k\left[\left\|R^{true}_{\alpha_k}(y_k)\right\|^2\right] + \tilde{\iota}_0^2\tilde{\delta}_k^2, \quad \forall k \geq 0.
        \end{align*}    
    \end{enumerate}
\end{lemma}
\bproof 
For the finite-sum problem \eqref{eq:deter_error_obj}, using \cref{cond:sampling}, we get,
\begin{align*}
    \|g_k - \nabla f(y_{k})\| 
    &\leq \tfrac{\eta_k}{2} \left\|R_{\alpha_k}(y_k) - R_{\alpha_k}^{true}(y_k) + R_{\alpha_k}^{true}(y_k)\right\| + \iota_0 \delta_k \\
    &\leq \tfrac{\eta_k}{2} \left\|R_{\alpha_k}(y_k) - R_{\alpha_k}^{true}(y_k)\right\| + \tfrac{\eta_k}{2} \left\|R_{\alpha_k}^{true}(y_k)\right\| + \iota_0 \delta_k \\
    &\leq \tfrac{\eta_k}{2} \left\|g_k - \nabla f(y_k)\right\| + \tfrac{\eta_k}{2} \left\|R_{\alpha_k}^{true}(y_k)\right\| + \iota_0 \delta_k,
\end{align*}
where the final inequality follows form \eqref{eq:projected_grad_error}.
Rearranging the final inequality yields the desired result for the finite-sum problem \eqref{eq:deter_error_obj}.

For the expectation problem \eqref{eq:stoch_error_obj}, from \cref{cond:sampling} and using Jensen's inequality, we get, 
\begin{align*}
    \Embb_k\left[\|g_k - \nabla f(y_k)\|^2\right] 
    &\leq \tfrac{\tilde{\eta}_k^2}{4} \Embb_k\left[\left\|R_{\alpha_k}(y_k) - R^{true}_{\alpha_k}(y_k) + R^{true}_{\alpha_k}(y_k)\right\|^2\right] + \tilde{\iota}_0^2\tilde{\delta}_k^2 \\
    &\leq \tfrac{\tilde{\eta}_k^2}{2} \Embb_k\left[\left\|R_{\alpha_k}(y_k) - R^{true}_{\alpha_k}(y_k)\right\|^2\right] + \tfrac{\tilde{\eta}_k^2}{2} \Embb_k\left[\left\|R^{true}_{\alpha_k}(y_k)\right\|^2\right] + \tilde{\iota}_0^2\tilde{\delta}_k^2 \\
    &\leq \tfrac{\tilde{\eta}_k^2}{2} \Embb_k\left[\|g_k - \nabla f(y_k)\|^2\right]  + \tfrac{\tilde{\eta}_k^2}{2} \Embb_k\left[\left\|R^{true}_{\alpha_k}(y_k)\right\|^2\right] + \tilde{\iota}_0^2\tilde{\delta}_k^2,
\end{align*}
where the second inequality follows from the identity $(a + b)^2 \leq 2a^2 + 2b^2$ and the final inequality follows from \eqref{eq:projected_grad_error}.
Rearranging the final inequality yields the desired result for the expectation problem \eqref{eq:stoch_error_obj}.
\eproof

\begin{lemma} \label{lem:deter_sample_set_bound}    
    Suppose Assumptions \ref{ass:base_smoothness} and \ref{ass:bounded_variance} hold in \cref{alg:stochastic_proximal_gradient} for the finite-sum problem \eqref{eq:deter_error_obj}. Let $g_k = \tfrac{1}{|S_k|}\sum_{\xi \in S_k} \nabla F(y_k, \xi)$, where $S_k \subseteq \Scal$ $\forall k \geq 0$. Then, \cref{cond:sampling} is satisfied $\forall k \geq 0$ if
    \begin{equation*}
        |S_k| = \left\lceil \tfrac{N}{\left(1 + \tfrac{\eta_k \left\|R_{\alpha_k}(y_k)\right\| + 2\iota_0 \delta_k}{2\sigma}\right)}\right\rceil.
    \end{equation*}
\end{lemma}
\bproof
In iteration $k \geq 0$, from the definition of $g_k$, the gradient error can be bounded as,
\begin{align*}
    \left\|g_k - \nabla f(y_k)\right\| 
    &= \tfrac{1}{|S_k|}\left\| \sum_{\xi \in S_k} \left(\nabla F(y_k, \xi) - \nabla f(y_k)\right)\right\| \\
    &= \tfrac{1}{|S_k|}\left\| \sum_{\xi \in \Scal / S_k} \left(\nabla F(y_k, \xi) - \nabla f(y_k)\right)\right\| \\
    &\leq \tfrac{1}{|S_k|}\sum_{\xi \in \Scal / S_k}\left\|\nabla F(y_k, \xi) - \nabla f(y_k)\right\| \\
    &\leq \tfrac{1}{|S_k|}\sum_{\xi \in \Scal / S_k} \sigma = \tfrac{N - |S_k|}{|S_k|} \sigma,
\end{align*}
where the second equality following from the definition of the finite-sum problem \eqref{eq:deter_error_obj} and the last inequality follows from \cref{ass:bounded_variance}.
From the defined sample size, we have $ |S_k| \geq \tfrac{N}{\left(1 + \tfrac{\eta_k \left\|R_{\alpha_k}(y_k)\right\| + 2\iota_0 \delta_k}{2\sigma}\right)}$. Therefore, the gradient error can be bounded as,
\begin{align*}
    \left\|g_k - \nabla f(y_k)\right\| 
    &\leq \tfrac{N}{|S_k|} \sigma - \sigma 
    \leq \left(1 + \tfrac{\eta_k \left\|R_{\alpha_k}(y_k)\right\| + 2\iota_0 \delta_k}{2\sigma}\right) \sigma - \sigma 
    = \tfrac{\eta_k \left\|R_{\alpha_k}(y_k)\right\| + 2\iota_0 \delta_k}{2},
\end{align*}
satisfying \cref{cond:sampling}.
\eproof

\section{Additional Numerical Experiments} \label{app:logistic_Exps}

\begin{figure}[H]
\centering
\begin{subfigure}{\textwidth}
    \includegraphics[width=\textwidth]{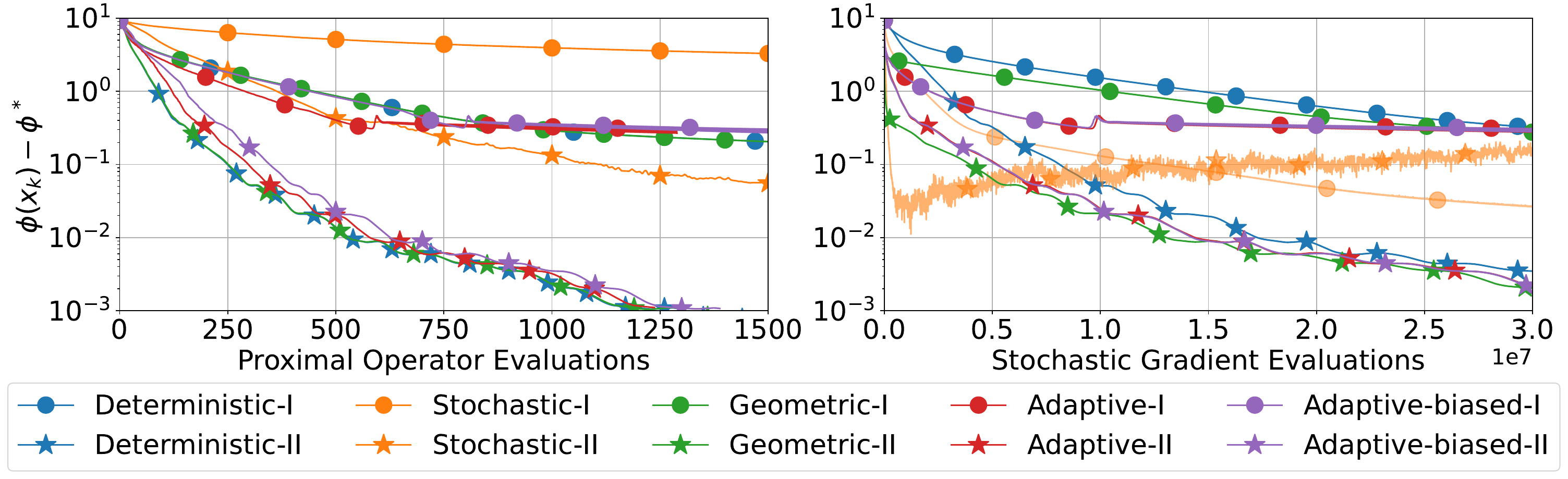} 
    \caption{a9a dataset ($d = 123$, $N = 32,561$, \cite{CC01a})}
    \label{fig:a9a}
\end{subfigure}
\begin{subfigure}{\textwidth}
    \includegraphics[width=\textwidth]{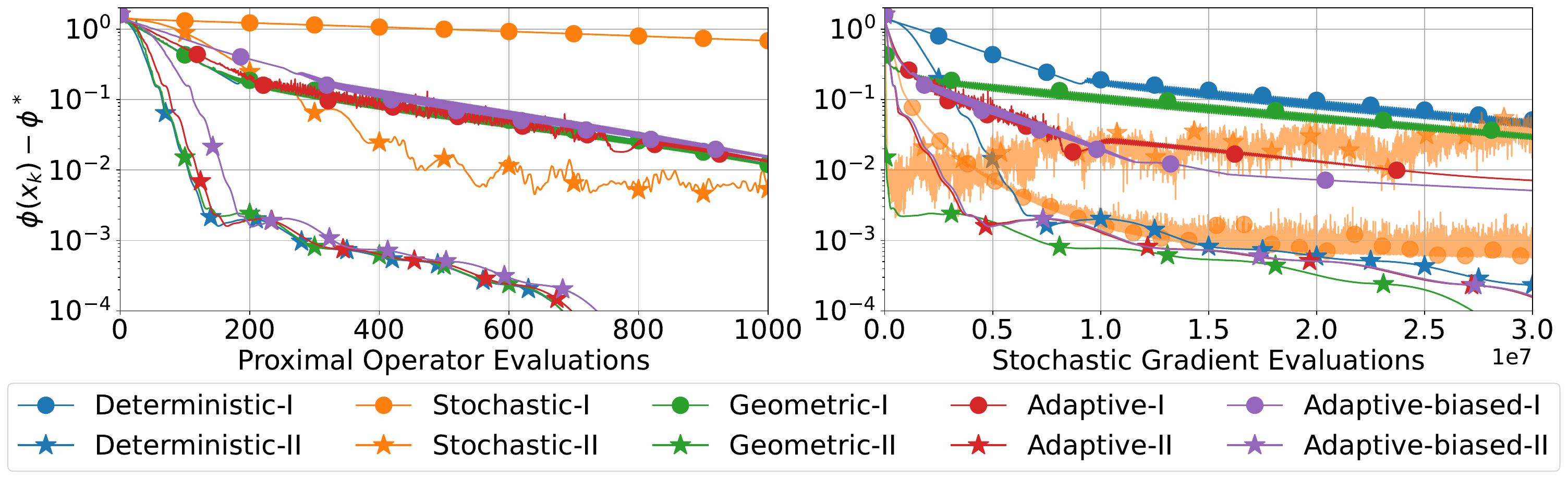} 
    \caption{ijcnn dataset ($d = 23$, $N = 49,990$, \cite{CC01a})}
    \label{fig:ijcnn}
\end{subfigure}
\caption{Optimality Gap ($\phi(x_k) - \phi^*$) with respect to the number of proximal operator evaluations and the number of stochastic gradient evaluations of \cref{alg:stochastic_proximal_gradient}, evaluated under \textbf{Option I} (Proximal Gradient) and \textbf{Option II} (Accelerated Proximal Gradient), using ``Deterministic", ``Stochastic", ``Geometric", ``Adaptive" and ``Adaptive-biased" sampling strategies for $l_1$-regularized binary classification logistic regression \eqref{eq:logistic_problem} on; (a) a9a dataset ($d = 123$, $N = 32,561$, \cite{CC01a}) and (b) ijcnn dataset ($d = 23$, $N = 49,990$).}
\label{fig:logistic_exps}
\end{figure}

In this section, we illustrate the performance of \cref{alg:stochastic_proximal_gradient} on $l_1$-regularized binary classification logistic regression problems of the form
\begin{equation} \label{eq:logistic_problem}
    \min_{x \in \Rmbb^d} \phi(x) = \frac{1}{N} \sum_{i=1}^N \log(1 + e^{-b_iA_ix}) + \frac{1}{N}\|x\|_1, 
\end{equation}
where $a_i \in \Rmbb^d$ is the feature vector (including one for the bias term) and $b_i \in \{0, 1\}$ is the label for each datapoint  $i \in \{1, 2, \dots, N\}$. Experiments were performed 
on the a9a dataset ($d = 123$, $N = 32,561$) and the ijcnn dataset ($d = 23$, $N = 49,990$)\cite{CC01a}.
The sequence $\{\beta_k\}$ for \textbf{Option II} is chosen as described in \cref{sec:convexity}, since the logistic regression binary cross-entropy loss is convex not strongly convex. All other implementation details are the same as in \cref{sec:exp}. The results are summarized in \cref{fig:logistic_exps}.

\end{document}